\documentclass[11pt,leqno]{amsart}
  
\usepackage{latexsym,enumitem}
\usepackage{amssymb}
\usepackage[cp850]{inputenc}
\usepackage{epsfig}
\usepackage{psfrag}
\usepackage{amsthm}
\usepackage{amscd}
\usepackage{amsmath}
\usepackage{amsfonts}
\usepackage{graphics,caption}
\usepackage[all]{xy}
\usepackage{etoolbox}
\patchcmd{\quote}{\rightmargin}{\leftmargin 2em \rightmargin}{}{}
\usepackage[sort,numbers]{natbib} 
\usepackage{mathtools}

\usepackage[english]{babel}
\usepackage{comment}
\usepackage{color}
\usepackage{hyperref}

\let\phi\varphi

\DeclareMathOperator{\eqdef}{\doteq\,} 
\newcommand{\f}[2]{\frac{#1}{#2}} 

\let\epsilon\varepsilon
\let\subset\subseteq
\newcommand{\be}{\begin{equation*}}
 \newcommand{\ee}{\end{equation*}}
\newcommand{\bpf}{\begin{dimo}}
\newcommand{\epf}{\end{dimo}}
\newcommand{\bdefi}{\begin{defin}}
\newcommand{\edefi}{\end{defin}}
\newcommand{\bthm}{\begin{thm}}
\newcommand{\ethm}{\end{thm}}
\newcommand{\blem}{\begin{lem}}
\newcommand{\elem}{\end{lem}}
\newcommand{\bcor}{\begin{cor}}
\newcommand{\ecor}{\end{cor}}

\newcommand{\bprop}{\begin{prop}}
\newcommand{\eprop}{\end{prop}}
\newcommand{\bese}{\begin{ese}}
\newcommand{\eese}{\end{ese}}
\newcommand{\brem}{\begin{rem}}
\newcommand{\erem}{\end{rem}}
\newcommand{\bpfc}{\begin{dimoclaim}}
\newcommand{\epfc}{\end{dimoclaim}}
\newcommand{\eps}{\epsilon}
\newcommand{\al}{\alpha}

\newcommand{\rar}{\rightarrow} 
\newcommand{\imm}{\looparrowright}

\newcommand{\abs}[1]{\left\lvert#1\right\rvert}						
\newcommand{\set}[1]{\left\{#1\right\}}					
	
\newcommand{\quotient}[2]{\left.\raisebox{.1em}{$#1\!$}\middle/\raisebox{-.1em}{$#2$}\right.}

\DeclareMathOperator{\emp}{\varnothing} 
\DeclareMathOperator{\N}{\mathbb N}			
\DeclareMathOperator{\R}{\mathbb R}			
\DeclareMathOperator{\C}{\mathbb C}			
\DeclareMathOperator{\sing}{\mathrm{sing}}		

\newcommand{\diffeo}{\xrightarrow{{}_{\simeq}}} 
\newcommand{\bH}{\mathbb H^3}
 \newcommand{\hyp}[1]{\quotient{\bH}{#1}}
  \newcommand{\iso}{\overset{\text{iso}}\simeq}

  \DeclareMathOperator{\ind}{ind}

\DeclareMathOperator{\Diff}{Diff}    

\newenvironment{quot}
{
	\vspace{-0.2cm}
	\vspace{0.2cm}
}

\theoremstyle{definition}
\newtheorem{d1}{Definition}[section] 

\newenvironment{defin}
{
	\begin{quot}
		\begin{d1}
		}
		{\end{d1}
	\end{quot}

}

\theoremstyle{definition}
\newtheorem{r1}[d1]{Remark}

\newenvironment{rem}
{
	\begin{quot}
		\begin{r1}
		}
		{\end{r1}
	\end{quot}
}

\theoremstyle{definition}
\newtheorem{e1}[d1]{Exercise}

\theoremstyle{definition}
\newtheorem{ese1}[d1]{Example}

\newenvironment{ese}
{
	\begin{quot}
		\begin{ese1}
	}
	{	
		\end{ese1}
	\end{quot}
}

\theoremstyle{definition}

\theoremstyle{definition}
\newtheorem{f2}[d1]{Fact}

\theoremstyle{definition}

\theoremstyle{definition}

\theoremstyle{definition}
\newtheorem{t1}[d1]{Theorem}

\newenvironment{thm}
{
	\begin{quot}
		\begin{t1}}
		{\end{t1}
	\end{quot}
}

\theoremstyle{definition}
\newtheorem*{T1*}{Theorem}

\newenvironment{teor*}
{
	\begin{quot}
		\begin{T1*}}
		{\end{T1*}
	\end{quot}
}

\newenvironment{dimo}
{\begin{proof}[Proof]
	}
	{\end{proof}}

\newenvironment{dimoclaim}{\emph{Proof of Claim:}\;}{\hfill$\square$}

	\theoremstyle{definition}
	\newtheorem{l1}[d1]{Lemma}
	
	\newenvironment{lem}
	{
		\begin{quot}
			\begin{l1}}
			{\end{l1}
		\end{quot}
	}
	\theoremstyle{definition}
	\newtheorem{p1}[d1]{Proposition}
	
	\newenvironment{prop}
	{
		\begin{quot}
			\begin{p1}}
			{\end{p1}
		\end{quot}
	}
	
	\theoremstyle{definition}
	\newtheorem{c1}[d1]{Corollary}
	
	\newenvironment{cor}
	{
		\begin{quot}
			\begin{c1}}
			{\end{c1}
		\end{quot}
	}

		\renewenvironment{abstract}
	{\list{}{\rightmargin\leftmargin}%
		\item[\textbf{Abstract:}]\relax}
	{\endlist}

\newtheorem{innercustomthm}{Theorem}
\newenvironment{customthm}[1]
  {\renewcommand\theinnercustomthm{#1}\innercustomthm}
  {\endinnercustomthm}
  
  \newtheorem{innercustomdef}{Definition}
\newenvironment{customdef}[1]
  {\renewcommand\theinnercustomdef{#1}\innercustomdef}
  {\endinnercustomthm}

  \newtheorem{innercustomprop}{Proposition}
\newenvironment{customprop}[1]
  {\renewcommand\theinnercustomprop{#1}\innercustomprop}
  {\endinnercustomthm}

 \newtheorem*{Theorem*}{Theorem}
 \newtheorem*{Proposition*}{Proposition}
 \newtheorem*{Lemma*}{Lemma}

\usepackage{lipsum} 
\usepackage{etoolbox}
\usepackage{ulem}
\usepackage{natbib}
\usepackage{overpic}
\usepackage{setspace}
\usepackage{extarrows}

\setlist{  
  listparindent=\parindent,
  parsep=\parskip,
}

\usepackage[margin=1.25in,  marginparwidth=1in]{geometry}
\makeatletter
\def\subsection{\@startsection{subsection}{2}%
  \z@{.5\linespacing\@plus.7\linespacing}{.5\linespacing}%
  {\normalfont\bfseries}}
\makeatother

\makeatletter
\def\subsubsection{\@startsection{subsubsection}{2}%
  \z@{.5\linespacing\@plus.7\linespacing}{.5\linespacing}%
  {\normalfont\bfseries}}
\makeatother

\makeatletter
\def\paragraph{\@startsection{paragraph}{4}
  \z@\z@{-\fontdimen2\font}
  {\normalfont\bfseries}}
\makeatother

\patchcmd{\subsection}{-.5em}{.5em}{}{}

\newtheorem*{theorem*}{Theorem}

\renewcommand{\tilde}{\widetilde}

\renewcommand{\bar}{\overline}

\renewcommand{\P}{\mathcal P}
\renewcommand{\S}{\mathbb S^1}
\renewcommand{\hat}{\widehat}

\newcommand{\rev}[1]{{\color{black}#1}}

\renewcommand{\eqdef}{:=}
	
\DeclareTextFontCommand{\emph}{\itshape}

\usepackage{xcolor} 
 \usepackage{tcolorbox}
\newtcbox{\hl}[1][red]{on line, arc=7pt,colback=#1!10!white,colframe=#1!50!black,
  before upper={\rule[-3pt]{0pt}{10pt}},boxrule=1pt, boxsep=0pt,left=6pt,
  right=6pt,top=2pt,bottom=2pt}
 
\begin{document}
	 
\title[On volumes and filling collections of multicurves]{On volumes and filling collections of multicurves}
\author[T. Cremaschi, J. A. Rodriguez-Migueles and Andrew Yarmola]{TOMMASO CREMASCHI, JOSE ANDRES RODRIGUEZ-MIGUELES AND ANDREW YARMOLA}

\maketitle
\begin{abstract} Let $S$ be a surface of negative Euler characteristic and consider a finite filling collection  $\Gamma$  of closed curves on $S$ in minimal position. An observation of Foulon and Hasselblatt shows that $PT(S) \setminus \hat{\Gamma}$ is a finite-volume hyperbolic 3-manifold, where $PT(S)$ is the projectivized tangent bundle and $\hat\Gamma$ is the set of tangent lines to $\Gamma$. In particular, $vol(PT(S) \setminus \hat{\Gamma})$ is a mapping class group invariant of the collection $\Gamma$. \rev{ When $\Gamma$ is a filling pair of simple closed curves, we show that  this volume is coarsely comparable to Weil-Petersson distance between strata in Teichm\"uller space}. Our main tool is the study of \rev{{\it stratified} hyperbolic links $\bar\Gamma$ in a Seifert-fibered space $N$ over $S$. For such links, the volume of $N\setminus\bar\Gamma$ is coarsely comparable to expressions involving distances in the pants graph.}
\end{abstract}

 \section{Introduction}
  
Let $S = \Sigma_{g,k}$ be a hyperbolic surface, that is an orientable smooth surface of genus $g$ with $k$ punctures and negative Euler characteristic. Associated to $S$ is the $3$-manifold $PT(S)$ --- the projectivized tangent bundle. For any finite collection  $\Gamma$ of  smooth essential closed curves on $S$, there is a \emph{canonical} lift $\hat \Gamma$ in $PT(S)$ realized by the set of tangent lines to $\Gamma$. Drilling $\hat \Gamma$ from $PT(S)$ produces a $3$-manifold $N_{\hat\Gamma} = PT(S) \setminus \hat \Gamma$ and any invariant of $N_{\hat\Gamma}$ naturally becomes a mapping class group invariant of $\Gamma$. When $\Gamma$ is filling and in minimal position, Foulon and Hasselblatt \cite{FH13} observed that $N_{\hat\Gamma}$ admits a complete hyperbolic metric of finite-volume, in particular $vol(N_{\hat\Gamma})$ is such an invariant.  \rev{ Recall that $\Gamma$ is filling if $S\setminus \Gamma$ is a collection of disks and  once-puncture disks.} One should think of $\Gamma$ as a weak version of a link diagram where ``over'' and ``under'' crossings are encoded by the tangent directions to $\Gamma$. \rev{An alternate perspective is that these canonical links are examples of Legendrian links for the natural contact structure on $PT(S)$.} The main result of this paper shows that if $\Gamma$ is a filling pair of multicurves, then $vol(N_{\hat\Gamma})$ is  \rev{ coarsely comparable} to the Weil-Petersson distance between strata of nodal surfaces corresponding to pinching different components of $\Gamma$. \rev{ Recall that a \emph{multicurve} is a disjoint union of non-parallel essential \emph{simple} closed curves.} 

Several upper and lower bounds for $vol(N_{\bar\Gamma})$ in terms of invariants of $\Gamma$ have been studied in recent literature, see \cite{BPS171,BPS172, M17}. For example, in \cite{CM19} the first and second author gave an upper bound which is linear in terms of the self-intersection number of $\Gamma$, a fact reminiscent of classical results in knot theory. In fact, they show that there exists $\Gamma$ for which this bound is asymptotically optimal for arbitrary lifts.  In \cite{BPS171}, it is shown that for every hyperbolic structure $X$ on $S$, there is a constant $C_X$ such that $vol(N_{\hat\Gamma}) \leq C_X \ell_X(\Gamma)$, where $\ell_X$ gives the length of the geodesic representative.  However, since $vol(N_{\hat\Gamma})$ is independent of the choice of $X$, this bound demonstrates some odd behaviour, see \cite{M17} for examples.

The best known lower bound appears in \cite{M17,CM19, M20}, where the bound is given in terms of the number of homotopy classes of arcs of $\Gamma$ after cutting $S$ open along any multicurve $\mathfrak m$\, and taking the maximum over such $\mathfrak m$. While this lower bound is shown to be sharp for some families of {\it non-simple} closed curves on the modular surface \cite{M20}, it is always at most $6(3g + n)(3g -3 + n)$ whenever $\Gamma$ is composed entirely of {\it simple} closed curves. Our main results helps address this ineffectiveness.

{\bf Bounds on filling pairs.} Let $\P(S)$ and $\mathcal C(S)$ denote the pants graph and curve graph of $S$, respectively. In addition, let $\mathcal T(S)$ be the Teichm\"uller space of $S$, $d_{WP}$ the Weil-Petersson metric on $\mathcal T(S)$, and let \rev{$\mathcal T_\gamma \subset \partial \mathcal T(S)$ denote the strata of noded surfaces where we pinch a multicurve $\gamma$. Here,  $\partial \mathcal T(S)$ is the boundary of the Weil-Petersson completion of $T(S)$.}

\rev{Given a set $P$ and two functions $f, g: P \to \mathbb R$, we say that $f$ and $g$ are \emph{coarsely comparable} if there are constants  $A \geq 1$ and $B \geq 0$ such that $g(x)/A - B \leq f(x) \leq A g(x) + B$ for all $x \in P$. With this notation in mind, we can now state our main result.}

\begin{customthm}{A}\label{pair} Let $S$ be a hyperbolic surface, $N = PT(S)$, and let $(\al, \beta)$ be a filling pair of essential \rev{multicurves} on $S$ in minimal position. Then:
\rev{\[ vol(N_{(\hat{\al}, \hat{\beta})}) \asymp d_{WP}(\mathcal T_\alpha , \mathcal T_\beta)\]
where $\asymp$ denotes coarsely comparability with constants depending only on $S$.}
\end{customthm}

\rev{In the proof, we actually demonstrate that  \[vol(N_{(\hat{\al}, \hat{\beta})}) \asymp \inf_{P_{\al},P_{\beta}}
 d_{\mathcal P(S)}(P_{\al},P_{\beta})\] where $P_{\al}, P_{\beta}$ are any pants decompositions of $S$ with $\al \subset P_{\al}, \beta \subset P_{\beta}$.  The connection to Weil-Petersson geometry is derived from the work of Brock \cite{Bro01}, where it is established that $\P(S)$ and $(\mathcal T(S),d_{WP})$ are quasi-isometric metric spaces.}

When $S = \Sigma_{1,1}$ or $\Sigma_{0,4}$, we obtain a more general result involving the curve graph $\mathcal C(S)$.

\begin{customthm}{B}\label{punctorus} Let $S = \Sigma_{1,1}$ or $\Sigma_{0,4}$, $N = PT(S)$, and let $\Gamma$ be a filling collection of non-parallel essential simple closed curves in minimal position. Then, there exists an ordering $\Gamma = \{\gamma_i\}_{i = 1}^n$ and universal constants $K_1 \geq 1$, $K_0 \geq 0$, such that:
$$  \frac{1}{2K_1}\left(\sum_{i =1}^n d_{\mathcal C(S)}(\gamma_i,\gamma_{i+1})\right) \leq vol(N_{\widehat{\Gamma}})\leq K_1 \left(\sum_{i = 1}^n d_{\mathcal C(S)}(\gamma_i,\gamma_{i+1})\right)+nK_0,$$
where $\gamma_{n+1} = \gamma_1$. \rev{ When $n= 2$, this a special case of Theorem \ref{pair}.}
\end{customthm}

\rev{To prove these results, we first study} the special case where $\al$ and $\beta$ are leaves of singular line fields that are obtained by \rev{rotations from} each other. In this setting, the line fields can be realized as disjoint surfaces in $PT(S)$, which allows us to use 3-manifold topology and geometry. We adapt Canary's Interpolation Theorem \cite{Ca1996} and techniques of Brock \cite{Bro01} to our setting by replacing the use of simplicial hyperbolic surfaces with a mild-generalization of collapsed simplicial ruled surfaces (CSRS) developed in \cite[Section 5]{BS2017}. Finally, we show \rev{that this special case is coarsely comparable} to the general setting using topological arguments. It is important to note that unlike the lower bound in  \cite{M17,CM19}, our bounds control the volume of the thick-part of $N_{\widehat{\Gamma}}$.

Previous results of this flavor have related pants graph distance to volumes of convex cores or related translation lengths of mapping classes to volumes of mapping tori, see \cite{Bro01, BM, field2021endperiodic}. \rev{Since our results deal with finite-volume hyperbolic 3-manifolds, they may automatically benefit from Dehn surgery techniques and their implications for volumes of hyperbolic 3-manifolds}. In particular, this observation is reminiscent of an announced result of Brock and \rev{Minsky}, referenced in \cite[Theorem 1.4]{BM}, which states that the Weil-Petersson length spectrum of moduli space is well-ordered and has order type $\omega^\omega$.


{\bf Length and self-intersection bounds.} We can relate our topological bounds back to bounds in terms of length and intersection number by constructing surfaces with specific geometries  via the following theorem.

\begin{customthm}{C}\label{ce}
Let $S  = \Sigma_{g,k}$ be a hyperbolic surface with $(g,k)$ equal to $(0,m+4)$, $(1,m+1)$ or $(g,2+gm)$ for $m\in\N$ and \rev{let} $N = PT(S)$. Then, there exists a sequence $\{(\alpha_n,\beta_n)\}_{n\in\mathbb N}$ of filling pairs of essential simple closed curves on $S$ in minimal position with the property that $d_{\mathcal C(S)}(\alpha_n,\beta_n)\nearrow \infty$  and a hyperbolic metric $X$ on $S$ such that:
$$vol(N_{(\hat{\alpha}_n,\, \hat\beta_n)})\asymp d_{\mathcal C(S)}(\alpha_n,\beta_n)\asymp \log\left(\ell_{X}(\alpha_n)+\ell_{X}(\beta_n)\right)$$
where $\asymp$ denotes \rev{coarsely comparability} with constant depending only on $S$.\end{customthm}

When  $S = \Sigma_{1,1}$, $\Sigma_{0,4}$, this construction gives rise to asymptotically optimal upper bounds.

\begin{customthm}{D}\label{length}Let $S = \Sigma_{1,1}$ or $\Sigma_{0,4}$, $N = PT(S)$, and let $(\al_n, \beta_n)$ be a sequence of filling pairs of essential simple closed curves in minimal position on $S$. Then there is a constant $C > 0$ and for any hyperbolic metric $X$ on $S$ there is a constant $C_X > 0$ such that:
\[  \limsup_{n \to \infty} \frac{vol(N_{(\hat{\al_n}, \hat{\beta_n})})}{\log(\iota(\al_n, \beta_n))} \leq C \quad \text{and} \quad \limsup_{n \to \infty} \frac{vol(N_{(\hat{\al_n}, \hat{\beta_n})})}{\log(\ell_{X}(\al_n) + \ell_{X}(\beta_n))} \leq C_X\]
Further, there are sequences where both equalities are attained.
\end{customthm}

\brem In reference to Theorem \ref{length}, we strongly suspect that $\log(\ell_X)$ is asymptotically \rev{sharp} in more generality for \emph{stratified canonical lifts}, which we define in the next section. For canonical lifts of \emph{non-simple} closed curves, there are examples where volume grows at least as fast as $\ell_X/\log(\ell_X)$, see \cite{M17}. Note, volumes of canonical lifts can grow at most as fast as $\ell_X$ by \cite{BPS171}. For arbitrary lifts, the story is different. In \cite{M17}, there are examples where volume grows at least as fast as $\ell_X$. Further, since $i(\Gamma, \Gamma) \leq K_X \ell_X(\Gamma)^2$ for any collection of essential closed curves by \cite{Basmajian2013}, the intersection number bound of \cite{CM19} gives a universal upper bound of $\ell_X^2$ for arbitrary lifts.
\erem

\subsection{Stratified links in Seifert-fibered spaces}

Many of the bounds mentioned in the previous section apply in more generality. We find that it is useful to consider these generalizations even though the tangent bundle setting appears to be the most interesting due to its canonical nature, connections to periodic orbits of the geodesic flow, \rev{Legendrian links}, and quadratic differentials. Our generalized setting is as follows.

Let $\overline{\Gamma}$ be a link in a Seifert-fibered space $N$ over $S$. In \cite{CM19}, the authors demonstrate that  $N_{\bar\Gamma}\eqdef  N\setminus\bar\Gamma$ is hyperbolic whenever $N_{\bar\Gamma}$ is acylindrical and $\bar\Gamma$ is a lift of a filling collection $\Gamma$ of essential closed curves in minimal position via the canonical Seifert projection $\pi:N\rar S$. Note, the acylindricity condition is needed only if one allows multiple copies of the same loop in $\Gamma$. We are interested in the relationship between $vol(N_{\bar\Gamma})$ and topological invariants of $\Gamma$. As before, we think of $\Gamma$ as a weak version of a link diagram, but now the crossing data is entirely forgotten. To have some control over how these crossings behave, we introduce the notion of \emph{stratification} for $\bar\Gamma$ as follows.

\bdefi\label{1}  A link $\bar\Gamma \subset N$ is \emph{stratified} if one can decompose $\overline{\Gamma} $ into sub-links $ \set{\bar\Gamma_j}_{i = 1}^n$ such that there is a disjoint properly embedded collection of incompressible surfaces $\{S_i\}_{i = 1}^n \subset N$ with $\bar\Gamma_i$ embedded in $S_i$ for each $i$. Otherwise, $\bar\Gamma$ is said to be \emph{unstratifiable}.

The $S_i$ are called \emph{stratification surfaces} and we assume that $\{S_i\}_{i = 1}^n$  is \emph{cyclically ordered around the $\S$-fiber}. We call the collection of pairs $\mathcal H =  \{(\bar\Gamma_i, S_i)\}_{i = 1}^n$ a \emph{stratification} of $\bar\Gamma$.
\edefi

\begin{figure}[htb]
\begin{overpic}[scale=0.9] {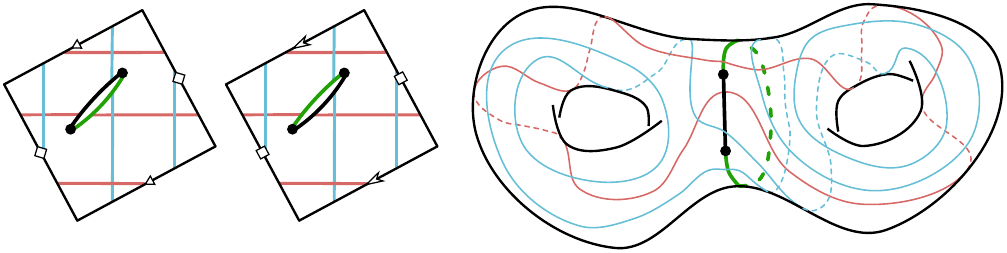}
\put(8.8,9){$\beta$}
\put(13.3,11.7){$\al$}
\end{overpic}
\caption{Filling curves $\al, \beta$ on $\Sigma_{2,2}$ with $\hat\al, \hat\beta$ stratified in $PT(\Sigma_{2,2})$.}
\label{slit-torus}
\end{figure}

\bese A motivating example is the slit torus construction in Figure \ref{slit-torus}. There, the horizontal curve $\al$ and the vertical curve $\beta$ fill $\Sigma_{2,0}$. If we puncture the vertices of the slit, $\hat\al$ and $\hat\beta$ become stratified in $PT(\Sigma_{2,2})$ by the horizontal and vertical line fields arising from the flat structure. In Section \ref{examples}, we extend this construction to a general setting.

\eese

To demonstrate that this definition arises naturally in several settings, we give abundant natural families of stratified canonical lifts of filling collections of multicurves in Section \ref{examples}.

\rev{In Section \ref{strat}, we show that a stratification $\mathcal H =  \{(\bar\Gamma_i, S_i)\}_{i = 1}^n$ of a hyperbolic link can be realized by \emph{horizontal} surfaces. Recall, a \emph{horizontal} surface in a Seifert-fibered space is one that is everywhere transverse to the $\S$-fiber. This allows us to conclude that $S_i \iso_f S_1$, were $\iso_f$ denotes fiber-wise isotopy.} We make a canonical choice for these homeomorphisms by cutting $N$ along $S_{\mathcal H} \eqdef  S_1$ and \rev{use} the ``downward'' \rev{fiber-wise} projection to define $\Gamma_i^{\mathcal H} \iso_f  \bar\Gamma_i$. Let $\psi_{\mathcal H}$ denote the ``bottom-to-top'' monodromy of the cut along $S_{\mathcal H}$, then \rev{ the Seifert-fibered space $N$ over $S$ } is naturally homeomorphic to the mapping torus $S_{\mathcal H} \times I / (x,0) \sim (\psi_{\mathcal{H}}(x), 1)$.  Lastly, define the \emph{complexity} of a stratification as $\kappa(\mathcal H) = 3g_{\mathcal H} - 3 + k_{\mathcal H}$, where $S_{\mathcal H} \cong \Sigma_{g_{\mathcal H},k_{\mathcal H}}$.

With $\P(S)$ denoting the pants graph of $S$, we build on the work of Brock \cite{Bro01} to show: 

\begin{customthm}{E}\label{geombounds}
Let $N$ be a Seifert-fibered space over a hyperbolic surface $S$ and $\bar\Gamma \subset N $ be a hyperbolic link. If $\mathcal H = \{(\bar\Gamma_i, S_i)\}_{i = 1}^n$ is a stratification of $\bar\Gamma$, then there exist constants $K_1 > 1$ and $K_0>0$, depending only on $\kappa(\mathcal H)$, and pants decompositions $P_{X_i},P_{Y_i}$ of $S_{\mathcal H}$ with $\Gamma_i^{\mathcal H} \subset P_{X_i}, \Gamma_{i+1}^{\mathcal H}  \subset P_{Y_i}$ for $i = 1, \ldots, n$, such that:
$$ \f 1 {nK_1}\left(\sum_{i=1}^n d_{\mathcal P(S_{\mathcal H})}(P_{X_i},P_{Y_i})\right) \leq vol(N_{\bar\Gamma})\leq K_1\left(\sum_{i=1}^n d_{\mathcal P(S_{\mathcal H})}(P_{X_i},P_{Y_i})\right)+nK_0$$
where $\Gamma_{n+1}^{\mathcal H} = \psi_{\mathcal H}(\Gamma_{1}^{\mathcal H}) \subset P_{Y_n}$.
\end{customthm}

The pants decompositions in Theorem \ref{geombounds} arise by taking the covers associated to the stratifing surfaces and choosing shortest pants decompositions on their conformal boundaries. Therefore, they intrinsically depend on the geometry of $N_{\bar\Gamma}$. However, we can obtain the following purely topological result.
\begin{customthm}{F}\label{combbounds}
Let $N$ be a Seifert-fibered space over a hyperbolic surface $S$ and $\bar\Gamma \subset N $ be a hyperbolic link. If $\mathcal H = \{(\bar\Gamma_i, S_i)\}_{i = 1}^n$ is a stratification of $\bar\Gamma$, then there exist constants $K_1 > 1$ and $K_0>0$, depending only on $\kappa(\mathcal H)$, such that:
 $$\frac 1 {nK_1} \inf_{\mathfrak P}\left( \sum_{i=1}^n d_{\mathcal P(S_{\mathcal H})}(P_i,P_{i+1})\right) \leq vol(N_{\bar\Gamma})\leq K_1\inf_{\mathfrak P}\left( \sum_{i=1}^n d_{\mathcal P(S_{\mathcal H})}(P_i,P_{i+1})\right)+nK_0$$
where $\mathfrak P = \{P_i\}_{i = 1}^n$ is any collection of pants decompositions with $\Gamma_i^{\mathcal{H}} \subset P_i$ and we let $P_{n+1} = \psi_{\mathcal H}(P_1)$. \end{customthm}

Our main interest in developing this technology is to apply it to the study of canonical link complements in $PT(S)$. In this setting, we can often remove the dependence on the stratification complexity $\kappa$ in Theorems \ref{geombounds} and \ref{combbounds} by using the following result.

 \begin{customprop}{G}\label{canonical} Let $S$ be a hyperbolic punctured surface and $\Gamma \subset S$ a collection of essential closed curves in minimal position. If $\hat{\Gamma}$ is stratified by horizontal surfaces and a component of $\Gamma$ is simple, then there exists a stratification $\mathcal{H}$ of $\hat\Gamma$ by sections of $\pi : PT(S) \to S$. In particular, all curves in $\Gamma$ are simple and non-zero in $H_1(S)$.
\end{customprop}
\subsection{Examples of stratified canonical links}

A lingering question -- the existence of stratified hyperbolic canonical links -- is addressed in Section \ref{examples}, where we give infinite families of stratified canonical links. One of our observations is that given a filling pair of multicurves, we can always add punctures to the surface to make them stratified. This observation is also crucial in proving Theorem \ref{pair}.

Stratified collections of multicurves are \rev{straightforward  to construct using flat metrics arising from both quadratic differentials $q \in Q(X)$ and from foliations of pseudo-Anosov diffeomorphisms. Let $cyl(q)$ denote the set of cylinder curves for $q$, then}

\begin{customthm}{H}\label{qdiff} Let $X$ be a complex structure on $\Sigma_{g,m}$. Then every $q \in Q(X)$ with $k$ singularities gives rise to a collection $cyl(q)$ of essential simple closed curves on $\Sigma_{g,k}$ such that $\{ \Gamma \subset cyl(q) \mid \Gamma \text{ is finite and } \hat{\Gamma} \text{ is stratified and  hyperbolic in } PT(\Sigma_{g,k})\}$ is infinite.
\end{customthm}

In Remark \ref{psuedo-anosovex}, we also build examples coming from pseudo-Anosov maps $\phi\in \mathrm{Mod}(\Sigma_{g,m})$ with $k$-singularities, such that there are infinitely many pairs $(\gamma_1, \gamma_2)$ of essential simple closed curves on $\Sigma_{g,k}$ where the canonical lift of $\Gamma_n = (\phi^n(\gamma_1), \gamma_2)$ is stratified and hyperbolic for $n \geq K$, which depends on $\Sigma_{g,k}$ and $d_{\mathcal{C}(\Sigma_{g,k})}(\gamma_1, \gamma_2)$. Further, there are constants $A > 1$ and $B > 0$, depending only on $\Sigma_{g,k}$, such that: \[\frac{1}{A}\,  d_{\mathcal C(\Sigma_{g,k})}(\phi^n(\gamma_1),\gamma_2)-B \leq vol (N_{\hat{\Gamma}_n}) \quad \text{for} \quad n \geq K.\]
\rev{Note, by Theorem \ref{pair}, the above lower bound holds for any pair $(\gamma_1, \gamma_2)$ as long as $n$ is large enough to guarantee that $(\phi^n(\gamma_1),\gamma_2)$ is filling. As we will see in Section \ref{examples}, many of these pairs are also stratified.}

As seen in Proposition \ref{canonical}, not all canonical lifts can be stratified. We give further non-trivial examples in the following proposition.

\begin{customprop}{I}\label{nogen}
Let $S$ be a hyperbolic punctured surface different from  $\Sigma_{1,1}$ or $\Sigma_{0,4}$. Then there exists a pair $(\alpha, \beta)$ of essential simple closed curves where $(\hat\alpha, \hat \beta)$ is unstratifiable.
\end{customprop}

\brem We conclude by pointing out that Seifert-fibered spaces $N$ which contain stratifiable links can be obtained as mapping tori $N \cong M^{\psi_{\mathcal H}}  \eqdef   S_{\mathcal H} \times I / (x,0) \sim (\psi_{\mathcal{H}}(x), 1)$ where $\psi_{\mathcal H}$ has finite order. Modifying the definition of a stratification by asking that all stratifying surfaces induce fibrations, our results should \rev{ extend} to mapping tori $M^\psi$ in which there are no restriction on the mapping class $\psi$.\erem

\paragraph{Outline.} In Section 2, we recall some background and in Section 3 we study the nature of stratification surfaces.  In Section 4, we prove the general volume bounds of Theorem \ref{geombounds} and Theorem \ref{combbounds}.  Section \ref{main proof} is devoted to the proof of our main Theorem \ref{pair}.  Section 6 contains constructions for stratified links and the proof of Theorems \ref{punctorus}, \ref{ce} and \ref{length}. 

\paragraph{Acknowledgments}The first and second author would like to thank Ian Biringer for inspiring conversations. The second author thanks Pekka Pankka for discussions on these topics. The second author gratefully acknowledges support from the Academy of Finland project 297258 ``Topological Geometric Function Theory" also from the grant 346300 for IMPAN from the Simons Foundation and the matching 2015-2019 Polish MNiSW fund.  The third author would like to thank David Gabai for helpful discussions and Princeton University for continued support.

 \section{Preliminaries}\label{setup}
 In the following sections we recall some facts and definitions about the topology of surfaces and 3-manifolds. For references, see \cite{He1976,Ha2007,Ja1980}.
 \begin{center}
 \emph{All objects are smooth and orientable throughout this paper unless otherwise stated.}
 \end{center}
 
 \subsection{3-manifolds} 
 
An orientable 3-manifold $M$ is said to be \emph{irreducible} if every embedded sphere $\mathbb S^2$ bounds a 3-ball. A map between manifolds is said to be proper if it sends boundaries to boundaries and pre-images of compact sets are compact. A connected properly immersed surface $\Sigma \imm M$ is \emph{$\pi_1$-injective} if the induced map on the fundamental groups is injective. Furthermore, if $\Sigma \hookrightarrow M$ is properly embedded, two-sided
, and $\pi_1$-injective we say that it is \emph{incompressible}. If $\Sigma \hookrightarrow M$ is a non $\pi_1$-injective two-sided surface, then by the Loop Theorem there is a compressing disk $\mathbb D^2\hookrightarrow M$ such that $\partial \mathbb D^2=\mathbb D^2\cap \Sigma$ and $\partial \mathbb D^2$ is non-trivial in $\pi_1(\Sigma)$. 

An irreducible 3-manifold $(M,\partial M)$ is said to have \emph{incompressible} \emph{boundary} if every map of a disk $(\mathbb D^2,\partial \mathbb D^2)\hookrightarrow (M,\partial M)$ is homotopic via a map of pairs into $\partial M$. Therefore, $(M,\partial M)$ has incompressible boundary if and only if each component $\Sigma \in \pi_0( \partial M)$ is incompressible (i.e. $\pi_1$-injective).  An orientable, irreducible and compact $3$-manifold is called \textit{Haken} if it contains a two-sided $\pi_1$-injective surface. A 3-manifold is said to be \emph{acylindrical} if every map $(\mathbb S^1\times I,\partial (\mathbb S^1\times I))\rar (M,\partial M)$ can be homotoped into the boundary via maps of pairs. 

An \emph{essential surface} $(\Sigma,\partial\Sigma)\hookrightarrow (M,\partial M)$ is an incompressible, $\partial$-incompressible properly embedded surface. An \emph{essential loop} $\gamma$ in a surface $\Sigma$ is a non-peripheral $\pi_1$-injective loop.

A compact $3$-manifold $M$ is \emph{hyperbolic} if its interior admits a complete metric of constant negative sectional curvature. Similarly, a link $L \subset M$ is \emph{hyperbolic} if $M \setminus L$ is a hyperbolic $3$-manifold. By the Geometrization  Theorem \cite{Per2003.1,Per2003.2,Per2003.3}, a compact 3-manifold is hyperbolic if and only if it is irreducible, contains no essential tori and $\pi_1$ is infinite.

\rev{
\subsection{Surfaces}

We let $\Sigma_{g,k}$ be a smooth orientable surface of genus $g$ with $k$ punctures. \rev{The projectivized unit tangent bundle of $S$ will be denoted by $PT(S)$ and for any smooth path $\al$ on $S$, we will use $\hat{\al} \subset PT(S)$ for the set of tangent lines to $\al$. We call $\hat{\al}$ the \emph{canonical lift} of $\al$}.

 A \emph{multicurve} on $S$ is a disjoint union of non-parallel essential simple closed curves. The set of multicurves is denoted by $\mathfrak M(S)$. Note, when we speak of collections of multicurves, we do allow for two multicurves in a collection to have parallel components.

In this paper, we let $\mathcal T(S)$ denote the Teichm\"uller space of $S$ and let $\overline{\mathcal T(S)}$ the Weil-Petersson completion, where we use $d_{WP}$ to denote the Weil-Petersson metric. The boundary $\overline{\mathcal T(S)}\setminus\mathcal T(S)$ is the set of \emph{noded surfaces}. Let $\alpha\in\mathfrak M(S)$ be a simple multi-curve, we denote by $S(\alpha)\subset \overline{\mathcal T(S)}\setminus\mathcal T(S)$  the strata of noded surfaces in which $\alpha$ is pinched and we have $S(\alpha)\cong \mathcal T(S\setminus\alpha)$. See \cite{wolpert_2003} for a reference on this material.

Lastly, we use $\mathrm{Mod}(S)$ to denote the mapping class group of $S$, which can be defined as $\mathrm{Mod}(S) := \mathrm{Diffeo}^+(S)/\mathrm{Diffeo}_0^+(S)^+$, where $\mathrm{Diffeo}_0^+(S)$ is the component of the identity.}

 \subsection{Curve and pants graph} 

We will make use of two important graphs associated to a surface, for details see \cite{Schleimer}. Let $S$ be a hyperbolic surface different from $\Sigma_{1,1},\Sigma_{0,4}$ and $\Sigma_{0,2}$. Then, the  \emph{curve graph} $\mathcal C(S)$ of $S$ is the metric graph obtained by taking $\mathcal S$, the set of free homotopy classes of essential simple closed curves in $S$, as vertices and letting $\alpha,\beta\in \mathcal S$ be connected by an edge if they can be realized disjointly on $S$. For $\Sigma_{1,1}$ and $\Sigma_{0,4}$, a similar definition exists where $d_{\mathcal C(S)} (\alpha,\beta)=1\Leftrightarrow \iota(\alpha,\beta)\leq 1$ or $\iota(\alpha,\beta)\leq 2$, respectively, where $\iota(\cdot, \cdot)$ is the geometric intersection number.

\begin{figure}
\begin{overpic}[scale=.3]{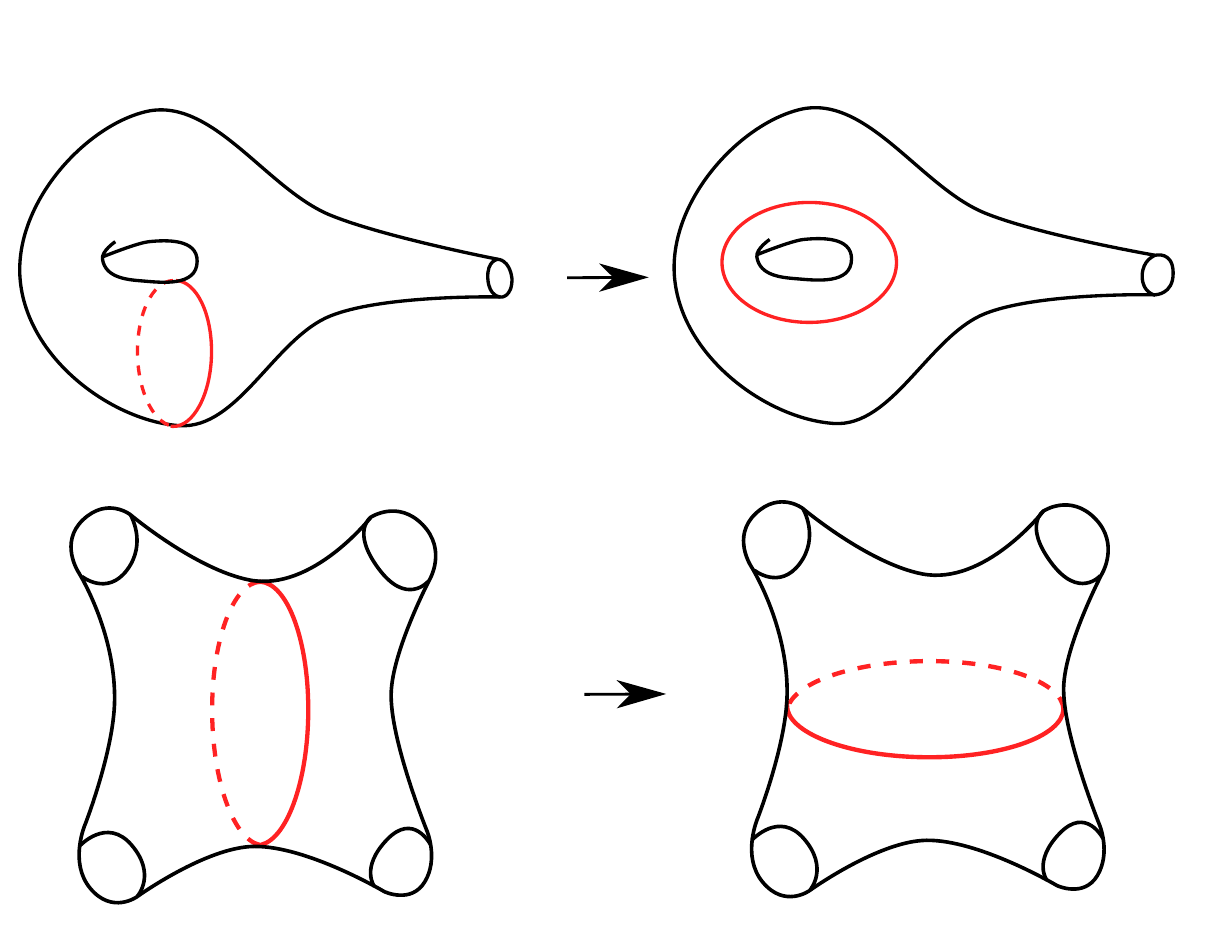}
\end{overpic}
\caption{Types of elementary moves between pants decompositions.}
\label{fig:pm}
\end{figure}

We will also use of the \emph{pants graph} $\mathcal P(S)$, see \cite{Schleimer} for details.  A multicurve is a \emph{pants decomposition} of $S$ if it has the maximal possible number of components, which is equal to $\kappa(S) = 3g -3 + k$, the complexity of $S$. The vertices of $\mathcal P(S)$ are isotopy classes of pants decompositions of $S$ and two pants decompositions have distance $1$ if and only if there is an elementary move taking one to the other. An \emph{elementary move} replaces a curve $\gamma$ in a pants decomposition with another curve $\gamma'$ that is disjoint from all other pants curves such that either (1) $\iota(\gamma, \gamma') = 2$ when $\gamma$ is in the boundary of two distinct pairs of pants or (2) $\iota(\gamma, \gamma') = 1$ when $\gamma$ is in the boundary of the same pair of pants,  \rev{  see Figure \ref{fig:pm}}. We let $d_{\mathcal P(S)}$ denote the induced metric distance in the pants graph.

\rev{\bdefi\label{ssp}
Suppose $Y$ is a non-annular, connected, incompressible subsurface of $S.$ The \emph{subsurface projection} $proj_Y: \mathcal C(S)\rightarrow  \mathcal C(Y)$ is a course map defined as follows. Given $\alpha\in \mathcal C(S)$, we isotope $\alpha$ to minimize the number of connected components of $\alpha\cap Y.$ The choice of $proj_Y(\al)$ then falls into several cases:
\begin{enumerate}
  \item If $\alpha\subset S\setminus Y$, we set $proj_Y(\alpha)=\emptyset$.
  \item If $\alpha\subset Y$, we set $proj_Y(\alpha)=\alpha$.
  \item If $i(\alpha,\partial Y)>0$, we pick an arc $\alpha' \subset \alpha \cap Y$ and let $V$ be a closed regular neighborhood of $\alpha' \cup\partial Y$. We then set $proj_Y(\alpha)$ equal to any essential and non-peripheral component $\alpha''$ of $\partial V$ in $Y$. 
\end{enumerate}
Note, we call $proj_Y$ a coarse map because the choice in $(3)$ is not unique, however, the distance in $\mathcal C(Y)$ between the different choices for a fixed $\al$ is at most $1$. For a collection $\mathfrak S$ of simple closed curves in $S$ let $d_Y(\mathfrak S)$ denote the diameter of $proj_Y(\mathfrak S)$ in $\mathcal C(Y)$. 
\edefi}

We can now state an important formula of Masur and Minsky:

\bthm[{\cite{MM}, Theorem 6.12}]\label{masurminsky} There is a constant $C_0(S)$ such that for any $C\geq C_0(S)$, there exist constants $c_1 > 1$ and  $c_0 > 0$ with the property
$$\frac 1 {c_1} d_{\mathcal P(S)}(P_1,P_2)-c_0\leq \sum_{Y\subset S} [d_Y(P_1 \cup P_2)]_C\leq c_1 d_{\mathcal P(S)}(P_1,P_2)+c_0$$
for any $P_1, P_2\in\mathcal P(S)$ and the sum is over all non-annular connected incompressible $Y \subset S$. \rev{ The notation $[a]_C$ is defined as $[a]_C = 0$ if $a  \leq C$ and $[a]_C = a$ otherwise.} \ethm

 \subsection{Line fields}
 \bdefi\label{lf} A \emph{line field} $\xi$ on $S$ is a section of $\pi: PT(S) \to S$. The punctures of $S$ are called singularities of $\xi$. We will often conflate $\xi$ with its image surface in $PT(S)$. Recall that line fields are always integrable by the Frobenious Theorem.
\edefi

Let $x$ be a singularity of a line field ${\xi}$ on $S$. Hopf associates to such a singularity a half-integer $\ind_{{\xi}}(x)\in \frac{1}{2}\mathbb{Z},$ defined as follows. Pick a small closed curve $\al :[0,1]\rightarrow S$ bounding a disk around $x$ with no other singularities on or inside $\al$. Further, make $\al$ small enough to be contained in some coordinate chart. Let $\vec{v}(0)$ be a vector at $\al(0)$ parallel to $\xi(\al(0))$ and pick a continuous extension of $\vec{v}(t)$ at $\al(t)$ such that $\vec{v}(t)$ is parallel $\xi(\al(t))$ for all $t \in [0,1)$. We can measure the number of total rotations of $\vec{v}(\al(t))$ relative to the local coordinates as we traverse $\al$. Hopf shows that this gives a half-integer $\ind_{{\xi}}(x)\in \frac{1}{2}\mathbb{Z}$ which is independent of the various choices involved. The interested reader can reference (Chapter III, \cite{H89}) for more information and pictures of line fields and possible singularities.

We recall the following well-known result in the theory of line fields over closed surfaces:

\bthm[Poincar\'e-Hopf]\label{P-H}
If ${\xi}$ is a line field on $\Sigma_{g,k}$ with singularities $x_1,...,x_k$, then
\[
    \sum_{i=1}^{n} \ind_{{\xi}}(x_i)=2-2g.
\]

\ethm

Lastly, we mention the operation of \emph{inflating} $\xi$ along a simple closed curve $\gamma$ with $\hat{\gamma} \subset \xi$. To inflate, we simply cut along $\gamma$ and replace it with a cylinder whose line field is parallel to the boundary. Gluing this cylinder back and smoothing produces a new smooth line field.

\section{Stratified links in Seifert-Fibered spaces} \label{strat}

Let $N$ be a Seifert-fibered space over a hyperbolic surface $S$ and $\bar\Gamma \subset N$ a link. In this paper, we restrict to the sub-class of links $\bar\Gamma$ that arise as lifts of a collection $\Gamma$ of essential closed curves in \emph{minimal position} on $S$. Drilling a regular neighborhood of $\bar\Gamma$, we obtain \rev{ $N_{\bar\Gamma} = N \setminus \bar{\Gamma}$}.  The following result of the first two authors tells us when $N_{\bar\Gamma}$ is hyperbolic.

\blem[\cite{CM19}]\label{hypfill} Let $\Gamma$ be a collection of essential closed curves on a hyperbolic surface $S$ in minimal position and let $N$ be a Seifert-fibered space over $S$. Then a lift $\bar\Gamma \subset N$ of $\Gamma$ is hyperbolic if and only if $\Gamma$ is filling on $S$ and $\bar\Gamma$ is acylindrical.
\elem

In our setting, where we work with curves in minimal position, acylindrical amounts to ``separating'' any two lifts of parallel components by lifts of transverse curves. To obtain our volume bounds on $N_{\bar\Gamma}$, we will need some additional control over the lifts $\bar\Gamma$. 

\begin{customdef}{1}\label{1}  A link $\bar\Gamma \subset N$ is \emph{stratified} if one can decompose $\overline{\Gamma} $ into sub-links $ \set{\bar\Gamma_j}_{i = 1}^n$ such that there is a disjoint properly embedded collection of incompressible surfaces $\{S_i\}_{i = 1}^n \subset N$ with $\bar\Gamma_i$ embedded in $S_i$ for each $i$. Otherwise, $\bar\Gamma$ is said to be \emph{unstratifiable}.

The $S_i$ are called \emph{stratification surfaces} and we assume that $\{S_i\}_{i = 1}^n$  is \emph{cyclically ordered around the $\S$-fiber}. We call the collection of pairs $\mathcal H =  \{(\bar\Gamma_i, S_i)\}_{i = 1}^n$ a \emph{stratification} of $\bar\Gamma$. \end{customdef}

Our first fact about stratifications follows from general properties of of Seifert-fibered spaces.

 \blem\label{horizontalsurface}
 Let $N$ be a Seifert-fibered space over $S$ and let $\bar\Gamma\subset N$ be a stratified hyperbolic link. Then, each stratifying surface $S_i$ is isotopic to an essential horizontal surface.
 \elem
 \bpf Since each $S_i$ is an essential properly embedded surface, Proposition \cite[1.11]{Ha2007} gives $S_i\iso F_i$ where $F_i$ is either a vertical torus or a horizontal surface. Assume that $F_1$ is a vertical torus. It follows that all the $S_i$'s must be vertical tori as they are disjoint. We then see that $\beta_i = \pi(S_i)$ is a simple closed curve on $S$ for each $i$ and that the $\beta_i$'s are all disjoint. Hence, $\bar\Gamma$ is a lift of $\{\beta_i\}_{i = 1}^n$. However, by Lemma \ref{hypfill}, $N \setminus \bar\Gamma$ cannot be hyperbolic as $\{\beta_i\}_{i = 1}^n$ \rev{ is not filling,} a contradiction. Thus, $F \iso S_i$ must be horizontal for all $i$. \epf
 As a direct application of Proposition \cite[1.11]{Ha2007}, we further have:
 
 \bcor\label{makehor} Let $N$ be a Seifert-fibered space over $S$ and let $\bar\Gamma\subset N$ be a stratified hyperbolic link. Then, there is an isotopy of $N$ such that all $S_i$ are horizontal surfaces. Moreover, if $\bar\Gamma$ is transverse to the fibers, we can assume that the isotopy is rel $\bar\Gamma$.\ecor 

\brem For $N$ to contain horizontal surfaces, $S$ must either be punctured or, more generally, the Euler number $e(N)$ must be zero.\erem

Given a stratification $\mathcal H =  \{(\bar\Gamma_i, S_i)\}_{i = 1}^n$, cut $N$ along $S_1$ to see that $S_i \iso_f S_1$ for all $i$, were $\iso_f$ denotes fiber-wise isotopy. Let $S_{\mathcal H} \eqdef  S_1$ and use the \rev{ fiberwise ``downward'' projection} to define $\Gamma_i^{\mathcal H} \iso_f  \bar\Gamma_i$. Note that $\Gamma_i^{\mathcal H}$ is a multicurve for each $i$. Let $\psi_{\mathcal H}$ be the ``bottom-to-top'' monodromy of the cut along $S_{\mathcal H}$, so $N$ is naturally homeomorphic to the mapping torus $S_{\mathcal H} \times I / (x,0) \sim (\psi_{\mathcal{H}}(x), 1)$ and $\psi_{\mathcal H}$ is a periodic element of the mapping class group of $S$. Lastly, define the \emph{complexity} of a stratification, where $S_{\mathcal H} \cong \Sigma_{g_{\mathcal H},k_{\mathcal H}}$, as $\kappa(\mathcal H) = 3g_{\mathcal H} - 3 + k_{\mathcal H}$.

\blem\label{minpos} Let $N$ be a Seifert-fibered space over $S$ and let $\bar\Gamma\subset N$ be a stratified link with stratification $\mathcal H =  \{(\bar\Gamma_i, S_i)\}_{i = 1}^n$ by horizontal surfaces. Then, up to isotopy of $\bar\Gamma$ in $\bigcup_{i = 1}^n S_i \subset N$, $\{\Gamma_i^{\mathcal H}\}_{i = 1}^n$ is in minimal position on $S_{\mathcal H}$.
\elem

\bpf Since $N \setminus S_{\mathcal H} \iso S_{\mathcal H} \times I$ and $\bar\Gamma_i$ is a multicurve on $S_i$ for each $i$, it follows easily that any Reidemeister type II or III move applied to $\{\Gamma_i^{\mathcal H}\}_{i = 1}^n$ can be realised by an isotopy in $\bigcup_{i = 1}^n S_i \subset S_{\mathcal H} \times I$, and therefore in $N$.\epf

\subsection{Stratified Canonical Links}

We now focus on the case where $N=PT(S)$. Care needs to be taken when dealing with homotopy classes of closed curves and canonical lifts. In particular, introducing self-intersections will change the isotopy class of the lift.

 \bdefi A homotopy $h:\mathbb S^1\times I\rar S$ between two self-transverse closed curves is \emph{transversal} if $h_t$ is self-transverse for all $t\in I$.
 \edefi
 
 A transversal homotopy $h_t$ between closed curves $\alpha$ and $\beta$ in $S$ induces an isotopy in $PT(S)$ between $\hat\alpha$ and $\hat\beta$, see \cite[Sec. 2]{M17}. Hass-Scott proved that any two self-transverse minimal representatives  are transversally homotopic, in particular to a unique closed geodesic for a fixed non-positively curves metric on $S$. Thus, we have:
 
 \blem\label{lem:trans} Let $\Gamma_1$, $\Gamma_2$ be two collections of essential closed curves on $S$ in minimal position such that $\phi(\Gamma_1)=\Gamma_2$ for some $\phi \in \Diff(S).$  Then, $N_{\widehat\Gamma_1}$ is homeomorphic to $N_{\widehat\Gamma_2}$.
 \elem  
 
 Therefore, when $\Gamma$ is a collection of essential closed curves in minimal position, it is enough to look at the geodesic representative of $\Gamma$ \rev{ in a hyperbolic metric}  on $S$ or even on a singular flat metric on $S$, provided the geodesics avoid the singularities.  Further, this guarantees that whether or not $\hat\Gamma$ is stratifiable is a well defined property of $\Gamma$. 
 
\bcor\label{cor:stratwell} Let $\Gamma_1$, $\Gamma_2$ be two collections of essential closed curved on $S$ in minimal position such that $\phi(\Gamma_1)=\Gamma_2$ for some $\phi \in \Diff(S).$  If $\hat\Gamma_1$ is stratified, then so is $\hat\Gamma_2$. 
 \ecor

We now turn to analyzing properties of stratifications of canonical lifts.
 
 \blem\label{nonnullhom}
Let $S = \Sigma_{g,k}$ with $k > 0$ and let  $\gamma\subset S$ be an essential simple closed curve. Then there exists a properly embedded, horizontal, incompressible surface $Z \subset PT(S)$ with $\hat\gamma \subset Z$  if and only if $[\gamma]\neq0\in H_1(S)$.
 \elem
 
 \bpf ($\Leftarrow$) Assume that $[\gamma] \neq 0$. Let $X$ be a complex structure on $S$ and consider a Jenkins-Strebel differential $q$ on $X$ associated to $\gamma$. Recall that $q$ is a meromorphic quadratic differential with simple poles at the punctures and has the property that all non-singular leaves of the foliation $\mathcal F_{h,q} = \ker \mathrm{Im}\sqrt{q}$ are closed and isotopic to $\gamma$, see \cite{S84} for details. One should think of $\mathcal F_{h,q}$ as a cylinder with intervals on the boundary identified to build $S$. 
 
 {\bf Case 1.} $\gamma$ is non-separating in $S$. Since all singular leaves of $\mathcal F_{h,q}$ are horizontal and finite length, the union of all singular leaves is a connected graph containing all singularities. Collapsing the edges of this graph via Whitehead moves, we can assume that all the singularities of \rev{  $\mathcal F_{h,q}$ live} at the punctures of $S$. 
 
 {\bf Case 2.} $\gamma$ is separating in $S$. By assumption, $[\gamma] \neq 0 \in H_1(S)$, so each component of $S \setminus \gamma$ contains a puncture. In this case,  the union of all singular leaves of $\mathcal F_{h,q}$ has two components corresponding to the pieces of $S \setminus \gamma$. Collapsing the edges of this graph via Whitehead moves, we assume that all the singularities of \rev{ $\mathcal F_{h,q}$  live} at the punctures of $S$.
 
Thus, the line field tangent to $\mathcal F_{h,q}$ is a section $S \to PT(S)$ whose image contains $\hat{\gamma}$.
 
 ($\Rightarrow$) Assume that $[\gamma] = 0$ and there is a properly embedded, horizontal, incompressible surface $Z \subset PT(S)$ with $\hat{\gamma} \subset Z$. Since $[\gamma] = 0$, it follows that $S \setminus \gamma$ has a component $Y$ that is a compact surface with $\partial Y = \gamma$ and genus $> 0$. By construction, the projection $p \eqdef \left.\pi\right|_Z : Z \to S$ is a covering map, so $\partial (p^{-1}(Y)) = p^{-1}(\gamma)$ and every component of $p^{-1}(Y)$ is a surface of genus $> 0$.
 
 From another perspective, let $U \subset Z$ be small enough such that $\left. p \right|_U$ is a homeomorphism. Then $U \subset PT(S)$ defines a line field on $p(U)$. Pulling back these line fields via $p$ over all such $U$ gives a line field $\eta$ on $Z$. Since $\hat\gamma \subset Z$ is a canonical lift, $\eta$ must be tangent to $\hat\gamma \subset Z$. Now, consider the restriction of $\eta$ on $p^{-1}(\gamma)$. Since $Z$ is embedded, $p^{-1}(\gamma)$ is a disjoint collection of loops in the torus $\pi^{-1}(\gamma) \subset PT(S)$. In particular, $\eta$ must be transverse to all components of $p^{-1}(\gamma)$ except for $\hat{\gamma}$. Now, let $Y_0$ be a component of $p^{-1}(Y)$. Then $\left. \eta \right|_{Y_0}$ is a line field on $Y_0$ that is transverse to $\partial Y_0 \setminus \hat\gamma$ and tangent to $\hat\gamma$ whenever $\hat\gamma \subset \partial Y_0$. In either case, we can cap off $\partial Y_0$ with disks extending $\eta$ to have singularities of index $+1$ in each disk. Let $Y_1$ be the closed surface obtained by capping off $Y_0$. Then, by Poincar\'e-Hopf  \ref{P-H}, we must have that \rev{ the number of boundary components of $\partial Y_0$ is $ \chi(Y_1)$.} However, $\chi(Y_1) \leq 0$ since $Y_1$ has genus $g > 0$ but  \rev{ the number of boundary components is greater than zero}, which is a contradiction. It follows that $\hat\gamma$ cannot live on $Z$.\epf 

 \begin{customprop}{G}\label{canonical} Let $S$ be a hyperbolic punctured surface and $\Gamma \subset S$ a collection of essential closed curves in minimal position. If $\hat{\Gamma}$ is stratified by horizontal surfaces and a component of $\Gamma$ is simple, then there exists a stratification $\mathcal{H}$ of $\hat\Gamma$ by sections of $\pi : PT(S) \to S$. In particular, all curves in $\Gamma$ are simple and non-zero in $H_1(S)$.
\end{customprop}

\bpf Let $\{S_i\}_{i = 1}^n$ be the given horizontal stratification surfaces for $\hat{\Gamma}$ and let $\gamma_1$ be the simple closed curve in $\Gamma$. By the proof of Lemma \ref{nonnullhom}, $[\gamma_1] \neq 0 \in H_1(S)$ and we can find a line field $\xi \subset PT(S)$ containing $\hat{\gamma}_1$. Moreover, we have that $\gamma_1$ is the \emph{unique} component of $\Gamma$ tangent to $\xi$. Further, up to inflating of $\xi$, we can assume that $\xi$ is parallel and tangent to $\gamma_1$ in a regular neighborhood of $\gamma_1$ on $S$. Let $F$ be any stratifying surface.

Up to fiber-wise isotopy, we make $F$ transverse to $\xi$ and, by deleting a small neighborhood of the punctures, we assume that $S$, $PT(S)$, $F$, and $\xi$ all have boundary. We will now think of $\gamma_1$ as our first component in the sub-link decomposition, so we have: $\Gamma_1=\set{\gamma_1}$.

{\bf Claim 1.} $F$ is obtained, up to fiber-wise isotopy, by surgery on parallel copies of $\xi$ over a collection $\mathfrak M$ of pairwise disjoint properly embedded essential arcs in $\xi$ such that $\mathfrak M \cap \hat\Gamma_1 = \emp$.

\bpfc Since $\xi$ and $F$ both contain $\hat\Gamma_1$, there is an open regular neighborhood $U \subset \xi$ of $\hat\Gamma_1$ and an isotopy of $N$ rel $\hat\Gamma$ such that $U \subset F$.  Note, $U$ is an annulus and, as $F$ is embedded, $\pi^{-1}(U) \cap F$ is a collection of parallel annuli over $U$. Further, we can assume that $U$ is small enough so that $\xi$ on $\pi(U)$ is parallel to $\Gamma_1$ and $\partial U$ is  \rev{ the set of tangent lines to (i.e. the canonical lift of)} $\partial \pi(U)$. Since $\xi$ has boundary, we can build $\xi$ from $U$ by only attaching $1$-handles $D_j$ along arcs $\mathfrak m_j, \mathfrak m_j'$ on $\partial U$. Notice that $\pi^{-1}(D_j)$ is solid torus and $\pi^{-1}(D_j) \cap F$ is a collection of parallel discs. It follows that $F$ is obtained by annual surgery along the collection $\mathfrak M = \bigcup_j \{\mathfrak m_j, \mathfrak m_j'\}$ in $\xi$. By construction, these arcs are disjoint from $\Gamma_1 \subset U$ and each $\mathfrak m \in \mathfrak M$ is the canonical lift of $\pi(\mathfrak m)$ since $\mathfrak m \in \partial U$.
\epfc

 {\bf Claim 2.} For $\hat{\gamma} \in \hat{\Gamma}$, either $\hat\gamma \subset \xi$ or $\hat\gamma \cap \xi = \emp$.
 
\begin{figure}[htb]
\begin{overpic}[scale=1.2] {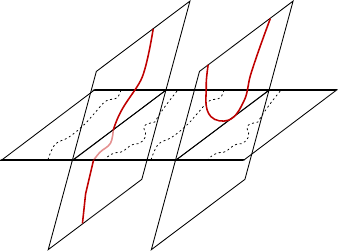}
\put(90,33){$\xi$}
\put(95,52){$\partial \xi$}
\put(57.5,69){$V \subset F$}
\put(27,15){$\hat\gamma$}
\end{overpic}
\caption{Local picture of $F \cap R \cap \pi^{-1}(\pi(W))$. Dotted lines are $\partial W$ and vertical annuli are part of $F$. The red arcs are parts of $\hat\gamma$.}
\label{R2}
\end{figure}

\bpfc Assume that $\hat\gamma \cap \xi \neq \emp$ for some $\hat{\gamma} \in \hat{\Gamma}$. Let $U$ and $\mathfrak M$ be as in the proof of Claim 1. Since $\hat{\gamma} \subset F$, we can assume that $\hat{\gamma}$ avoids $U$ by shrinking it as necessary.  By the construction in Claim 1, $\xi\setminus U$ is a disjoint union of $1$-handles $h_i$'s with $h_i\cong [0,1]_s\times [0,1]_t$ and $(s,\set{0,1})\subset \mathfrak M$. Moreover, the fibre structure induced by $\cup_t [0,1]_s\times\set{t}$ is the canonical lift of their projection $\pi([0,1]_s\times\set{t})$. 

By transversality  $\mathfrak N \eqdef  (F \setminus U) \cap \xi$ is a disjoint union of arcs, each of which is an essential arc in some $h_i\setminus \mathfrak M$. Therefore, up to isotopy, since each essential arc $\alpha$ in $h_i\setminus \mathfrak M$ is a canonical lift of $\pi(\mathfrak m)$, we can assume the same for all elements of $\mathfrak N$. Let $W \subset \xi$ be a regular neighborhood of $\mathfrak N$ in $\xi$ and let $R \subset N$ be a regular neighborhood of $\xi$ in $N$ such that $F \cap R \cap \pi^{-1}(\pi(W))$ is a collection of disjoint strips, see Figure \ref{R2}.

\begin{figure}[h]
\begin{overpic}[scale=0.8] {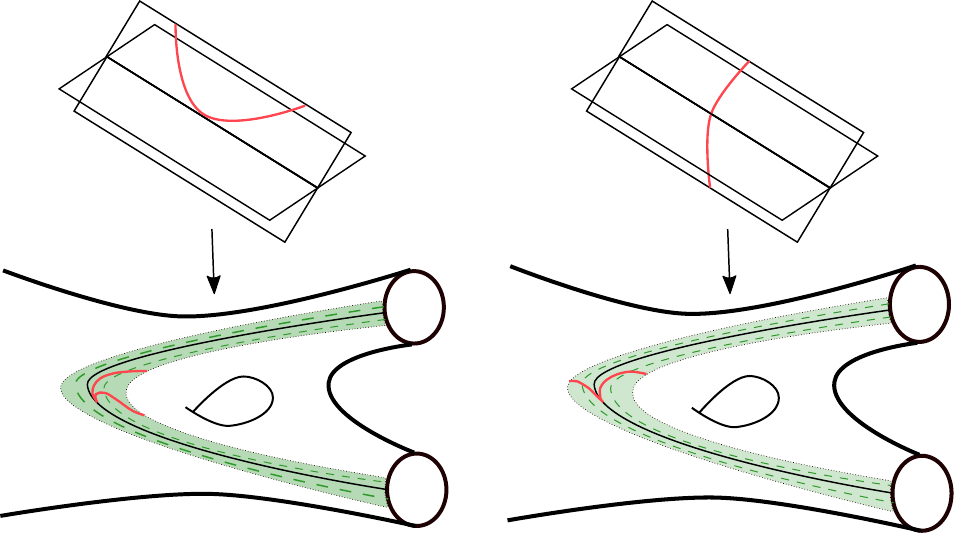}
\end{overpic}
\caption{The two possible projection of $\delta$ to $S$. Neither is smooth.} \label{thm5v2}
\end{figure}

Since $\hat\gamma \subset F$, the arcs of $\hat\gamma$ that intersect $\xi$ must pass into these strips. Let $V \subset F$ be such a strip containing an arc $\delta$ of $\hat\gamma$. By making $R$ small enough, we can assume that either $\delta$ crosses $V$ or $\delta$ dips though $V$ to be become tangent to $\xi$, see Figure \ref{thm5v2}.
 
Since $F$ is horizontal, $\pi(V)$ is a neighborhood $\pi(\mathfrak n)$ for some arc $\mathfrak n  \in \mathfrak N$. Because $\mathfrak n$ is the canonical lift of $\pi(\mathfrak n)$, it follows that $\xi$ is an almost parallel line field in $\pi(V)$ tangent to $\pi(\mathfrak n)$. Consider $\pi(\delta)$, which must be a smooth arc with canonical lift $\delta$. However, in the two possible cases for $\delta \subset V$, the projection $\pi(\delta)$ will have a cusp at $\pi(\delta \cap\mathfrak n)$ tangent to $\pi(\mathfrak n)$, which contradicts smoothness, see Figure \ref{thm5v2}.\epfc

By Claims 1 and 2, we can cut away the annular surgery that produces $F$ from parallel copies of $\xi$ without cutting $\hat\Gamma$. Thus, $\hat\Gamma$ lives on parallel copies of $\xi$, which give us the stratification by sections. It follows that all curves in $\Gamma$ must be simple and, by Lemma \ref{nonnullhom}, they are all non-zero in $H_1(S)$.
\epf

We now describe the obstructions for pairs of simple closed loops to being stratifiable.

 \blem\label{pucture-strat} Let $S$ be a hyperbolic punctured surface and $\alpha,\beta$ a filling pair of essential simple closed curves on $S$ in minimal position. Then $\hat\alpha,\hat\beta$ is stratified if and only if every simply connected component of $S\setminus (\alpha\cup\beta)$ is a rectangle.
\elem

\begin{proof}  Let $\{D_i\}_{i = 1}^n$ be the connected components of $S \setminus \al \cup \beta$. Note that $D_i$ is either a $2m_i$-gon or a punctured $2m_i$-gon for some $m_i \in \mathbb N$ as $\al$ and $\beta$-arcs have to alternate on $\partial D_i$.

 \begin{figure}[h]
\begin{overpic}[scale=0.5]{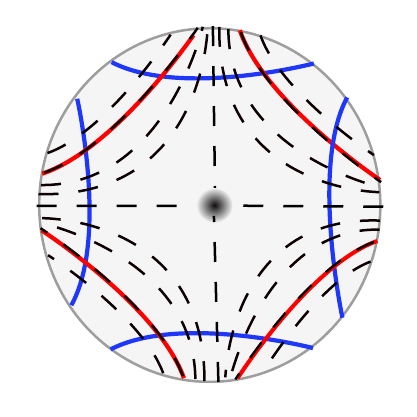}
\end{overpic}
\caption{The $2m_i$-gon around $x_i$ and the line field ${\xi},$ where red arcs belong to $\alpha$ and blue ones to $\beta.$ }\label{ind}
\end{figure}

Assume that $\hat\alpha,\hat\beta$ is stratified and let, by Proposition \ref{canonical}, $\xi$ be the line field on $S$ tangent to $\al$. When $D_i$ is punctured, let $x_i$ be the puncture in $D_i$. Otherwise, let $x_i$ be any point in $D_i$. Since \rev{ $\xi$ is tangent} to $\al$ and transverse to $\beta$, we can use $\partial D_i$ to compute $\ind_{\xi}(x_i) = (2- m_i)/2$, see Figure \ref{ind}. If $D_i$ does not contain a puncture, Poincar\'e-Hopf \ref{P-H} tells us that $\ind_{\xi}(x_i) = 0$, so $m_i = 2$ and $D_i$ is a rectangle.

\begin{figure}[htb]
\begin{overpic}[scale=0.7]{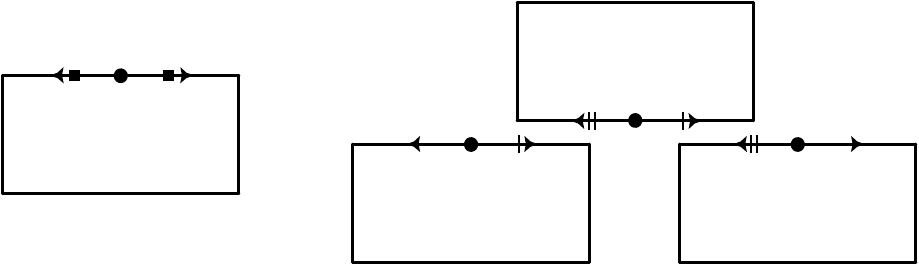}
\put(10,5){$\al$-arc}
\put(-10,13){$\beta$-arc}
\end{overpic}
\caption{Examples of a right angled $2$-gon and $6$-gon with cone points.}
\label{flat-gons}
\end{figure}

Conversely, assume that all simply connected components of $S \setminus \al \cup \beta$ are rectangles. We will build a flat metric on $S$ with singularities exactly at the punctures. If $D_i$ is a rectangle, give $D_i$ the flat metric of a unit square. If $D_i$ is a punctured $2m_i$-gon, we can can realise $D_i$ as a flat right-angled $2m_i$-gon with sides of a length $1$ and a cone point of angle $m_i \pi$ at the puncture, see Figure \ref{flat-gons}. Gluing these flat metrics along the edges gives a flat metric on $S$. Up to rotating by $\pi/2$, we can assume that $\al$ is a horizontal curve in this metric and $\beta$ is vertical. The horizontal and vertical line fields give a stratification $\hat\al$, $\hat\beta$.\end{proof}

\bcor\label{pair-punc} Let $\alpha,\beta$ a filling pair of multicurves in minimal position on $\Sigma_{g,m}$. Then puncturing all $k$ non-rectangular simply connected components of $\Sigma_{g,m} \setminus \al \cup \beta$ makes $\hat{\al}$, $\hat\beta$ stratified in $PT(\Sigma_{g,m+k})$ while keeping $\al$, $\beta$ filling on $\Sigma_{g,m+k}$. Further, $k  \leq 2g - 2$.
\ecor

\begin{proof} \rev{Aside form }the bound on $k$, this is a direct corollary of the proof of Lemma \ref{pucture-strat}. To see that $k  \leq 4 g - 4$, we fill in the $m$ punctures and let $\mathcal{D}_j$ denote the collection of regions of $\Sigma_g \setminus (\alpha\cup\beta)$ that are bounded by a $j$-gon. Since $\al$ and $\beta$ are simple, $j=2i$ for some $i\geq 2.$ Because the number of regions is finite, we let $2M$ be the maximal number of sides of a component of $\Sigma_g \setminus (\alpha\cup\beta)$. We extend the graph $\alpha\cup\beta$ to a triangulation on $\Sigma_g$ by doing barycentric subdivision on each complementary polygon. With $V$, $E$, $F$ the number of vertices, edges, and faces of this triangulation, we see that

 \[2g-2=-F+E-V=-\left(\sum_{i=2}^{M} 2i \, \sharp\mathcal{D}_{2i} \right) + \left(\sum_{i=2}^{M} \left(\frac{2i}{2} + 2 i \right)\sharp\mathcal{D}_{2i} \right)-\left(\sum_{i=2}^{M} \left(\frac{2i}{4} + 1 \right)\sharp\mathcal{D}_{2i} \right)=\] \[ =\left(\sum_{i=2}^{M} \frac{i-2}{2}\sharp\mathcal{D}_{2i} \right)\geq \sum_{i=3}^{M} \frac{1}{2}\sharp\mathcal{D}_{2i} \geq \frac{k}{2} .\]\end{proof}

\brem We can extend Corollary \ref{pair-punc} to include pairs $\al$, $\beta$ of essential multicurves curves where instead of minimal position we assume that $\al \cup \beta$ only requires Reidemeister type II moves to get to minimal position. Then, puncturing the associated embedded bigons, we can stratify the curves. However, the bound on $k$ will also depend on the number of bigons. \rev{ In the more general setting of non-mininal position families of curves and their associated canonical lift complements, one could consider sequences of bigons as ``twist regions,'' which play an important role in contolling the geometry in the classical setting of hyperbolic links. In fact, the studying of hyperbolic links whose diagrams project onto surfaces with genus has been of recent interest, see \cite{kalfagianni2020alternating, adams2021generalized, bavier2021guts, champanerkar2020volumish, will2020homological}}.\erem

 \section{Volume Bounds}\label{volumebounds}
 
Here, we will prove the two main volume estimates \rev{ in Theorem \ref{geombounds} and Theorem \ref{combbounds}.}

 \bdefi Let $ccvol(M)$ denote the convex-core volume of a hyperbolic 3-manifold $M$. Recall, that $ccvol(M)<\infty$ if and only if $M$ is geometrically finite. Let $\overline M$ denote the compactification of $M$ obtained by the Tameness Theorem \cite{AG2004,CG2006}.  \rev{ Consider a collection $P\subset \partial \overline M$ comprised of pairwise disjoint essential simple closed curves on $\partial \overline M$ together with all the tori components of $\partial \overline M$.} Define $GF(\overline M,P)$ to be the space of geometrically finite structures on $M$ for which $\pi_1(P)$ is parabolic for each structure in $GF(\overline M,P)$.
\edefi

\brem By \cite{MT1998}, $GF(\overline M,P)$ is parametrized by the Teich\"uller space $\mathcal T(\partial \overline M\setminus P)$. This parametrization is realized via the conformal boundary at infinity of \rev{$Y \in GF(\overline M,P)$.}
\erem

  \blem\label{lem1} For $M \in GF(\overline S \times I , \partial \overline S \times \{1/2\} \cup \alpha\times\{0\}\cup\beta\times\{1\})$, where $\alpha$ and $\beta$ are multicurves, there are constants $K_1 \geq 1$ and $K_0 \geq 0$, depending only on  $S$, and  pants decompositions $P_X, P_Y$ of $S$ with  $\alpha\subset P_X,\beta\subset P_Y$, such that:
  $$\frac 1 {K_1} d_{\mathcal P(S)}(P_X,P_Y)-K_0\leq ccvol (M)\leq K_1 d_{\mathcal P(S)}(P_X,P_Y)+K_0.$$
 \elem

\bpf
By the Bers Density Theorem \cite{Ma2007}, there is a sequence $Q_n\cong Q(X_n,Y_n)$ of quasi-Fuchsian manifolds with geometric limit $M$ where $X_n, Y_n \in \mathcal T(S)$. This sequence corresponds to pinching $\alpha$ on $X_n$ and $\beta$ on $Y_n$. \rev{By \cite[Theorem 1.2]{Bro01} and the discussion at the end of \cite{Bro01} for punctured surfaces}, there are constants $K_1 \geq 1$ and $K_ 0 \geq 0$, depending only on $S$, such that for all $n\in\mathbb N$:
 $$\frac 1 {K_1} d_{\mathcal P(S)}(P_{X_n},P_{Y_n})-K_0\leq ccvol (Q_n)\leq K_1 d_{\mathcal P(S)}(P_{X_n},P_{Y_n})+K_0$$
where the $P_{X_n},P_{Y_n}\in\mathcal P(S)$ are the Bers pants decompositions corresponding to the conformal structures $X_n,Y_n$, respectively (see  \cite[Theorem 12.8]{FM2011}  for a reference). Since $Q_n\rightarrow M$ geometrically we have that $ccvol(Q_n) \to ccvol(M)$. Further, $X_n\rightarrow X$ and $Y_n\rightarrow Y$ for some nodal surfaces $X$ and $Y$, with $\alpha$ and $\beta$ pinched, respectively. Thus, for large enough $n$, the pants decompositions $P_{X_n},P_{Y_n}$ are constant and equal to $P_X, P_Y$, which contain $\alpha$ and $\beta$, respectively. \epf

Let $\mathcal H = \{(\bar\Gamma_i, S_i)\}_{j = 1}^n$ be a stratification of a hyperbolic link $\bar\Gamma\subset N$. As every $S_i$ is two-sided and $N$ is orientable, we can make a canonical choice of parallel copy $\Sigma_i$ of $S_i$ obtained as the ``top'' boundary of a regular neighborhood of $S_i$. The choice can be made with all $\Sigma_i$'s disjoint. Let  $\pi_i : N_i \to N_{\bar\Gamma}$ be the cover of $N_{\overline \Gamma}$ corresponding to $\pi_1(\Sigma_i)$ \rev{with the pull-back metric}. The metric on $N_i$ is the pull-back of the hyperbolic metric on $N_{\bar\Gamma}$ via $\pi_i$.

 \blem\label{convexcores} If a hyperbolic link $\bar\Gamma\subset N$ admits a stratification $\mathcal H = \{(\bar\Gamma_i, S_i)\}_{i = 1}^n$, then 
 $$vol(N_{\bar\Gamma})\leq \sum_{i=1}^n ccvol(N_i).$$
 \elem
 \bpf Since $\bar\Gamma$ is a hyperbolic link, we can assume that each $S_i$ is an essential horizontal surface carrying $\bar\Gamma_i$, see Section \ref{strat} for details. It follows that $Y_i \eqdef  S_i \setminus \bar\Gamma_i$ is an incompressible, nonseparating surface in $N_{\bar\Gamma}$, so by results of \cite{CHMR17}, it has a least area representative $Y'_i$ in $N_{\bar\Gamma}$. Moreover, by \cite{FHS1983}, since $Y'$ has an embedded representative, it is itself embedded. Further, still by \cite{FHS1983}, $Y_i'$ and $Y_j'$ are disjoint for $i \neq j$ as they have disjoint homotopic representatives. 
 
Consider, the components $A_i$ of $N_{\bar\Gamma} \setminus \bigcup_{i = 1}^n Y'_i$. Notice that $\partial A_i = Y_i' \cup Y_{i+1}'$, with $n+1=1$, lifts homeomorphically to $N_i$, the cover constructed from $\pi_1(\Sigma_i)$.  By \cite{FHS1983}, since these lifts are least area surfaces, they are contained in $CC(N_i)$, the convex core of $N_i$. Thus, $A_i$ lifts to $\tilde A_i\subset CC(N_i)$ and $vol(A_i)\leq ccvol(N_i)$. As the $A_i$'s decompose $N_{\bar\Gamma}$, the result follows. \epf

 \begin{customthm}{E}\label{geombounds}
Let $N$ be a Seifert-fibered space over a hyperbolic surface $S$ and $\bar\Gamma \subset N $ a hyperbolic link. If $\mathcal H = \{(\bar\Gamma_i, S_i)\}_{i = 1}^n$ is a stratification of $\bar\Gamma$, then there exist constants $K_1 > 1$ and $K_0>0$, depending only on $\kappa(\mathcal H)$, and pants decompositions $P_{X_i},P_{Y_i}$ of $S_{\mathcal H}$ with $\Gamma_i^{\mathcal H} \subset P_{X_i}, \Gamma_{i+1}^{\mathcal H}  \subset P_{Y_i}$ for $i = 1, \ldots n$, such that:\small
$$ \f 1 {nK_1}\left(\sum_{i=1}^n d_{\mathcal P(S_{\mathcal H})}(P_{X_i},P_{Y_i})\right) \leq vol(N_{\bar\Gamma})\leq K_1\left(\sum_{i=1}^n d_{\mathcal P(S_{\mathcal H})}(P_{X_i},P_{Y_i})\right)+nK_0$$\normalsize
where $\Gamma_{n+1}^{\mathcal H} = \psi_{\mathcal H}(\Gamma_{1}^{\mathcal H}) \subset P_{Y_n}$.
\end{customthm}
 
We now prove Theorem \ref{geombounds} by proving two lemmas each of which shows one of the bounds.

\blem\label{geomub} If $\mathcal H = \{(\bar\Gamma_i, S_i)\}_{i = 1}^n$ is a stratification of a hyperbolic link $\bar\Gamma\subset N$, then there exist constants $K_1 > 1$ and $K_0>0$ depending only on $\kappa(\mathcal H)$, and pants decompositions $P_{X_i},P_{Y_i}$ of $S_{\mathcal H}$ with $\Gamma_i^{\mathcal H} \subset P_{X_i}, \Gamma_{i+1}^{\mathcal H}  \subset P_{Y_i}$ for $i = 1, \ldots n$, such that:
$$vol(N_{\bar\Gamma})\leq K_1\left(\sum_{i=1}^n d_{\mathcal P(S_{\mathcal H})}(P_{X_i},P_{Y_i})\right)+nK_0$$\normalsize
where $\Gamma_{n+1}^{\mathcal H} = \psi_{\mathcal H}(\Gamma_{1}^{\mathcal H}) \subset P_{Y_n}$.
\elem
\bpf
By Lemma \ref{convexcores}, we know that:
$$vol(N_{\bar\Gamma})\leq \sum_{i=1}^n ccvol(N_i)$$
By the Covering Theorem \cite{Ca1996}, all the $N_i$'s are geometrically finite and by the Tameness Theorem, \cite{AG2004,CG2006}, we know that the manifold compactification $\overline N_i$ of $N_i$ is homeomorphic to $F_i\times I$ for some surface $F_i$. Moreover, since the surface $\Sigma_i\subset N_{\bar\Gamma}$ lifts homeomorphically to $N_i$ and induces a cusp preserving homotopy equivalence, it must be that $F_i \cong \Sigma_i\cong S_{\mathcal H}$. Therefore, \[N_i\in GF(\overline{S}_{\mathcal H} \times I, \partial \overline{S}_{\mathcal H} \times \{1/2\} \cup \Gamma_i^{\mathcal H}\times\{0\}\cup\Gamma_{i+1}^{\mathcal H}\times\{1\})\] and Lemma \ref{lem1} implies that for all $i$:
$$ccvol(N_i)\leq K_1 d_{\mathcal P(S_{\mathcal H})}(P_{X_i},P_{Y_i})+K_0.$$
Hence, the result follows. \epf

\blem\label{geomlb} If $\mathcal H = \{(\bar\Gamma_i, S_i)\}_{i = 1}^n$ is a stratification of a hyperbolic link $\bar\Gamma\subset N$, then there exist constants $K_1 > 1$ and $K_0>0$ depending only on $\kappa(\mathcal H)$, and pants decompositions $P_{X_i},P_{Y_i}$ of $S_{\mathcal H}$ with $\Gamma_i^{\mathcal H} \subset P_{X_i}, \Gamma_{i+1}^{\mathcal H}  \subset P_{Y_i}$ for $i = 1, \ldots n$, such that:
$$ \frac{1}{n K_1} \left(\sum_{i=1}^n d_{\mathcal P(S_{\mathcal H})}(P_{X_i},P_{Y_{i+1}})\right) \leq vol(N_{\bar\Gamma})$$
where $\Gamma_{n+1}^{\mathcal H} = \psi_{\mathcal H}(\Gamma_{1}^{\mathcal H}) \subset P_{Y_n}$. Moreover, $K_1$ and $P_{X_i},P_{Y_i}$ are the same as in Lemma \ref{geomub}.
\elem
\bpf
The idea is to follow the proof of \cite[4.1]{Bro01} and use short geodesics to bound the volume of $N_{\bar\Gamma}$ from below.

Fix $0<L<\min\{\mu_3,L_2\}$, where $\mu_3$ is the 3-dimensional Margulis constant and $L_2$ is the Bers constant for $S_{\mathcal H}$. Denote by $\mathcal G_L=\mathcal G_L(N_{\bar\Gamma})$ the set of homotopy classes of geodesics in $N_{\bar\Gamma}$ of length less than $L$. Since $N_{\bar\Gamma}$ is a finite volume hyperbolic 3-manifold, by Lemma \cite[4.8]{Bro01}, there is a constant $C_1 > 1$ depending only on $L$, such that:
$$\frac{\vert \mathcal G_L\vert}{C_1}  < vol(\mathcal N_{\epsilon}(N_{\bar\Gamma}^{\geq \epsilon})) < vol(N_{\bar\Gamma}),$$
where $\mathcal N_{\epsilon}(N_{\bar\Gamma}^{\geq \epsilon})$ is the $\epsilon$-neighborhood of the $\epsilon$-thick part of $N_{\bar\Gamma}$, with $0 < \epsilon < L/2$. As in Lemma \cite[4.2]{Bro01}, we need to bound $\vert \mathcal G_L\vert$ by pants distances.  Let $\mathcal S$ be the collection of homotopy classes of essential simple closed curves in $S_{\mathcal H}$. Given a proper embedding $h: S_{\mathcal H} \hookrightarrow M$, define
$$\mathcal S_L^h(M)=\left\{\alpha\in \mathcal S \mid \ell_{M}(h(\alpha))<L \right\}.$$
Since $S_{\mathcal H} \iso_f  S_i \iso_f  \Sigma_i $, we have well defined embeddings $h_i:S_{\mathcal H} \hookrightarrow N_i$. Then, by Lemma \cite[4.2]{Bro01} and the geometric limit argument of Lemma \ref{lem1},
$$d_{\mathcal P(S_{\mathcal H})}(P_{X_i},P_{Y_{i+1}})\leq K_1\vert \mathcal S_L^{h_i}(N_i)\vert$$
where $K_1$ is a constant depending only on $S_{\mathcal H}$. This gives:
$$\sum_{i=1}^nd_{\mathcal P(S_{\mathcal H})}(P_{X_i},P_{Y_{i+1}})\leq K_1\sum_{i=1}^n\vert \mathcal S_L^{h_i}(N_i)\vert.$$
Moreover, for every $\gamma \in\bigcup_{i=1}^n \mathcal S_L^{h_i}(N_i)$, we know that
 $$\ell_{N_{\bar\Gamma}}(\pi_i(h_i(\gamma)))\leq L,$$
 where $\pi_i : N_i \to N_{\bar\Gamma}$ are the covering maps. Thus, we get that $\bigcup_{i=1}^n \pi_i(\mathcal S_L^{h_i}(N_i))\subset\mathcal G_L$. To get a count, we need to see how many elements of $\bigcup_{i=1}^n \mathcal S^{h_i}_L(N_i)$ can be homotopic in $N_{\bar\Gamma}$.
  
  Since all elements of $\mathcal S_L(N_i)$ are homotopically distinct in $N_i$, their projections are distinct in $N_{\bar\Gamma}$. Thus, if $\pi_i(h_i(\alpha))\sim \pi_j(h_j(\beta))$ for some $\alpha, \beta, i,$ and $j$, then $i\neq j$. Therefore, 
   $$\frac 1 n\sum_{i=1}^n\vert \mathcal S_L^{h_i}(N_i)\vert\leq \vert\mathcal G_L\vert,$$
which gives
$$\frac 1 n\sum_{i=1}^nd_{\mathcal P(S_{\mathcal H})}(P_{X_i},P_{Y_{i+1}})<\frac {K_1}n\sum_{i=1}^n\vert \mathcal S_L^{h_i}(N_i)\vert\leq K_1\vert\mathcal G_L\vert$$
and we are done. \epf

\brem\label{pairwisefilling}

 If we assume that each pair $\Gamma_i^{\mathcal H}, \Gamma_{i+1}^{\mathcal H}$ fills $S_{\mathcal H}$, then $\pi_i(h_i(\alpha))\sim \pi_j(h_j(\beta))$ in $N_{\bar\Gamma}$ implies that $j$ must be equal to either $i+1$ or $i-1$, but not both. This is because if two loops $\alpha \in N_{i-1}$ and $\beta \in N_{i+1}$ are homotopic in $N_{\bar\Gamma}$, then $N_{\bar\Gamma}$ would not be acylindrical, contradicting hyperbolicity. Thus, if we let $\mathcal U_L = \bigcup_{i=1}^n \pi_i(\mathcal S_L^{h_i}(N_i))\subset \mathcal G_L$, we get that every element of $\mathcal U_L$ contains at most two projections of elements of the $\mathcal S_L^{h_i}(N_i)$.  Therefore,
  $$\frac 1 {2}\sum_{i=1}^n\vert \mathcal S_L^{h_i}(N_i)\vert\leq \vert\mathcal G_L\vert.$$
Hence,
 $$\frac 1 2\sum_{i=1}^nd_{\mathcal P(S_{\mathcal H})}(P_{X_i},P_{Y_{i+1}})<\frac {K_1}2\sum_{i=1}^n\vert \mathcal S_L^{h_i}(N_i)\vert\leq K_1\vert\mathcal G_L\vert.$$
\erem

By combining Lemmas \ref{geomub}, and \ref{geomlb}, we obtain Theorem \ref{geombounds}.

Before we can prove Theorem \ref{combbounds}, we will need the following quasi-isometry result about pants graphs of essential subsurfaces.

\begin{figure}[h]
\begin{overpic}[scale=0.5]{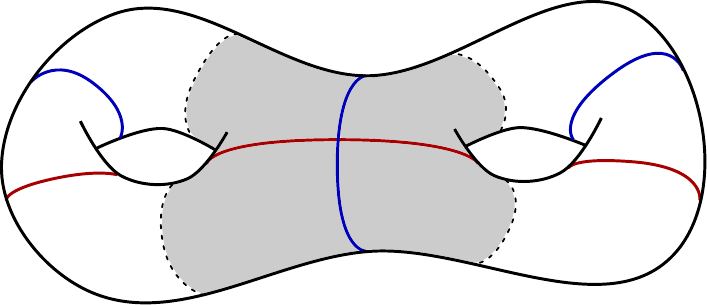}
\put(33,10){$V$}
\end{overpic}
\caption{\rev{Subsurface $V$ for pants decompositions in Proposition \ref{pantqi}, where the two pants decompositions  $P_1$, $P_2$ are in different colors.}}\label{non-annular}
\end{figure}

\bprop\label{pantqi} Let $P_1, P_2$ be pants decompositions of a hyperbolic surface $S$ that are in minimal position. Let $Y$ be the minimal incompressible subsurface containing $P_1 \cup P_2$ and define $V$ to be the union of all non-annular components of $Y$, see Figure \ref{non-annular}. Then there are constants $A > 1$ and $B> 0$, depending only on $S$, such that:
$$\frac{1}{A} d_{\mathcal P(V)}(P_1', P_2')  - B \leq  d_{\mathcal P(S)}(P_1, P_2) \leq  d_{\mathcal P(V)}(P_1', P_2')$$
where $P_i' = V \cap P_i$ are pants decompositions of $V$.
\eprop
\bpf Let $Q = P_1 \cap (S \setminus V) = P_2 \cap (S \setminus V)$ be the cores of the annular components that we removed. Consider the inclusion $\mathcal P(V)\hookrightarrow \mathcal P(S)$ given by $P'\mapsto P'\cup Q$. Then trivially
$$ d_{\mathcal P(S)}(P_1, P_2) \leq  d_{\mathcal P(V)}(P_1', P_2'),$$
which gives us the right-hand side of the desired inequality since $P_i' \cup Q = P_i$.

For the lower bound, recall that Theorem \ref{masurminsky} gives a constant $C_0(S)$ such that for any $C \geq C_0(S)$ there exist $c_1 > 1 $, and $c_0 > 0$ such that:
$$\frac 1 {c_1} d_{\mathcal P(S)}(P_1,P_2)-c_0\leq \sum_{Y\subset S} [d_Y(P_1 \cup P_2)]_C\leq c_1 d_{\mathcal P(S)}(P_1,P_2)+c_0$$
where the sum is over all non-annular connected incompressible subsurfaces $Y\subset S$. We will apply this formula to both $S$ and all possible subsurfaces $V$ to get the constants $A,B$.

{\bf Claim.} There exists $C>0$ such that if $[d_Y(P_1 \cup P_2)]_C\neq 0$ then $Y\subset V$.

\bpfc Let $proj_Y$ denote the subsurface projection of multicurves in $S$ to multicurves in $Y$. If $Y\not\subset V$ then $Q\cap Y\neq\emp$. Thus, $proj_Y(P_1)$ and $proj_Y(P_2)$ are both collections of coarsely defined pairwise disjoint simple closed curves containing $proj_Y(Q)$. Hence, the diameter $d_Y(P_1 \cup P_2)\leq 2$ and so the result follows by taking $C>2$. \epfc

Since the set of possible topological types of $V \subset S$ is finite, we can take \[C_0 = \max_{V} \{ C_0(V) \mid V \subset S \text{ is incompressible with no annular components}\}\]

Then, taking $C \geq \max(C_0, 3)$, for each $V$ we get $c_0(V)$ and $c_1(V),$ with the property that:
\begin{align*}
\frac 1 {c_1(V)} d_{\mathcal P(V)}(P_1',P_2')-c_0(V) &\leq  \sum_{Y\subset V} [d_Y(P_1' \cup P_2')]_C 
\xlongequal{\bf Claim}\\ \xlongequal{\bf Claim} \sum_{Y\subset S} [d_Y(P_1 \cup P_2)]_C &\leq c_1(S)  d_{\mathcal P(S)}(P_1,P_2)+c_0(S)
 \end{align*}
Thus,
$$ \frac 1 {c_1(V)c_1(S)} d_{\mathcal P(V)}(P_1',P_2')-\frac{c_0(V)+c_0(S)}{c_1(S)}\leq  d_{\mathcal P(S)}(P_1,P_2)$$
To make this hold for all $V$, we take the constants $$A\eqdef  \max_V  c_1(V)c_1(S)\qquad B\eqdef  \max_V\frac{c_0(V)+c_0(S)}{c_1(S)}.$$ Moreover, note that $Q$ was completely arbitrary and had no influence in the constants.\epf

We now prove our combinatorial bounds.

\begin{customthm}{F}\label{combbounds}
Let $N$ be a Seifert-fibered space over a hyperbolic surface $S$ and $\bar\Gamma \subset N $ a hyperbolic link. If $\mathcal H = \{(\bar\Gamma_i, S_i)\}_{i = 1}^n$ is a stratification of $\bar\Gamma$, then there exist constants $K_1 > 1$ and $K_0>0$, depending only on $\kappa(\mathcal H)$, such that:
 $$\frac 1 {nK_1} \inf_{\mathfrak P}\left( \sum_{i=1}^n d_{\mathcal P(S_{\mathcal H})}(P_i,P_{i+1})\right) \leq vol(N_{\bar\Gamma})\leq K_1\inf_{\mathfrak P}\left( \sum_{i=1}^n d_{\mathcal P(S_{\mathcal H})}(P_i,P_{i+1})\right)+nK_0$$
where $\mathfrak P = \{P_i\}_{i = 1}^n$ are any pants decompositions with $\Gamma_i^{\mathcal{H}} \subset P_i$ and $P_{n+1} = \psi_{\mathcal H}(P_1)$.
 \end{customthm}
 \bpf 
Taking an infimum in Lemma \ref{geomlb} gives the lower bound. The rest of the proof will focus on the upper bound.

Let $M = S_{\mathcal H} \times \S$ with the projection $\pi_{\mathcal H} : M \to S_{\mathcal H}$. Since $\phi_{\mathcal H}$ is periodic, we get a natural finite cover $q: M \to N$, where the degree $m$ is controlled by $\kappa(\mathcal H)$. Drilling the full preimage $\bar\Gamma' \eqdef  q^{-1}(\overline \Gamma)$ gives an $m$-degree cover of $N_{\overline \Gamma}$, which is hyperbolic, so $m \cdot vol(N_{\overline \Gamma}) = vol(M \setminus  \bar\Gamma')$.

Let $\mathfrak P = \{P_i\}_{i = 1}^n$ be a collection of pants decompositions of $S_i$ with $\bar\Gamma_i \subset P_i$. Our goal will be to use $\mathfrak P$ to build a hyperbolic link $\mathcal L\subset M$ containing $\bar\Gamma'$. Since $\bar \Gamma' \subset \mathcal L$, we know that $vol(M \setminus  \bar\Gamma') < vol(M \setminus \mathcal L)$, so it is enough to bound $vol(M \setminus \mathcal L)$.

\begin{figure}[htb]
\begin{overpic}[scale=0.4]{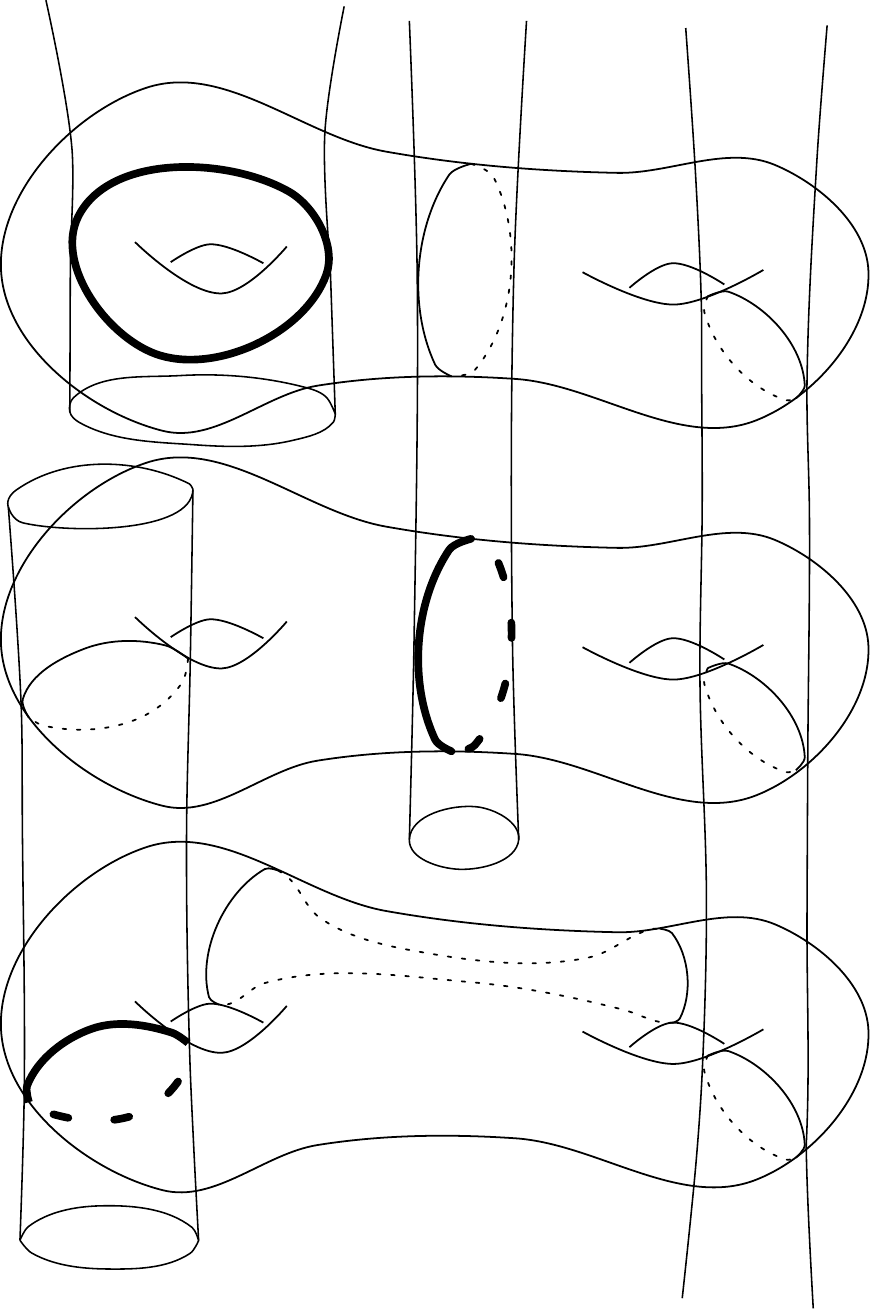}
\put(70,49){$S_i'$}
\put(70,78){$S_{i+1}'$}
\put(70,20){$S_{i-1}'$}
\put(33,79){$\gamma_2$}
\put(18.5,45.5){$\eta(A_{\gamma_2})$}
\put(13,84){$\gamma_1$}
\put(32.5,93){$A_{\gamma_2}$}
\put(13,97){$A_{\gamma_1}$}
\put(-22,15){$\gamma_3 = \eta(A_{\gamma_3})$}
\put(6,-1){$A_{\gamma_3}$}
\put(54.5,93){$A_{\gamma_4}$}
\end{overpic}
\caption{Examples of the cylinders $A_{\gamma_i}$. The curves $\gamma_1$, $\gamma_2$ and $\gamma_3$ are labeled, while $\gamma_4$ is outside the picture. The thick curves are elements of $\mathcal L$, while all the other curves are not included.}
\label{cylinders}
\end{figure}

{\it Step 1: Build $\mathcal L$}. Start with $\mathcal L = \bar\Gamma'$, we will add elements as follows. The lifts of the stratification and the pants decompositions can be ordered as 
\[\mathcal H'=\set{S_1, \ldots, S_n, \phi_{\mathcal H}(S_1), \ldots, \phi_{\mathcal H}(S_n), \phi_{\mathcal H}^2(S_1), \ldots, \phi_{\mathcal H}^{m-1}(S_n)}=\set{S_i'}_{i=1}^{n'}\]
\[\mathfrak P' = \set{P_1, \ldots, P_n, \phi_{\mathcal H}(P_1), \ldots, \phi_{\mathcal H}(P_n), \phi_{\mathcal H}^2(P_1), \ldots, \phi_{\mathcal H}^{m-1}(P_n)}=\set{P_i'}_{i=1}^{n'}\] giving a notion of height, where $n' = m \cdot n$. Over every component $\gamma$ of $\mathfrak P'$, there is a torus $T_\gamma = \gamma \times \S$. Let $O_\gamma = T_\gamma \cap \mathfrak P'$ be the curves in $\mathfrak P'$ that are contained in $T_\gamma$.  Since $M \setminus \bar\Gamma'$ is hyperbolic, $T_\gamma$ is hit transversely by $\bar\Gamma'  \subset \mathfrak P'$ at least once. Cut $T_\gamma$ into cylinders along the stratification surfaces where $\mathfrak P' \setminus O_\gamma$ hits $T_\gamma$ transversely, see Figure \ref{cylinders}. Whenever such a cylinder contains an element of $O_\gamma$, we add the element of lowest height in that cylinder to $\mathcal L$. Repeating this for all components of $\mathfrak P'$ gives $\mathcal L$. For each component $\gamma$ of $\mathfrak P'$, let $A_\gamma \subset T_\gamma$ be the cylinder containing $\gamma$ and $\eta(A_\gamma) \in\pi_0( \mathcal L)$ be the loop added to $\mathcal L$ from $A_\gamma$.

We still need to show that $\mathcal L$ is hyperbolic. Since $\bar\Gamma'$ is hyperbolic and stratified in $M$, Lemma \ref{minpos} implies that $\pi_{\mathcal H} (\bar\Gamma') \subset S_{\mathcal H}$ is a collection of essential closed curves in minimal position, possibly after an isotopy in $M$.  By Lemma \ref{hypfill}, it follows that $\pi_{\mathcal H} (\bar\Gamma')$ must be filling. As $\mathcal L$ contains $\bar\Gamma'$ and is stratified by the horizontal surfaces $\{S_i'\}_{i = 1}^{n'}$, we get that $\pi_{\mathcal H} (\mathcal L)$ is also in minimal position and filling after maybe an isotopy in $M$. Since $\mathcal L$ is acylindrical by construction, we apply Lemma \ref{hypfill} again to see that it must be hyperbolic. 

\begin{figure}[htb]
\begin{overpic}[scale=0.9] {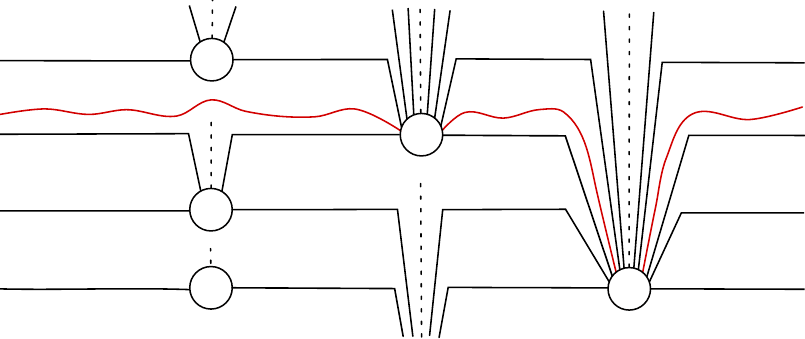}
\put(103,14.5){$S_{i-1}'$}
\put(103,33.5){$S_{i+1}'$}
\put(103,24.5){$S_i'$}
\put(103,5){$S_{i-2}'$}
\end{overpic}
\caption{Diagram of a surgery to the $Y_t$'s to make them properly embedded. The circles correspond to elements of $\mathcal L$, the dotted lines are the annuli $A_\gamma$, and the arcs leaving the circles are the surgered components. The surfaces are as in Figure \ref{cylinders} with the addition of $S_{i-2}'$. The red wavy arcs represent the $F_j$'s obtained from $\Sigma_i$.}
\label{surger}
\end{figure}

Let $r = m \cdot n \cdot \kappa(\mathcal H)$ and let $\mathcal Y = \{Y_t\}_{t = 1}^r$ be the connected components of $\bigcup_{i = 1}^{n'} (S_i' \setminus P_i')$. Note that each $Y_t$ is a pairs of pants. We can make $Y_t$ properly embedded in $M \setminus \mathcal L$ as any boundary component $\gamma$ of $Y_t$ that is not in $\mathcal L$ lies in the cylinder $A_\gamma$. Indeed, each $A_\gamma$ has a unique representative $\eta(A_\gamma) \in \mathcal L$ isotopic to $\gamma$ in $M \setminus \mathcal L$, so we can isotope $Y_t$ in a neighborhood of $A_\gamma$ to be properly embedded and $\S$-fiber transverse. Further, we can do this isotopy such that the collection $\mathcal Y$ is properly embedded and disjoint, see Figure \ref{surger}.

{\it Step 2: Decompose $M \setminus \mathcal L$}. Let $\Sigma_i$ be a positive push-off of $S_i'$ such that $\Sigma_i$ meets $\mathcal Y$ minimally. Then $\Sigma_i\cap \mathcal Y$ is a collection of simple closed curves and cutting $\bigcup_{i = 1}^{n'} \Sigma_i$ along $\mathcal Y$ gives a collection $\mathcal F=\set{F_j}_{j = 1}^s$ of non-annular, connected, incompressible subsurfaces of the $\Sigma_i$'s and a \rev{collection} of annuli. We can make $\mathcal F$ properly embedded in $M \setminus \mathcal L$ by continuing each $\partial F_j$ parallel to the cylinder in $\mathcal Y$ that cut it off. Further, we can do this such that $\mathcal F$ is disjoint from $\mathcal Y$, see Figure \ref{surger}. Notice that $\Sigma_i$ is cut exactly along $Q_i = \{\gamma \mid \gamma \in P_i' \text{ and } \gamma \in P_{i+1}'\}$, the isotopic simple closed curves shared between $P_i$ and $P_{i+1}$ when realized on $\Sigma_i$. This follows from the surgery construction since $\gamma \in P_{i+1}'$ is surgered along $A_\gamma$ if and only if it is isotopic to some curve in $P_i'$.

Let $\mathcal M = \{ M_j\}_{j = 1}^s$ be the connected component decomposition of $M \setminus \mathcal L$ cut along $\mathcal Y$. By construction, there is a bijection between $M_j$'s and $F_j$'s given by $F_j \subset M_j$. Further, for each $j$, we have that $M_j \simeq F_j \times I$. Notice that $\partial M_j \setminus (\partial F_j \times I)$ decomposes into two copies of $F_j$ where each copy carries a natural pants decompositions arising form the cut $\mathcal Y$. These pants decompositions are precisely $P(F_j)^+ = F_j \cap P_{i+1}'$ and $P(F_j)^- = F_j \cap P_i'$ since the $F_j$'s are components of $\Sigma_i \setminus Q_i$. Further, by the definition of $Q_i$, $P(F_j)^\pm$ must be transverse on $F_j$. These observations imply

{\bf Claim 1:} Either $M_j \simeq \Sigma_{0,3} \times I$ or $M_j \simeq F_j \times I$ with $\kappa(F_j) > 0$ and $\partial M_j \setminus (\partial F_j \times I)$ carries the transverse pants decompositions $P(F_j)^\pm$.

When $M_j \simeq \Sigma_{0,3} \times I$, the boundary is a pair of isotopic $Y_t$'s, so we can simply remove one of them from the collection $\mathcal Y$. Let $\mathcal Y^\vee = \{ Y_t \}_{t = 1}^{r'}$ be obtained from $\mathcal Y$ by selecting a unique representative from each isotopy class in $\mathcal Y$ rel $\partial (M \setminus \mathcal L)$. Similarly, we get $\mathcal M^\vee = \{M_j\}_{j = 1}^{s'}$ and $\mathcal F^\vee = \{F_j\}_{j = 1}^{s'}$ by cutting $M$ along $\mathcal Y^\vee$ and dropping all $F_j$ with $\kappa(F_j) = 0$, respectively. Note, we keep the indexing such that $F_j \subset M_j$ and $M_j \simeq F_j \times I$. We now have that $F_{j_1}$ and $F_{j_2}$ are isotopic in $M \setminus \mathcal L$ rel boundary if and only if $j_1 = j_2$.

{\it Step 3: Make the decomposition geometric}. Each $Y_t \in \mathcal Y^\vee$ has a minimal area representative $Y_t^\circ$, which by \cite{CHMR17} is a totally geodesic embedded \rev{thrice punctured} sphere. By  \cite{FHS1983}, since the $Y_t$'s are disjoint, their minimal representatives are disjoint or coincide, the latter happening precisely when they are homotopic. Further, since the $Y_t$'s are properly embedded essential surfaces in $M \setminus \mathcal L$, \rev{being homotopic is equivalent to being isotopic} by using Waldhausen's results \cite{Wa1968}. By construction, none of the representatives coincide and we get the collection $\mathcal Y^\circ = \{Y_t^\circ\}_{t = 1}^{r'}$ of disjoint totally geodesic thrice punctured spheres. Let $\mathcal M^\circ$ be the connected component decomposition of $M \setminus \mathcal L$ cut along $\mathcal Y^\circ$.

{\bf Claim 2:} Each element of $\mathcal M^\circ$ is isotopic to some $M_j \in \mathcal M^\vee$. Further, with $M_j^\circ \in \mathcal M^\circ$ isotopic to $M_j$, we have that $M_j^\circ$ is the convex core of some manifold in
\[GF(F_j \times I, \partial F_j \times \{1/2\} \cup P(F_j)^- \times \{0\} \cup P(F_j)^+ \times \{1\}).\]

\bpfc The isotopies of $Y_t$ to $Y_t^\circ$ can all be combined to an isotopy of $M \setminus \mathcal L$ taking $\mathcal Y^\vee$ to $\mathcal Y^\circ$, since the $\mathcal Y^\vee$ is a collection of disjoint, non-isotopic, properly embedded surfaces. Thus, each $M_j \in \mathcal M^\vee$ is taken to some element of $\mathcal M^\circ$. Since $\mathcal Y^\circ$ is a collection of totally geodesic thrice punctured spheres, we see that $\partial M_j^\circ$ is totally geodesic and maximally cusped, so $M_j^\circ$ is convex. By construction, the cusps are the elements of $P(F_j)^\pm$.
\epfc

{\it Step 4: Compute volume bound}. For each $\Sigma_i$, let $V_i$ be the union of the $F_j$'s with $\kappa(F_j) > 0$ that arise from $\Sigma_i$ cut along $Q_i$. Notice that $V_i$ is precisely the subsurface of $\Sigma_i$ obtained from the minimal incompressible subsurface containing $P_i \cup P_{i+1}$ by removing all annular components. By applying Lemma \ref{lem1} and Claim 2 to each $M_j^\circ$ ,we get that:
\[vol(M \setminus \mathcal L)=\sum_{j=1}^{r'}  vol (M_j^\circ)\leq \sum_{j=1}^{r'}\left(K_1\cdot d_{\mathcal P(F_j)}(P(F_j)^-,P(F_j)^+)+K_0\right)\]
Grouping all $F_j$'s to form $V_i$ and applying Proposition \ref{pantqi} gives, 
\[vol(M \setminus \mathcal L) \leq\sum_{i=1}^{n'} \left(K_1\cdot d_{\mathcal P(V_i)}(P(V_i)^-,P(V_i)^+)+K_0\right) \leq\sum_{i=1}^{n'} \left(K_3\cdot d_{\mathcal P(S_{\mathcal H})}(P_i',P_{i+1}')+K_2\right) \]
Since $P_{m \cdot k + i}' = \psi_{\mathcal H}^k(P_i)$ and $n' = m \cdot n$, we obtain 
\begin{align*}m \cdot vol(N_{\bar\Gamma})=vol(M \setminus  \bar\Gamma') &\leq vol(M \setminus \mathcal L) \\
& \leq\sum_{i=1}^n\left(\sum_{k = 1}^m \left(K_3\cdot d_{\mathcal P(S_{\mathcal H})}(\psi_{\mathcal H}^k(P_i),\psi_{\mathcal H}^k(P_{i+1}))+K_2\right)\right)\\
&=\sum_{i=1}^n\left(m\cdot \left(K_3\cdot d_{\mathcal P(S_{\mathcal H})}(P_i\,P_{i+1})+K_2\right)\right)\\
&=m\cdot\sum_{i=1}^n \left(K_3\cdot d_{\mathcal P(S_{\mathcal H})}(P_i,P_{i+1})+K_2\right)\end{align*}
where we use the invariance of $d_{\mathcal P(S_{\mathcal H})}$ under the mapping class group. By taking the infimum over $\mathfrak P$, we conclude 
\[vol(N_{\bar\Gamma})\leq K_3\inf_{\mathfrak P}\left( \sum_{i=1}^n d_{\mathcal P(S_{\mathcal H})}(P_i,P_{i+1})\right)+nK_2.\]\epf

 \section{Bounds for unstratifiable pairs and proof of Theorem \ref{pair}}\label{main proof}
 
In this section we prove our main result:
\begin{customthm}{A}\label{pair} Let $S$ be a hyperbolic surface, $N = PT(S)$, and let $(\al, \beta)$ be a filling pair of multicurves on $S$ in minimal position. Then there are constants $K_1 > 1$ and $K_0 > 0$, depending only on $S$, with
\[  \frac 1 {K_1} \inf_{P_{\al},P_{\beta}} d_{\mathcal P(S)}(P_{\al},P_{\beta}) - K_0\leq vol(N_{(\hat{\al}, \hat{\beta})}) \leq K_1 \inf_{P_{\al},P_{\beta}} d_{\mathcal P(S)}(P_{\al},P_{\beta}) + K_0,\]
where $P_{\al}, P_{\beta}$ are any pants decompositions of $S$ with $\al \subset P_{\al}, \beta \subset P_{\beta}$. Further, when $S$ has punctures, this bound is asymptotically optimal. Moreover, we have that $vol(N_{(\hat{\al}, \hat{\beta})})$ is also \rev{coarsely comparable} to the Weil-Petersson distance \rev{between} the strata $S(\alpha),S(\beta)\subset \bar{\mathcal T(S)}$.
\end{customthm}


We will prove the upper and lower bound in Proposition \ref{upperboundpair} and Theorem \ref{pairoffillingvol} respectively. We first start with a \rev{couple of lemmas} relating the pants distance in punctured surface and the corresponding filled surface.

\blem\label{pants0} Let $S$ be a hyperbolic surface and let $\alpha, \beta \subset S$ be a filling pair of multicurves in minimal position. If $S'$ is obtained from $S$ by removing a point $x$, then for any pants decomposition with $\alpha \subset Q_1$, $\beta \subset Q_2$, there are extensions to $Q_1'$, $Q_2'$ in $S'$ with
$$\frac 1 {2}d_{\mathcal P(S')}( Q_1', Q_2')\leq d_{\mathcal P(S)} (Q_1, Q_2)\leq  d_{\mathcal P(S')}(Q_1',Q_2').$$\elem
\bpf Since the set of simple closed curves $\mathcal S$ on $S$ has measure zero, we can assume there \rev{exists a point }$x \in S\setminus \mathcal S$. The right-hand inequality is obtained by taking a geodesic in $\mathcal P(S')$ between any choice of extensions $Q_1\subset Q_1'$, $Q_2\subset Q_2'$ and filling in the punctures. We now show that
\begin{equation}\label{dist}d_{\mathcal P(S')}(Q_1',Q_2')\leq 2\, d_{\mathcal P(S)} ( Q_1, Q_2)\end{equation}
for appropriately chosen extensions $Q_1'$ and $Q'_2$.  Let $P_1,\dotsc, P_n$ be a geodesic in $\mathcal P(S)$ between $Q_1=P_1$ and $Q_2=P_n$. We will inductively build a path in $\mathcal P(S')$ of length at most $2\, d_{\mathcal P(S)} (Q_1, Q_2)$ that will join some extensions $Q_1', Q_n'$ of $Q_1,Q_n$, respectively.

First, we build $P_1'$. Since $x\in S\setminus\mathcal S$, there is a well defined component $Y_i$ of $S \setminus P_i$ that contains $x$. Let $\delta_1$ to be an arc connecting $x$ to some $\beta_1 \subset \partial Y_1$ and let $\alpha_1$ be the loop obtained by taking the non-$\beta_1$ boundary component of a regular neighborhood of $\delta_1 \cup \beta_1$ in $S'$. Since $\alpha_1$ is essential in $S'$ and disjoint from $P_1$, we get a pants decomposition $P_1'\eqdef  \alpha_1\cup P_1$ of $S'$.

At each step of the induction, $P_i'$ will be obtained from an arc $\delta_i \subset Y_i$ that terminates in a curve $\beta_i \in \partial Y_i$ such that $P_i' \eqdef  P_i  \cup \al_i$ as above. Recall that to get from $P_i$ to $P_{i+1}$, we replace a curve $\gamma_i \in P_i$ by curve $\gamma_{i+1}$ such that $\gamma_{i+1}$ is disjoint from $P_i \setminus \gamma_i$ and either $\iota(\gamma_i, \gamma_{i+1}) = 2$, whenever $\gamma_i$ lies in the boundary of two distinct pairs of pants, or $\iota(\gamma_i, \gamma_{i+1}) = 1$, whenever $\gamma_i$ lies in the boundary of the same pair of pants. The former is called an $S$-move and the latter an $A$-move.

\begin{figure}[htb]
\begin{overpic}[scale=.5]{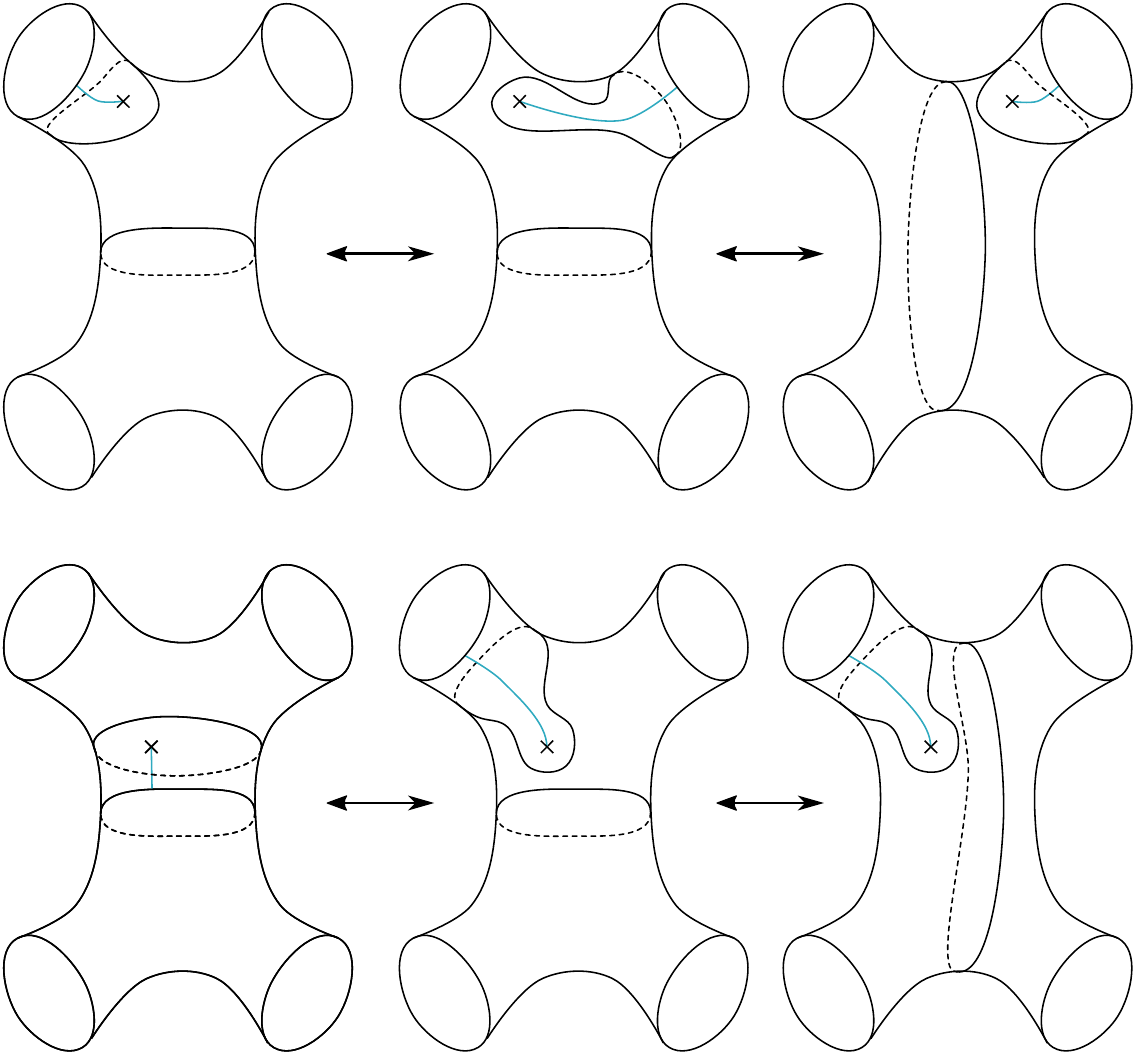}
\put(-13,70){$(i)$}
\put(-4,90){$\beta_i$}
\put(3.5,70){$\gamma_i$}
\put(80.5,53){$\gamma_{i+1}$}
\put(-13,20){$(ii)$}
\put(-5,20){$\gamma_i = \beta_i$}
\end{overpic}
\caption{The two possible obstructions to extending \rev{$S$-moves}.}
\label{fig:puncture-pants}
\end{figure}

{\bf Case.} $S$-move. By induction, we are given $P_i'$ and $\delta_i$. If $\iota_{S'}(\gamma_{i+1} ,\al_i)=0$, then we just take $\delta_{i +1} = \delta_i$ and $P_{i+1}'=P_{i+1}\cup \alpha_{i+1}$, so $d_{\mathcal P(S')}( P_i', P_{i+1}')=1$. Assume that $\iota_{S'}(\gamma_2 ,\al_i) \neq 0$. If $\beta_i \neq \gamma_i$, then we are in case $(i)$ from Figure \ref{fig:puncture-pants}. Here, $\beta_i \neq \gamma_i$ and we must move $\delta_i$ to $\delta_{i+1} \subset Y_{i+1}$. Let $P_i''$ be the pants decomposition where $\delta_{i+1}$ connects $x$ to the unique shared boundary component between $Y_i$ and $Y_{i+1}$. This operation takes one elementary move. Making the flip from $\gamma_i$ to $\gamma_{i+1}$ in $P_i''$ takes us to $P_{i+1}'$. This gives $d_{\mathcal P(S')}( P_i', P_{i+1}')=2$. Similarly, in the case where $\gamma_i = \beta_i$, we must move $\delta_i$ out of the way, which takes one move, and then flip $\gamma_i$ to $\gamma_{i+1}$. See Figure \ref{fig:puncture-pants} $(ii)$. We again get $d_{\mathcal P(S')}( P_i', P_{i+1}')=2$.

{\bf Case.} $A$-move. Essentially the same as in the case of $S$-moves and is left for the reader.

Building inductively gives a path in $\mathcal P(S')$ between $Q_1' = P_1'$ to $Q_2' = P_n'$ of length at most $2 \, d_{\mathcal P(S)}(Q_1,Q_2)$. Hence,
$$\label{dist}d_{\mathcal P(S')}(Q_1',Q_2')\leq 2 d_{\mathcal P(S)} ( Q_1, Q_2).$$
Notice that the extension $Q_2'$ depends on the initial choice of $Q_1$ and the chosen path between $Q_1$ and $Q_2$.
 \epf

\bcor\label{pants} Let $S$ be a hyperbolic surface and let $\alpha, \beta \subset S$ be a filling pair of multicurves in minimal position. If $S'$ is $S$ punctured at $m$ points in $S \setminus \al \cup \beta$, then
$$\frac 1 {2^m}\inf_{P_1',P_2'}d_{\mathcal P(S')}( P_1', P_2')\leq\inf_{P_1,P_2} d_{\mathcal P(S)} ( P_1, P_2)\leq \inf_{P_1',P_2'} d_{\mathcal P(S')}(P_1',P_2'),$$
where $P_1, P_2$ are pants decompositions of $S$ with $\al \subset P_1$, $\beta \subset P_2$ and $P_1', P_2'$ are pants decompositions of $S'$ with $\al \subset P_1'$, $\beta \subset P_2'$.\ecor
\bpf  For $m=1$, we use Lemma \ref{pants0} and take an infimum. Applying this inductively as we add one puncture at a time gives the desired result. \epf

\bprop\label{upperboundpair} Let $S$ be a hyperbolic surface, $N = PT(S)$, and let $(\al, \beta)$ be a filling pair of multicurves on $S$ in minimal position. Then there are constants $K_1 > 1$ and $K_0 > 0$, depending only on $S$, with
\[ vol(N_{(\hat{\al}, \hat{\beta})}) \leq K_1 \inf_{P_{\al},P_{\beta}} d_{\mathcal P(S)}(P_{\al},P_{\beta}) + K_0,\]
where $P_{\al}, P_{\beta}$ are any pants decompositions of $S$ with $\al \subset P_{\al}, \beta \subset P_{\beta}$. Further, when $S$ has punctures, this bound is asymptotically optimal.
\eprop

\bpf Consider $S\setminus (\alpha\cup\beta)$ and put a puncture $x_1,\dotsc, x_m$ in every disk region that is not a square to obtain a new surface $S'$. Note that, by Euler characteristic, $m$ is a bounded function of $\chi(S)$. By Lemma \ref{pucture-strat}, the canonical lift of $\Gamma = (\al, \beta)$ to $PT(S')$ is stratified by sections. Moreover, $\Gamma$ is still filling $S'$ and in minimal position, so $M\eqdef  PT(S')\setminus \hat\Gamma$ is hyperbolic by Lemma \ref{hypfill}. Since $M$ is obtained from $N_{\hat\Gamma}$ by drilling the vertical fibers corresponding to $x_1,\dotsc, x_k$, we see that $vol (N_{\hat\Gamma}) \leq vol (M)$. Applying Theorem \ref{combbounds} to $M$ gives the bound
$$vol (M)\leq  K_1\inf_{P_\alpha',P_{\beta}'}d_{\mathcal P(S')}(P_\alpha',P_{\beta}')+K_0,$$
which by Lemma \ref{pants} is equivalent to:
$$vol (M)\leq  2^m\cdot K_1\inf_{P_\alpha,P_{\beta}}d_{\mathcal P(S)}(P_\alpha,P_{\beta})+K_0.$$
When $S$ is punctured, this bound is asymptotically optimal as Theorems \ref{qdiff} \rev{produces} sequences $\Gamma_n = (\al_n, \beta_n)$ with $\hat\Gamma_n$ stratified, so \rev{there is no need} to add any extra punctures. \epf

We now prove the lower bound in \rev{Theorem \ref{pair}. Again}, we start with a few topological Lemmas and then move to adapt Canary's interpolation Theorem \cite{Ca1996} and some of Brock's results \cite{Bro01} to our setting showing that having funnels is not an obstruction.

\blem\label{incompsurfacegluing} For $n \geq 1$, let $H\eqdef  S_{g,n}\times I$ and consider $A\eqdef  \partial S_{g,n}\times I\subset \partial H$ in which we denote by $A_1,\dotsc, A_n$ each annular component. Let $V_1,\dotsc, V_n$ be solid tori and let $B_i\subset\partial V_n$ be regular neighborhoods of $(1,p_i)$ slopes with $\abs{p_i}>1$. Then the gluing
\[H'\eqdef  H\cup_f \left(\cup_{i=1}^n V_i\right) \quad \text{where} \quad  f: A\diffeo \cup_i B_i \quad \text{is given by} \quad  f\vert_{A_i}:A_i\diffeo B_i,\]
has incompressible boundary.
\elem 
\bpf  We recall that up to isotopy, we have the following facts:
\begin{enumerate}
\item each compressing disk in $S_{g,n}\times I$ intersects $A$ in $k$ vertical arcs with $k\geq 2$; 
\item each compressing disk in $V_i$ intersects a $(1,p)$ curve in $p$ points.
\end{enumerate}

\rev{Let $(D,\partial D)\hookrightarrow (H,\partial H')$} then, up to an isotopy of $D$, we can assume that $D\pitchfork A$ and that $\abs{\pi_0(\mathcal A)}$, where $\mathcal A = D\cap A$, is minimised. By an innermost argument we can find an arc $\alpha\in \mathcal A$ such that $\alpha\subset D$ co-bounds with $\beta\subset\partial D$ a disk $D'$ with $\mathcal A\cap\text{int}(D')=\emp$. Thus, $D'$ is either contained in some $V_i$ or in $S_{g,n}\times I$ but in either case we contradict (1) or (2). \epf

The following is just a reformulation of some previously proven results:

\bcor\label{vector field} \rev{Let $\alpha,\beta$ be a filling pair of multicurves in minimal position on a $\Sigma_{g,0}$. Then the number of non-square complementary regions $D_i$ is bounded in terms of $g$ and there exists a vector field $X$ on $\Sigma_{g,0}$ with a single singularity $x_i$ in each $D_i$ such that:
\begin{itemize}
\item $\abs{ind_{x_i}(X)}= (n_i/2)-2$ for $n_i$ the number of sides of $D_i$;
\item $\hat\alpha$ is parallel to $X$ and $\hat\beta$ is orthogonal to $X$.
\end{itemize}}
\ecor
\bpf 
The bound on the number of regions comes from Corollary \ref{pair-punc}. The required vector field comes from the construction in Lemma \ref{pucture-strat} and Corollary \ref{pair-punc}. \epf

We now adapt Canary's interpolation Theorem \cite{Ca1996} and some of Brock's results \cite{Bro01} to our setting. To do so we will replace the technology of simplicial hyperbolic surfaces with \emph{collapsed simplicial ruled surfaces} (CSRS) developed in \cite[Section 5]{BS2017}. We will need to use a mild generalisations of CSRS, specifically we will need to adapt it to surfaces $f:S_{g,k}\rar M$ with $k$-boundary components that will correspond to funnels in $M$. Moreover, in our special setting we will have that geodesics $\Upsilon$ corresponding to $\partial S_{g,k}$ will be on the boundary of the convex core $\partial CC(M)$.

\bdefi
Given a hyperbolic 3-manifold $M\cong S_{g,n}\times\mathbb R$ we say that a proper $\pi_1$-injective map $f:S_{k,m}\rar M$ is \emph{type preserving} if it maps cusps of $S_{k,m}$ to parabolic cusps or geometrically finite funnels of $M$.
\edefi

and:
\bdefi\label{SRS} Let $f:S_{g,n}\rar M$ be a proper map and $\Gamma\subset\pi_1(S_{g,n})$ be the peripheral loops of $S_{g,n}$ corresponding to funnels. Let $\tau$ be a triangulation of $S_{g,n}$ such that $\Gamma\subset \tau^1$. The data $(S,f,\tau,\Gamma)$ is a \emph{pre-simplicial ruled surface} (SRS) if for every triangle $\triangle\in \tau$ we have that:
\begin{itemize}
\item[(1)] every edge of $\triangle$ is mapped to a geodesic segment or ray in $M$;
\item[(2)] there is a foliation of $\triangle$ by segments beginning at some vertex $v\in\triangle$ and terminating on the opposite edge that are mapped to geodesic segments or rays in $M$;
\item[(3)] every funnel is triangulated as a square (with an edge at infinity) with a diagonal in which two opposite sides are identified.
\end{itemize}
When the extrinsic cone angle at every vertex $\alpha\in\tau^0$ is at least $2\pi$ we say that $(S,f,\tau,\Gamma)$ is a \emph{simplicial ruled surface}. \edefi

Following \cite{BS2017}, we define a weight system on $(S_{g,k},\tau)$ by: $w: E(\tau)\rar\R_{\geq0}$ with the property that on any edge $e$ of a funnel $w(e)=1$. We, say that a map $f:S\rar M$ minimises $(\tau,w)$ if it minimises: 
$$ L(f,\tau,w)\eqdef \sum_{e\in E(\tau)} w(e)(\ell_N(f(e)))^2$$

Since on funnels any two properly homotopic minimisers $f_0,f_1$ for $(\tau,w)$ are equivalent we have that the proof of Fact 5.4 in \cite{BS2017} follows verbatim. Thus, we have:

\blem\cite[Fact 5.4]{BS2017}\label{fact5.4} Suppose $f_i:S_{g,n}\rar M$, $i=0,1$, are properly homotopic, both minimisers for $(\tau,w)$ and $(f_i)_*(\pi_1(S_{g,n}))$ is not cyclic. Then, $f_0=f_1$ on $\tau$ up to orientation-preserving reparameterisation of each edge.
\elem

 Finally, we say that $(S',\tau')$ is a collapse of $(S,\tau)$ if we have a map $c:S\rar S'$ such that for any triangle $\triangle\in \tau$ the restriction $c\vert_\triangle$ is one of the following:
\begin{itemize}
\item an affine isomorphism on a triangle $\triangle'\in\tau'$;
\item an affine map onto an edge $e'\in\tau'$ that is an isomorphism on two edges of $\triangle$ and is mapping the third into $\partial e'$;
\item a constant map.
\end{itemize}
Such a surface will be referred as a \emph{collapsed simplicial ruled surface} (CSRS).

We can then show, following \cite[5.5]{BS2017}, the existence of minimisers in our setting.
\blem\label{simplhypsurf}
Consider $M \in GF( S_{g,n}\times I)$, possibly with funnels, and let $g:S_{g,n}\rar M$ be a proper type-preserving homotopy equivalence where $f(\partial S_{g,n})$ \rev{is comprised of loxodromics whose geodesic representatives live on $\partial CC(M)$ and correspond to funnels. Then, for any triangulation $\tau$ of $S_{g,n}$ with the $n$ peripheral loops $\Gamma\subset \tau^1$  and any weight $w$ on $\tau$ there is a collapse: $c:(S,\tau)\rar(S',\tau')$ and $g:S'\rar M$ such that $g\circ c$ is properly homotopic to $f$ and minimizes $L(f,\tau,w)$.}\elem

\bpf First note that by the definition of $\tau$ and $w$ by tightening $f$ around $\Gamma$ we have a proper homotopy of $f$ so that for any triangle $\triangle\in \tau$ contained in the funnel we have that $f\vert_\triangle$ satisfies properties (1)-(3) of Definition \ref{SRS}. Then, since our minimizers are fixed in the funnels by following the proof of the corresponding Lemma \cite[5.5]{BS2017} one obtains the required result. \rev{Moreover, note that all collapses necessarily happen outside funnels since $f\vert_\triangle$ is fixed in the funnels.}\epf

\brem\label{cont}
Similarly to Lemma \cite[5.6]{BS2017} we have that after picking a `ruling vertex' $rv(\triangle)$ in every $\triangle$ in $\tau$ the map $g\circ c$ in Lemma \ref{simplhypsurf} is unique if we also require that:
\begin{itemize}
\item $g\circ c\vert_e$, for $e$ an edge in $\tau$, is a constant speed geodesic; 
\item for any triangle $\triangle$ in $\tau$ the restriction $g\circ c\vert_\triangle$ is a constant speed geodesic on every affine line $rv(\triangle)$ to the opposite side.
\end{itemize}
Moreover, the map: $w\mapsto g\circ c$ is \rev{continuous if we topologise the target space with uniform convergence}. We will now adapt the interpolation result to our setting.\erem

\bthm[Interpolation, \cite{BS2017}]\label{interpolation} Let $M \in GF(S_{g,n}\times I, \alpha \times\set 0\cup\beta\times\set1)$ such that $\Gamma=\partial S_{g,n}$ are geodesics realised on $\partial CC(M)$. Let $X_\alpha\cong S_{g,n}\setminus\alpha$ and $X_\beta\cong S_{g,n}\setminus\beta$ be the conformal boundaries of $M$ separated along the geodesic lifts of $\Gamma$. Moreover, let $P_\alpha $ and $P_\beta $ be Bers pant's decomposition of $X_\alpha$ and $X_\beta$ respectively. Then, there exists CSRS $f_\alpha,f_\beta:S_{g,n}\rar M$ realising $P_\alpha$ and $P_\beta$ respectively and a continuous family of CSRS $f_t:S_{g,n}\rar M$ interpolating between them.
\ethm 
\bpf Let $\mathcal T_0$, $\mathcal T_1$ be triangulations of $S_{g,n}$ such that $\Gamma\cup P_\alpha$ and $\Gamma\cup P_\beta$ are contained in the $1$-skeleton of $\mathcal T_0$ and $\mathcal T_1$ respectively and let $w_0,w_1$ be weights on $E(\mathcal T_0)$, $E(\mathcal T_1)$ taking value $1$ on $\Gamma\cup P_\alpha$ and $\Gamma\cup P_\beta$ respectively. Then, by Lemma \ref{simplhypsurf} with $g:S_{g,n}\hookrightarrow M$ the natural inclusion we have CSRS realising $f_0,f_1: S_{g,n}\rar M$ minimising $L(g,\mathcal T_0,w_0)$ and $L(g,\mathcal T_1,w_1)$ respectively. Moreover, by the fact that $\Gamma\subset \partial CC(M)$ we have that $f_0,f_1$ are totally geodesic embeddings in $M\setminus CC(M)$ with the same image. To complete the Theorem we need to build the continuous family of CSRS interpolating between $f_0$ and $f_1$. We will achieve so by applying the interpolation result of \cite[Section 5]{BS2017}.

\begin{figure}[h]
\begin{overpic}[scale=0.5] {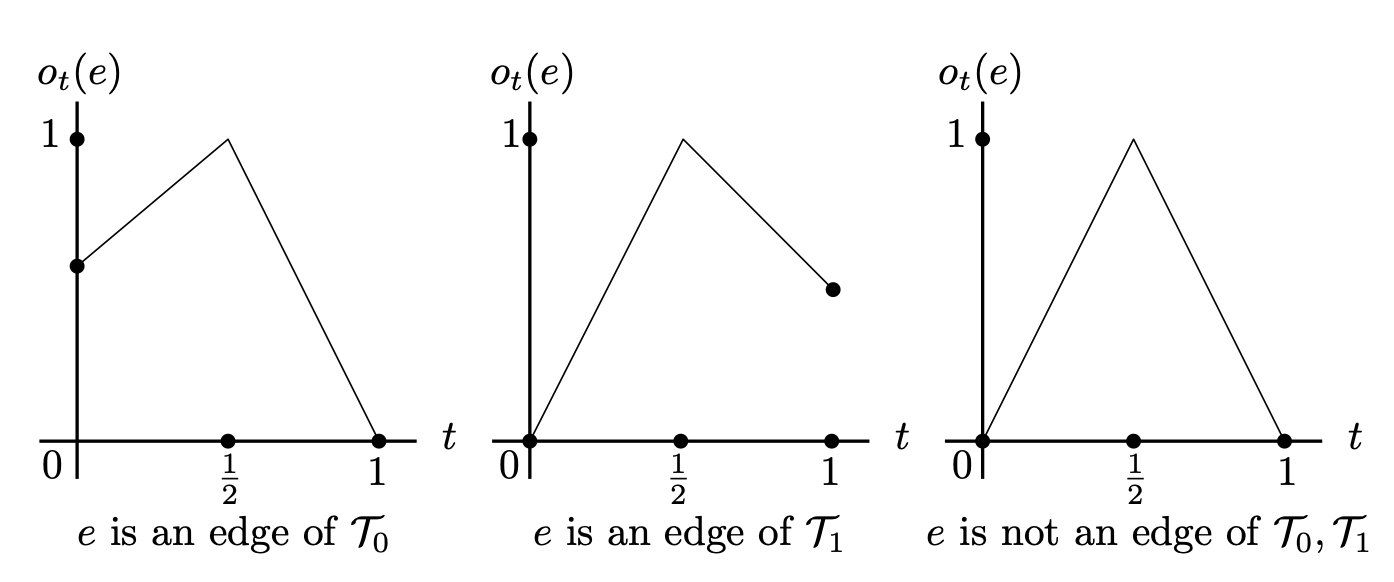}
\end{overpic}
\caption{Picture taken from \cite[Figure 9]{BS2017}.}\label{weights}
\end{figure}

Let $F\subset S_{g,n}$ be the funnels in which we can assume, without loss of generality, $\mathcal T_0\vert_F=\mathcal T_1\vert_F$. By an isotopy of $\mathcal T_0$, $\mathcal T_1$ in $S_{g,n}\setminus F$ we can assume that the 1-skeleta are transverse. Then, let $\mathcal T$ be a common refinement of $\mathcal T_0$ and $\mathcal T_1$ such that $\mathcal T\vert_F=\mathcal T_0\vert_F=\mathcal T_1\vert_F$. Then, define $\mathcal T_0',\mathcal T_1'\subset \mathcal T^1$ to be the edges induced by $\mathcal T_0$ and $\mathcal T_1$ respectively. The weights $w_0$ and $w_1$ induce weights $o_0'$ on $\mathcal T_0'$ and $o_1'$ on $\mathcal T_1'$ by adding the convention that if, say, we subdivide an edge $e\in\mathcal T_0$ into $n$ edges $e_j$'s the corresponding weight on $e_i\subset \mathcal T_0'$ is $o_0'(e_j)=\frac{w_0(e)}n$, similarly for $\mathcal T_1'$. Moreover, up to an isotopy of $f_0,f_1$, we can also assume that $\ell_M(f_\alpha(e_j))=\ell_M(f_\alpha(e_i))$, $i,j\leq n$. Then, we define the weights $o_t$ on $\tau$ as in Figure \ref{weights}.

Fix an affine structure on $\mathcal T$ such that $f_0$ and $f_1$ have constant speed on edges of $\mathcal T_0'$ and $\mathcal T_1'$ respectively. By arbitrarily picking ruling vertices and applying Lemma \ref{simplhypsurf} and Remark \ref{cont} we obtain for $t\in(0,1)$, which is when $o_t>0$ on every edge of $\mathcal T$, a continuous family $F_t:S_{g,n}\rar M$ of CSRS. Moreover, since $\iota_*(\pi_1( \tau_\alpha'))$ is not parabolic (similarly for $\tau_\beta')$, by Lemma \cite[5.7]{BS2017} we have that:
$$\lim_{t\rar0} F_t\vert_{\mathcal T_0'}=f_0 \vert_{\mathcal T_0'} \qquad \lim_{t\rar 1} F_t\vert_{\mathcal T_1'}=f_1\vert_{\mathcal T_0'}.$$
Moreover, we have that all $F_t$ are the same on the funnels. Thus, by picking some very small $t_0>0$ we have a homotopy, with very small trace, from $F_{t_0}$ to a map $F'_{t_0}$ such that $F'_{t_0}$ agrees with $f_0$ on $\mathcal T_0'$ and in the funnels. Then, by working on the compact surfaces we conclude as in the Interpolation Theorem \cite{BS2017}. \epf

Finally, the proof of Lemma \cite[5.2]{BS2017} carries through \rev{in our setting} so for any SRS $f:S_{g,n}\rar M$ let $S$ be \rev{the surface with the pull-back metric induced from $f$} then, there is a hyperbolic surface $X_{hyp}$ and a $K$-Lipschitz homeomorphism $X_{hyp}\rar S$ where $K=K(\abs{\chi(S_{g,n})})$. Thus, given an interpolation between SRS or CSRS we obtain a path in the corresponding Teichm\"uller space $\mathcal T(S_{g,n})$.

Given a 3-manifold with non-empty boundary $(M,\partial M)$ we will use $D_\partial (M)$ to denote the manifold obtained by doubling $M$ alongs its boundary: $D_\partial (M)\eqdef  M\cup_{\partial M} M$.

\bprop \label{doublevolumebound}Let $S_{g,n}$ be a hyperbolic surface and $\Gamma=\alpha\cup\beta$ be a pair of filling simple closed curves in minimal position. Let $\bar\Gamma$ be stratified lifts of $\Gamma $ in the interior of $S_{g,n}\times I$. Then, there exists constants $K_1\geq 1$ and $K_0\geq 0$ depending only on $S_{g,n}$ such that: 
$$\frac 1 {K_1}  \inf_{P_\alpha,P_\beta}d_{\mathcal P(S_{g,n})}(P_{\alpha},P_{\beta}) -K_0\leq vol(D_\partial (S_{g,n}\times I)\setminus \bar \Gamma)$$
where $P_\alpha$ and $P_\beta$ are any pants decompositions containing $\alpha,\beta$ respectively and $\bar\Gamma$ is deleted in only one copy of $S_{g,n}\times I$ in $D_\sigma(S_{g,n}\times I)$.\eprop
\bpf Let $S_i\cong S_{g,n}$, $i=1, 2$, be the push-off of the stratifying surfaces so that $S_i$ separates $\bar \alpha$ from $\bar\beta$.

In $D_\partial (S_{g,n}\times I)\setminus \bar \Gamma$ we have a properly embedded separating surface $S$ homeomorphic to $D_\partial (S_{g,n})\setminus \bar\Gamma$. Note that because $\alpha\cup\beta$ are filling then $S$ is incompressible in $D_\partial (S_{g,n}\times I)\setminus \bar \Gamma$.

\vspace{0.3cm}

\textbf{Claim 1:} The surface $S$ is totally geodesic.

\vspace{0.3cm}

\bpfc The surface $S$ splits $D_\partial (S_{g,n}\times I)\setminus \bar \Gamma$ into two components both homeomorphic to $S_{g,n}\times I\setminus(\alpha\times\set 0\cup\beta\times\set 1)$. Thus, we have an involution: $\iota: D_\partial (S_{g,n}\times I)\setminus \bar \Gamma\rar D_\partial (S_{g,n}\times I)\setminus \bar \Gamma$ exchanging the two components whose fixed set is $S$. \epfc

We define by $C_1$ and $C_2$ the totally geodesic submanifolds of $D_\partial (S_{g,n}\times I)\setminus \bar \Gamma$ containing the stratifying surfaces $S_1$ and $S_2$ respectively. We also note that the inclusions $S_i\hookrightarrow C_i$ induce homotopy equivalences.

For $i=1,2$ the cover $M_i$ of $D_\partial (S_{g,n}\times I)\setminus \bar \Gamma$ corresponding to $\pi_1(S_i)$ is homeomorphic to $S_{g,n}\times\R$ and the induced hyperbolic structure is in $AH\left(S_{g,n}\times I,\alpha \times\set 0\cup \beta\times\set 1\right)$ where the end is ``geometrically finite" with funnels corresponding to $\partial S_{g,n}$. Moreover, since $\alpha$ and $\beta$ fill we have that $M_i$ has no accidental parabolics.
 
  We now want to count $L$-short curves in $CC(M_i)$. Note that since $C_i$ lifts homeomorphically to $M_i$ and the maps:
 $$S_i\hookrightarrow C_i\overset{\tilde\iota}\hookrightarrow M_i$$
induce homotopy equivalences. Since $C_i$ has totally geodesic boundary we get that $\tilde\iota(C_i)=CC(M_i)$.

 By the proof of Lemma \ref{geomlb} we have that:
$$\frac{\vert \mathcal G_L^i\vert}{C_1}  < vol(\mathcal N_{\epsilon}(M_i^{\geq \epsilon}))\qquad \frac{\vert \mathcal G_L\vert}{C_1}\leq vol(D_\partial (S_{g,n}\times I)\setminus \overline \Gamma)$$
for a constant $C_1>1$ depending on $L$ and $\mathcal G_L^i$ the set of homotopy classes of  loops in $M_i$ of length less than $L$ and:
$$\frac 1 2\left( \vert \mathcal G_L^2\vert+ \vert \mathcal G_L^1\vert\right)\leq \vert \mathcal G_L\vert$$
thus it suffices to bound $\vert\mathcal G_L^i\vert$, for $i=1,2$, in terms of pants distance.

\vspace{0.3cm}

\textbf{Claim 2:} There exists constants $K_1\geq 1, K_0\geq 0$ depending on $g$ such that for Bers pants decompositions. $P_\alpha$, $P_\beta$ containing $\alpha,\beta$ we have:
$$d_{\mathcal P(S_{g,n})}(P_\alpha,P_\beta)-K_0\leq K_1 \vert \mathcal G_L^i\vert$$

\vspace{0.3cm}

\bpfc Let $E\eqdef  \set{\eta_j}_{j=1}^n$ be the geodesic representatives of $\partial S_i$ in $M_j$. Since $\partial C_i$ is totally geodesic and each $\eta_j$ is homotopic, in $C_i$, into $\partial C_i$ we have that:
$$E\subset \partial C_i.$$
Thus, in the conformal boundary $\Omega_i$ corresponding to $M_i$ there are natural choices of arcs $\tilde \eta_i$ connecting points in the limit set such that $\tilde E\cup \Lambda_i$ splits $\hat {\C}$ into two disks $ \Delta_1$ and $\Delta_2$ that are invariant under the action of the Kleinian group $G_i$ such that $M_i\cong\hyp{G_i}$. Then, we obtain two conformal structures $X_i\eqdef  \quotient{\Delta_1}{G_i}$ and $Y_i\eqdef  \quotient{\Delta_2}{G_i}$, $i=1,2$, on $S_{g,n+2}$ such that they both have totally geodesic boundary, corresponding to $E$, and have cusps corresponding to $\alpha$ and $\beta$ respectively.

Let $L$ be the Bers constant, \cite[12.8]{FM2011}, for $S_{g,n+2}$ and let $P_\alpha, P_\beta$ be Bers pants decompositions of $S$. Note that they both contain $E$. We will now count $L$-short loops in $M_i$.

By Lemma \ref{simplhypsurf} we have CSRS $Z_\alpha\subset M_i, Z_\beta\subset M_i$ adapted to $P_{\alpha}\setminus \alpha$ and $P_\beta\setminus\beta$. Let $\tau_\alpha$ and $\tau_\beta$ be their associated triangulations and recall that by a result of Bers for any loop $\gamma$ in $P_\alpha\cup P_\beta\setminus \set{\alpha\cup\beta}$ we have that:
$$\ell_{M_i}(\gamma)\leq 2L$$

Then, by the Interpolation Theorem \ref{interpolation} we have a continuous path:
$$h_t: S_{g,n}\rar M_i$$
interpolating between $Z_\alpha$ and $Z_\beta$. We denote by $Z_t^h$ the corresponding structures in $\mathcal T(S_{g,n},E)$. Thus, we can find pants decompositions $P_1,\dotsc, P_N$ such that we have an open cover $U_1,\dotsc, U_n$ of $Z_t$ with:
$$U_i\subset V(P_i)$$
where $V(P)\subset \mathcal T(S_{g,n},E)$ is the collection of hyperbolic structures in which the loops corresponding to $P$ have length at most $L$. Thus, by applying Lemmas \cite[Lemmas 3.3 and 4.3]{Bro01} we get a constant $K_0$, depending only on $S_{g,n}$, such that:
$$d_{\mathcal P(S_{g,n})}(P_1,P_N)\leq K_0\abs{\mathcal S_g}$$
for $\mathcal S_g\subset \mathcal G_L^i$ a sub-set of $L$-short loops. Thus, by Lemma \cite[4.7]{Bro01} we have:
$$d_{\mathcal P(S_{g,n})}(P_\alpha,P_\beta)-2B'\leq d_{\mathcal P(S_{g,n})}(P_1,P_N).$$
Combining the two equations we obtain the desired result. \epfc

Thus, we have that:

$$\frac 1 {2K_1C_1}d_{\mathcal P(S_{g,n})}(P_\alpha,P_\beta) -\frac{K_0}{C_1K_1}\leq \frac{\abs{\mathcal G_L}}{C_1}\leq vol(D_\partial (S_{g,n}\times I)\setminus \overline \Gamma)$$

for $P_\alpha,P_\beta$ the Bers pants decompositions containing $\alpha,\beta$. By taking the \rev{infimum} the required result follows.\epf

We can finally prove the lower bound in our main Theorem:

\bthm\label{pairoffillingvol} Let $S=S_g$ be a hyperbolic surface, $N = PT(S)$, and let $(\al, \beta)$ be a filling pair of multicurves on $S$ in minimal position. Then there are constants $K_1 \geq 1$, $K_0\geq 0$ depending only on $S$, such that:
\[  \frac1 {K_1} \inf_{P_{\al},P_{\beta}} d_{\mathcal P(S)}(P_{\al},P_{\beta})-K_0\leq vol(N_{(\hat{\al}, \hat{\beta})}) ,\]
where $P_{\al}, P_{\beta}$ are any pants decompositions of $S$ with $\al \subset P_{\al}, \beta \subset P_{\beta}$. Further, when $S$ has punctures, this bound is asymptotically optimal.
\ethm
\bpf  

Let $D_1,\dotsc, D_n$, $n\leq N=N(g)$, be the non-square complementary region, let $X\in\mathfrak X(S)$ as in Corollary \ref{vector field} and let $x_i\in D_i$ be the singular points. The lower bound will come by cutting $PT(S)\setminus (\hat\alpha\cup\hat\beta)$ along an incompressible surfaces and then applying a Theorem of Agol-Storm-Thurston \cite{agol2005lower} to get a lower bound on the volume by bounding one of the components.

Let $F_i$ be the fibers in $PT(S)$ corresponding to $x_i$ and note that $PT(S_g)\setminus \cup_{i=1}^n F_i=PT(S_{g,n})$. Thus, we have that:
$$PT(S_g)=PT(S_{g,n})\cup\left(\cup_{i=1}^n V_i\right)$$
for $V_i$ solid tori. Moreover, by Lemma \ref{pucture-strat} we have that $\hat\alpha,\hat\beta$ are stratified by sections in $PT(S_{g,n})$. Let $\sigma:S_{g,n}\hookrightarrow PT(S_{g,n})$ be the section induced by the vector field $X$ coming from Corollary \ref{vector field}. By using the meridian $m_i$ coming from $V_i$ and $F_i$ as longitude we get coordinates on each boundary component of $PT(S_{g,n})$.

\textbf{Claim 1:} Let $\partial x_i$ be the boundary component of $PT(S_{g,n})$ corresponding to $x_i$. Then, $\sigma\cap \partial x_i$ is a $(1,p_i)$ curve in which $\abs{p_i}\geq 1$ with equality if and only if $x\in D_i$ is an hexagon.

\bpfc Since $\sigma$ is induced by a vector field the slope is equal to the index of the vector field and we are assuming that we have no hexagon the result follows from Corollary \ref{vector field}. \epfc

Let $N\eqdef  N_\epsilon (\sigma)\overset{\psi}{\cong} S_{g,n}\times I$, $A\eqdef  \psi^{-1}(\partial S_{g,n}\times I)$ let $M$ be its complement in $PT(S)$ so that $\hat\alpha,\hat\beta\subset M$. Then, we get the following decomposition of $PT(S_g)$:
$$PT(S_g)=\left(N \cup _A \cup_{i=1}^n V_i\right)\cup_\Sigma M$$
where $\Sigma=\partial M\cong S_{2g+n-1}$.

If we have hexagonal regions we have that the foliation $\xi$ induced by the vector field $X$ on $S_g$ is not orientable. If we pass to the orientation double cover $UT(S_g)\rar PT(S_g)$\rev{, where $UT(S_g)$ the unit tangent bundle,} we get that the lift $\tilde\xi$ of $\xi$ is a connected immersed surface in $UT(S_g)$ such that the projection: $p:UT(S_g)\rar S_g$ induces a degree two branched cover $p:\tilde\xi\rar S_g$ with branch points given by $\sing(X)=\sing(\xi)=\set{x_i}_{i=1}^n$.

Let $\hat \al^\pm$ and $\hat \beta^\pm$ be the canonical lifts of $\al$ and $\beta$ in $UT(S_g)$, with two lifts for each curve corresponding to the two orientations. Notice that $UT(S_g) \setminus (\hat \al^\pm \cup  \hat \beta^\pm)$ double covers  $PT(S_g)  \setminus (\hat \al \cup  \hat \beta)$ and that both manifolds are hyperbolic. Thus, since $UT(S_g) \setminus (\hat \al^\pm \cup  \hat \beta^\pm)$ is twice the volume of $PT(S_g)  \setminus (\hat \al \cup  \hat \beta)$ to obtain a lower bound for $PT(S_g)  \setminus (\hat \al \cup  \hat \beta)$ it suffices to find one for $UT(S_g) \setminus (\hat \al^\pm \cup  \hat \beta^\pm)$. Moreover, we get that $\tilde\xi$ is an orientable foliation hence $\alpha^\pm\cup\beta^\pm$ have no hexagonal complementary regions.

Let $F_1,\dotsc, F_n$ be the pre-images under $UT(S_g)\rar PT(S_g)\rar S_g$ of $x_1,\dotsc, x_n$. By deleting thickened neighborhoods $U_i$ of the $F_i$'s in $UT(S_g)$ we get that $UT(S_g)\setminus \left( \cup_{i=1}^n U_i\right)$ is homeomorphic to $UT(S_{g,n})$. Moreover, in $UT(S_{g,n})$ the surface $\tilde\xi\cong S_{2g-1,n}$, since the double cover of $S_g$ has genus $2g-1$, is properly embedded, carries the loops $\hat \al^\pm$ and $\hat \beta^\pm$ and induces a homeomorphism:
$$UT(S_{g,n})\cong S_{2g-1,n}\times I/(x,0)\sim (\phi(x),0),\qquad \phi\in \mathrm{Mod}(S_{2g-1,n})$$
Let $\sigma$ be a fiber of $UT(S_{g,n})$ disjoint from $\hat \al^\pm\cup\hat \beta^\pm$ and let $N\eqdef  N_\epsilon(\sigma)$ and $M$ denote its complement in $UT(S_{g,n})$. Then, we get the following decomposition:
$$UT(S_g)=UT(S_{g,n})\cup \left( \cup_{i=1}^n U_i\right)=\left(N\cup \left( \cup_{i=1}^n U_i\right)\right)\cup M$$
as before we let $\Sigma\eqdef  \partial M$. Note that $\Sigma=D_\partial(S_{2g-1,n})\cong S_{4g+n-3}$. By \textbf{Claim 1}, by picking the meridian $m_i$ and the longitude $F_i$ as coordinates for $\partial U_i$ we get that each slope on $\partial U_i$ is $(1,p_i)$ with $\abs{p_i}>1$.

With an abuse of notation we will now use $X$ for either $PT(S_g)$ or $UT(S_g)$ and use $\hat\Gamma$ for $\hat\alpha\cup\hat\beta$ or $\hat\alpha^\pm\cup\hat\beta^\pm$ respectively. Thus, in either case we get an embedded closed surface $\Sigma\hookrightarrow X$ splitting $X$ into two submanifolds $M\cong F\times I\setminus \hat\Gamma$ and $N\cup_A\cup_{i=1}^n V_i$, for $V_i$ solid tori and $N\cong F\times I$ where $F\cong S_{g,n}$ or $S_{2g-1,n}$ depending on whether $X\cong PT(S_g)$ or $UT(S_g)$.

\textbf{Claim 2:} The surface $\Sigma$ is incompressible in $X\setminus\hat\Gamma$.

\bpfc Since $\hat\Gamma$ forms a filling system it is obvious that $\Sigma=\partial M$ is incompressible in $M\setminus\hat\Gamma$. Moreover, since by \textbf{Claim 1} we get that each component of $A$ is attached to a neighborhood $(1,p_i)$, $\abs{p_i}>1$, slope in $V_i$ by Lemma \ref{incompsurfacegluing} we get that $\Sigma$ is also incompressible in $N \cup _A \cup_{i=1}^n V_i$ and the result follows. \epfc

By a theorem of Agol-Storm-Thurston \cite{agol2005lower} we get that:
$$\frac{v_3}{2}||D_{\Sigma}(M\setminus\hat\Gamma)||\leq vol (X\setminus\hat\Gamma)$$
where $||M||$ is the simplicial volumes of $M$. The manifold $D_{\Sigma}(M\setminus \hat\Gamma )$ is not hyperbolic as it is not acylindrical. However, the acylindrical piece of \rev{$D_{\Sigma}(M\setminus \hat\Gamma )$} is homeomorphic to:
$$D_\partial (F\times I)\setminus (\hat\alpha\cup\hat\beta)\qquad \text{ where} \qquad F\cong S_{g,n} \text{ or } S_{2g-1,n}$$
or alternatively: $D_{\Sigma}^{hyp}(M\setminus \hat\Gamma)$ is equal to $D_{\Sigma}(M\setminus \hat\Gamma)$ in which we do trivial Dehn filling on one of the copies of $\hat\Gamma$. Finally, we are in the setting of Proposition \ref{doublevolumebound} and we get constant $K_1\geq 1$, $K\geq 0$ depending only on $g$, since $n\leq N(g)$, such that:
$$ \frac 1 {K_1} \inf_{P_\alpha,P_\beta}d_{\mathcal P(F)}(P_{\alpha},P_{\beta})-K_0 \leq vol(D_{\Sigma}^{hyp}(M\setminus \hat\Gamma)=||D_{\Sigma}(M\setminus \hat\Gamma)||$$
for $P_\alpha$, $P_\beta$ any pants decomposition containing $\alpha$ and $\beta$ respectively. Combining this with the previous inequality we get:
$$\frac 1 {K_1} \inf_{P_\alpha,P_\beta}d_{\mathcal P(F)}(P_{\alpha},P_{\beta})-K_0 \leq vol(X\setminus\hat\Gamma)$$

By Theorem \cite[1.13]{AWT2017} we get a constant $C=C(S_g)$ such that $\mathcal P(S_{2g-1})$ is quasi-isometric to $\mathcal P(S_g)$ and so by Lemma \ref{pants} we have that:
$$\frac 1 {2 C K_1} \inf_{P_\alpha,P_\beta}d_{\mathcal P(S_{g})}(P_{\alpha},P_{\beta})-\frac{C}{K_1}\leq \frac 1 {K_1} \inf_{P_\alpha,P_\beta}d_{\mathcal P(F)}(P_{\alpha},P_{\beta})$$
completing the proof. Note that if $F=S_{g,n}$ then $C=\frac 1 2$ and the $-\frac C{K_1}$ term disappears.\epf
Theorem \ref{pairoffillingvol} concludes the first part of Theorem \ref{pair}. The quasi-isometry with WP-geometry follows from the following discussion.

Since the Weil-Petersson metric $d_{WP}$ gives us a metric on $\overline{\mathcal T(S)}$ in which different strata are at finite distance we have that:

$$\alpha,\beta\in \mathfrak M(S):\quad d_{WP}(S(\alpha),S(\beta))=\inf_{\substack{X\in S(\alpha)\\ Y\in S(\beta)}} d_{WP}(X,Y)$$

By \cite{Bro01} we have:

\bthm There exists a map $q: \mathcal P(S)\rar \mathcal T(S)$ and constants $K_1\geq 1$, $K_0\geq 0$ and $K_3\geq 0$ such that $q$ satisfies:
$$\frac{ d_{\mathcal P(S)}(P_1,P_2)}{K_1}-K_2\leq d_{WP}(q(P_1),q(P_2))\leq K_1 d_{\mathcal P(S)}(P_1,P_2)+K_2$$
and the map $q$ is $dK_3$-surjective. Moreover, all the constants only depend on $S$.\ethm
The map $q$ is defined so that $q(P)$ is any point $X\in V(P)$. Thus, one can think of the quasi-inverse map $Q:\mathcal T(S)\rar \mathcal P(S)$ mapping $X$ to $P$ such that a $X\in V(P)$. Hence, one obtains the following immediate corollary:
\bcor\label{wpcorollary}
There is q.i.-map, whose constants only depend on $S$, such that for any pair of simple multi-curves $\alpha,\beta\in\mathfrak (S)$:
$$ d_{WP}(S(\alpha),S(\beta))=\inf_{\substack{X\in S(\alpha)\\ Y\in S(\beta)}} d_{WP}(X,Y)\overset{q.i.}{\simeq} \inf_{\substack{\alpha\in P_\alpha\\ \beta\in P_\beta}} d_{\mathcal P(S)}(P_\alpha,P_\beta) $$
\ecor

\section{Applications to $PT(S)$}\label{application}

\subsection{Canonical lift complements for $\Sigma_{1,1}$ and $\Sigma_{0,4}$}
 
 In this section, we study canonical lift complements of filling collections of essential simple closed curves in $S = \Sigma_{1,1}$ or $\Sigma_{0,4}$. Before proceeding, let us remark on the following useful correspondence between essential simple closed curves in the surfaces $\Sigma_{1,0}, \Sigma_{1,1}$ and $\Sigma_{0,4}$.
 
\blem{\cite[Proposition 2.6]{FM2011}}\label{sim} The inclusion map of $\Sigma_{1,1}$ in $\Sigma_{1,0},$ and the hyperelliptic involution from $\Sigma_{1,0}$ to $\Sigma_{0,4}$, induce bijections between the sets of basepoint-free homotopy classes of essential simple closed curves of $\Sigma_{1,0}$, $\Sigma_{1,1}$, and $\Sigma_{0,4}$. Further, this extends to transversal homotopy classes of collections of essential simple closed curves in minimal position.
 \elem

 We will also make use of the following fact, which follows from all the graphs in question being isomorphic to the Farey graph, see \cite{Schleimer}.
\blem If $S = \Sigma_{1,1}$ or $\Sigma_{0,4}$, then there is a natural isometry between $\mathcal P(S)$ and $\mathcal C(S)$.
 \elem
 
Our goal for this section is to prove:
   
\begin{customthm}{B}\label{punctorus} Let $S = \Sigma_{1,1}$ or $\Sigma_{0,4}$, $N = PT(S)$, and let $\Gamma$ be a filling collection of non-parallel essential simple closed curves in minimal position. Then, there exists constants $K_1>1$ and $K_0>0$, depending only on $S$, and an ordering $\Gamma = \{\gamma_i\}_{i = 1}^n$, such that:
$$  \frac{1}{2K_1}\left(\sum_{i =1}^n d_{\mathcal C(S)}(\gamma_i,\gamma_{i+1})\right) \leq vol(N_{\widehat{\Gamma}})\leq K_1 \left(\sum_{i = 1}^n d_{\mathcal C(S)}(\gamma_i,\gamma_{i+1})\right)+nK_0,$$
where $\gamma_{n+1} = \gamma_1$.
\end{customthm}

\bpf
Our main tool will be to build stratifying surfaces in $PT(\Sigma_{1,0})$ using flat structures and then push them to $PT(S)$. Note that since $\Gamma$ is filling and in minimal position, $N_{\widehat{\Gamma}}$ is hyperbolic by \cite{FH13}.

\vskip .2cm

 \textbf{Claim.} Given a filling collection $\Gamma_0$ of non-parallel essential simple closed curves in minimal position on $\Sigma_{1,0}$, there is a stratification of $\hat\Gamma_0$ in $PT(\Sigma_{1,0}).$

\bpfc Let $T$ be the flat square torus structure on $\Sigma_{1,0}$. Since any two minimal position representatives of $\Gamma_0$ give rise to isotopic canonical lifts by Corollary \ref{cor:stratwell}, we can assume that $\Gamma_0 = \{\gamma_i\}_{i = 1}^n$, where the $\gamma_i$ are flat geodesics ordered by increasing slope. For each $\gamma_i$, there is a line field on $T$ with the same slope as $\gamma_i$. These line fields can be viewed as disjoint and properly embedded incompressible surfaces in $PT(\Sigma_{1,0})$ that stratify $\hat\Gamma_0$.
\epfc

Let $\Gamma$ be a filling collection of non-parallel simple closed curves in minimal position on $S$.

\textbf{Case} $S=\Sigma_{1,1}.$  Let $\Gamma_0$ be the image of $\Gamma$ under the inclusion. By Lemma \ref{sim}, a transversal homotopy of $\Gamma_0$ in $\Sigma_{1,0}$ arrises from a transversal homotopy of $\Gamma$ in $\Sigma_{1,1}$. Thus, the transversal homotopy from $\Gamma_0$ to its geodesic representative on $T$ can be realized relative to the puncture of $\Sigma_{1,1}$. Remove the fiber over the puncture from $PT(\Sigma_{1,0})$ and puncture the stratifying surfaces of $\hat\Gamma_0$ to stratify $\hat\Gamma$.

\textbf{Case} $S=\Sigma_{0,4}.$ Let $T$ be the flat unit square torus structure on $\Sigma_{1,0}$ and consider the hyperelliptic involution obtained by 180\textdegree rotation. Remove the four fixed points of this involution and the corresponding fibers from $PT(\Sigma_{1,0})$. This gives the following commutative diagram
$$
\xymatrix{\hat\Gamma_1\subset PT(\Sigma_{1,4}) \ar[r]^{\hat{\mathfrak h}} \ar[d]^{} & PT(S) \ar[d]^{}\supseteq\hat\Gamma \\
\Gamma_1\subset\Sigma_{1,4} \ar[r] ^{\mathfrak h}& S\supseteq\Gamma
}$$
where $\hat{\mathfrak h}$ is the fiber-wise covering map induced by the hyperelliptic quotient map $\mathfrak h$. Let $\Gamma_1$ be a choice of one lift for each component of $\Gamma$ under $\mathfrak h$. Let $\Gamma_0$ on $\Sigma_{1,0}$ be obtained from $\Gamma_1$ by filling the punctures. By Lemma \ref{sim}, we can assume that the transversal homotopy taking $\Gamma_0$ to its geodesic representative on $T$ is relative to the four punctures. Thus, by puncturing the stratifying surfaces for $\hat\Gamma_0$, we get stratifying surfaces of $\hat\Gamma_1$. Since these stratifying surfaces are line fields on $T$, they cover disjoint and property embedded incompressible surfaces under $\hat{\mathfrak h}$, which are stratifying surfaces for $\hat\Gamma$.

In either case, by Theorem \ref{geombounds} and Remark \ref{pairwisefilling}, there are constants $K_1,K_0$, depending only on $S$, with
$$\frac{1}{2 K_1}\left(\sum_{i=1}^n d_{\mathcal P(S)}(P_{X_i},P_{Y_{i+1}})\right)\leq vol(N_{\hat\Gamma}) \leq K_1 \left(\sum_{i=1}^n d_{\mathcal P(S)}(P_{X_i},P_{Y_{i+1}})\right)+n K_0,$$
where $\{(P_{X_i},P_{Y_{i+1}})\}_{i = 1}^n$ are the Bers pants decompositions coming from the conformal structures of the covers corresponding to the stratifying surfaces. Since any such pants decomposition contains only one loop, the pants distance is the same as the curve distance. Moreover, this loop is forced to be the rank one cusp induced by $\gamma_i$. Therefore, we see that
$$d_{\mathcal P(S)}(P_{X_i},P_{Y_{i+1}})=d_{\mathcal C(S)}(\gamma_i,\gamma_{i+1})$$
and obtain our final result:
$$  \frac{1}{2 K_1}\left(\sum_{i = 1}^n d_{\mathcal C}(\gamma_i,\gamma_{i+1})\right)\leq vol(N_{\widehat{\Gamma}})\leq K_1 \left(\sum_{i = 1}^n d_{\mathcal C}(\gamma_i,\gamma_{i+1})\right)+nK_0.$$ \epf

\subsection{Examples of stratified canonical links on higher complexity surfaces}\label{examples}

In this section, we give various methods to constructing sequences of stratified canonical links. First, we build pairs $\Gamma_n$ of essential simple closed curves on hyperbolic surfaces $\Sigma_{g,k}$, with $(g,k)$ equal to $(0,m+4)$, $(1,m+1)$ or $(g,2+g\,m)$ for $m\in\N$, such that $\hat\Gamma_n$ is stratified and the volumes behave quasi-isometrically with respect to curve graph distance. This amounts to Theorem \ref{ce}, were we also show that the lower-bound of \cite{M17} behaves poorly for these links. In the second half of this section, we demonstrate how stratified canonical links naturally arise from quadratic differentials and pseudo-Anosov diffeomorphisms.

\subsubsection{Stratified filling pairs of simple geodesics in higher complexity surfaces}

Our goal for this subsection is to show

\begin{customthm}{C}
Let $S  = \Sigma_{g,k}$ be a hyperbolic surface with $(g,k)$ equal to $(0,m+4)$, $(1,m+1)$ or $(g,2+gm)$ for $m\in\N$ and let $N = PT(S)$. Then, there exists a sequence $\{(\alpha_n,\beta_n)\}_{n\in\mathbb N}$ of filling pairs of essential simple closed curves on $S$ in minimal position with the property that $d_{\mathcal C(S)}(\alpha_n,\beta_n)\nearrow \infty$  and a hyperbolic metric $X$ on $S$ such that:
$$vol(N_{(\hat{\alpha}_n,\, \hat\beta_n)})\asymp d_{\mathcal C(S)}(\alpha_n,\beta_n)\asymp \log\left(\ell_{X}(\alpha_n)+\ell_{X}(\beta_n)\right)$$
where the quasi-isometry constants only depend on $S$. Further, for any subsurface decomposition $\mathcal W$ of $S$, one has
$$\sum_{\Sigma\in \mathcal W}\sharp\{\text{homotopy classes of } \alpha_n \text{ and } \beta_n\text{-arcs in }\Sigma\} <  6(\kappa(S)+3)\kappa(S) \quad \text{for all} \quad n \in \mathbb N.$$
\end{customthm}

Note that the homotopy class count is exactly the best known lower bound given in \cite{M17}.

\bdefi\label{dec} 
Given a hyperbolic surface $S$, a \emph{subsurface decomposition of $S$} is the set of components obtained after cutting $S$ along a multicurve.
 \edefi

\blem\label{toro}
There exists a hyperbolic metric $X^o$ on $\Sigma_{1,1}$ and a sequence $\{(\alpha^o_n,\beta^o_n)\}_{n\in\mathbb N}$ of filling pairs of essential simple closed curves in minimal position with the property that $d_{\mathcal C(S)}(\alpha^o_n,\beta^o_n)\nearrow \infty$ and
$$vol(M_{(\hat\alpha^o_n, \hat\beta^o_n)})\asymp d_{\mathcal C(\Sigma_{1,1})}(\alpha^o_n,\beta^o_n)\asymp \log \left(\ell_{X^o}(\alpha^o_n)+\ell_{X^o}(\beta^o_n)\right),$$
\rev{where the constants in the coarse equivalences are uniform}.\elem

\bpf
By Theorem \ref{punctorus}, we only need to show the second quasi-isometry. Consider the hyperbolic metric $X^o$ on $\Sigma_{1,1}$ obtained by gluing the opposite sides of a regular ideal quadrilateral such that the orthogonal arcs become closed geodesics. We call these geodesics $a$ and $b$. It is easy to see that $a,b$ have the same length, intersect once at a right angle, and are the systoles of $X^o$. Treating their intersection point as a basepoint, we get $\pi_1(\Sigma_{1,1}) = \langle a, b\rangle$. We will use $+$ to denote composition in $\pi_1(\Sigma_{1,1})$.

Choose $\alpha^o_n$ to be $a$ and $\beta^o_n$ to be the geodesic representative of $u_n a + u_{n+1} b$, where $\{u_n\}_{n\in\mathbb N}$ is the Fibonacci sequence. By Lam\'e's Theorem \cite{Ho76}, we have that:
$$ d_{\mathcal C(\Sigma_{1,1})}(\alpha^o_n,\beta^o_n)=n-1$$
Let $s \eqdef \ell_{X^o}(a) = \ell_{X^o}(b) =\text{arccosh}\left(\frac{3+\sqrt{2}}2\right)$. Since $u_{n}\leq u_{n+1}\leq 2u_{n}$ and $\varphi^{n-2}\leq u_n\leq \varphi^{n-1},$ where $\varphi$ is the golden ratio, we conclude that:
$$ s \, \varphi^{n-1} \leq s\, u_{n+1}\leq \ell_{X^o}(\beta^o_n) \leq s\, (u_{n} + u_{n+1}) \leq3 s \, u_n \leq 3 s \, \varphi^{n-1}.$$
the second inequality comes from the fact that by cutting $\beta^o_n$ along $b$ gives us $u_{n+1}$ arcs each of which has length at least $s$. Therefore,
$$ \log(\varphi)d_{\mathcal C(\Sigma_{1,1})}(\alpha^o_n,\beta^o_n)-\log(4s) \leq \log (\ell_{X^o}(\alpha^o_n)+\ell_{X^o}(\beta^o_n))\leq \log(\varphi)d_{\mathcal C(\Sigma_{1,1})}(\alpha^o_n,\beta^o_n)+\log(4s)$$
which concludes the proof.\epf

We can now prove Theorem \ref{ce}.

\begin{proof}[Proof of Theorem \ref{ce}]

\begin{enumerate}[leftmargin=*]

\item[$(a)$]  \textbf{Case} $S=\Sigma_{1,k}$, $k>0$. Let $a,b, \alpha^o_n, \beta^o_n$, and $X^o$ be as in Lemma \ref{toro}. Notice that $a$ and $b$ intersect at a right angle. Consider, the $k$-fold cover $h : Y_k \rar X$ obtained by cutting along $b$ and gluing $k$-copies of $X \setminus b$ cyclically. Then, $\alpha_n \eqdef  h^{-1}(\alpha^o_n)$ is a simple closed geodesic that projects to $k$-times $\alpha^o_n = a$.

\paragraph{Claim:} Up to subsequence, the full lift $\beta_n \eqdef  h^{-1}(\beta^o_n)$ is connected.

\bpfc Since $u_n=\iota(\beta^o_n, b)$, the lift is connected whenever $gcd(u_n, k)=1$. Since $\set{u_n}_{n\in\N}$ is Fibonacci sequence, it is a simple fact that there is a subsequence $\{u_{n_i}\}_{i = 1}^\infty$ with $gcd(u_{n_i}, k)=1$ for all $i$.
\epfc

Thus, we can assume $\{(\alpha_n,\beta_n)\}_{n \in \mathbb N}$ is a sequence of filling pairs of essential simple closed curves on $\Sigma_{1,k}$.

\paragraph{Claim:} For all $n\in\N$ the canonical links $(\hat\alpha_n,\hat\beta_n)$ are stratified.

\bpfc Let $(S_n^1,S_n^2)$ be the sections stratifying $(\hat\alpha^o_n, \hat\beta^o_n)$ in Theorem \ref{punctorus}. The covering map $h :\Sigma_{1,k}\rar \Sigma_{1,1}$ induces a $k$-fold covering map $H : PT(\Sigma_{1,k})\rar PT(\Sigma_{1,1})$ mapping $\hat\alpha_n,\hat\beta_n$ to $\hat\alpha^o_n$, $\hat\beta^o_n$, respectively. Notice that the preimage of every section $\sigma:\Sigma_{1,1}\rar PT(\Sigma_{1,1})$ under $h$ is a section of $PT(\Sigma_{1,k})$.

The preimage $F_n^i \eqdef  H^{-1}(S_n^i)$ is a $\pi_1$-injective surface in $PT(\Sigma_{1,k})$ containing $\hat\alpha_n$ for $i = 1$ and  $\hat \beta_n$ for $i = 2$. Moreover, the $F_n^i$'s are essential, properly embedded, and $F_n^1 \cap F_n^2 = \emp$.  Finally, each $F_n^i$ is the image of a section, so homeomorphic to $\Sigma_{1,k}$.\epfc

Let $N \eqdef  PT(\Sigma_{1,k})$ and $N^o \eqdef  PT(\Sigma_{1,1})$. We get a $k$-fold cover
$N_{(\hat\alpha_n,\hat\beta_n)} \to N^o_{(\hat\alpha^o_n, \hat\beta^o_n)}$ and, therefore, $vol(N_{(\hat\alpha_n,\hat\beta_n)})=k\, vol(N^o_{(\hat\alpha^o_n, \hat\beta^o_n)})$. By Theorem \ref{punctorus}  and  Lemma \ref{toro},
$$vol(N^o_{(\hat\alpha^o_n, \hat\beta^o_n)})\asymp d_{\mathcal C(\Sigma_{1,1})}(\alpha^o_n,\beta^o_n) \asymp \log (\ell_{X^o}(\alpha^o_n)+\ell_{X^o}(\beta^o_n)).$$
Combining these equations we obtain:
$$vol(N_{(\hat\alpha_n,\hat\beta_n)})\asymp \log (\ell_{X^o}(\alpha^o_n)+\ell_{X^o}(\beta^o_n)) \asymp  \log (\ell_{X}(\alpha_n)+\ell_{X}(\beta_n)),$$
where $X$ is the pull-back of $X^o$ under the covering map. Note, the quasi-isometry constants only depend on $k$. Further, since covers induce quasi-isometric embeddings between curve graphes \cite{KS09}, we have that $d_{\mathcal C(\Sigma_{1,1})}(\alpha_n,\beta_n)\asymp d_{\mathcal C(\Sigma_{1,k})}(\alpha_n,\beta_n)$. Hence,
$$vol(N_{(\hat\alpha_n,\hat\beta_n)})\asymp d_{\mathcal C(\Sigma_{1,k})}(\alpha_n,\beta_n)\asymp \log (\ell_{X}(\alpha_n)+\ell_{X}(\beta_n))$$

\item[$(b)$] \textbf{Case} $S=\Sigma_{0,k}$ with $k\geq 4$. As in Theorem \ref{punctorus}, we can consider the hyperelliptic cover from $\Sigma_{1,4}$ to $\Sigma_{0,4}$ of degree $2,$  and induce a covering map from $\Sigma_{1,2(k-2)}$ to $\Sigma_{0,k}$ of degree $2,$ by removing $2(k-4)$ points from $\Sigma_{1,4}$ in a symmetric with respect to the cover.

By \cite{KS09}, covering maps between hyperbolic surfaces yield quasi-isometric embeddings between their corresponding curve graphes. Thus, we can take $\{(\alpha_n,\beta_n)\}_{n\in\mathbb N}$ on $\Sigma_{0,4}$ to be the image of the sequence from part $(a)$ for $\Sigma_{1,4},$ which is a sequence of filling pairs of simple closed curves, because the original covering from $\Sigma_{1,4}$ to $\Sigma_{0,4}$ sends simple closed curves to simple closed curves by Lemma \ref{sim}. Note that $(\alpha_n,\beta_n )$ is stratified by Theorem \ref{punctorus} for all $n \in \mathbb N$, so by part $(a)$ we have that:
$$vol(N_{(\hat\alpha_n,\hat\beta_n)})\asymp d_{\mathcal C(\Sigma_{0,k})}(\alpha_n,\beta_n)\asymp \log (\ell_{X}(\alpha_n)+\ell_{X}(\beta_n)),$$
where the hyperbolic metric $X$ is push-forward under the covering map of the metric given in $(a)$ for $\Sigma_{1,2(k-2)}$.

\item[$(c)$]  \textbf{Case} $S=\Sigma_{g,k}$ with $k=2+gm$, $m\in\N$, $g >1$.  First, we work with a twice-punctured surface of genus $g$ and then we explain how to generalise to $k=2+g \, m$, $m\in\N$. 

Let $\{(\alpha^1_n,\beta^1_n)\}_{n\in\mathbb N}$ be the sequence of filling closed curves on $\Sigma_{1,2}$ from part $(a)$ arising from the cover $h :\Sigma_{1,2}\to\Sigma_{1,1}$. Moreover, assume that the induced $h(\beta^1_n)= u_n a + u_{n+1}b\subset \Sigma_{1,1}$ have the property that $gcd(u_n, g)=1$. Consider the simple arc $\delta$ connecting the two punctures with the property that $h(\delta)$ on $\Sigma_{1,1}$ corresponds to the $a + b$ loop. Thus, $\iota(h(\beta^1_n),h(\delta))=u_n+u_{n+1}=u_{n+2}$ and $\iota(h(\alpha^1_n),h(\delta))=1$.

For $\Sigma_{g,2}$, consider a rotation through the two punctures in a symmetrically arranged way as in Figure \ref{co}. The quotient gives a degree $g$ cover of $\Sigma_{1,2}$. Let $\alpha_n$, $\beta_n$ be the preimages of $\alpha^1_n,\beta^1_n$ under this cover, respectively. Since
$$gcd(\iota(\beta^1_n,\delta),g)=gcd(\iota (h(\beta^1_n),h(\delta)),g)=1,$$
$\alpha_n$, $\beta_n$ are connected and give a filling pair of essential simple closed curves on $\Sigma_{g,2}$.

 \begin{figure}[h]
\centering
\includegraphics[scale=.3] {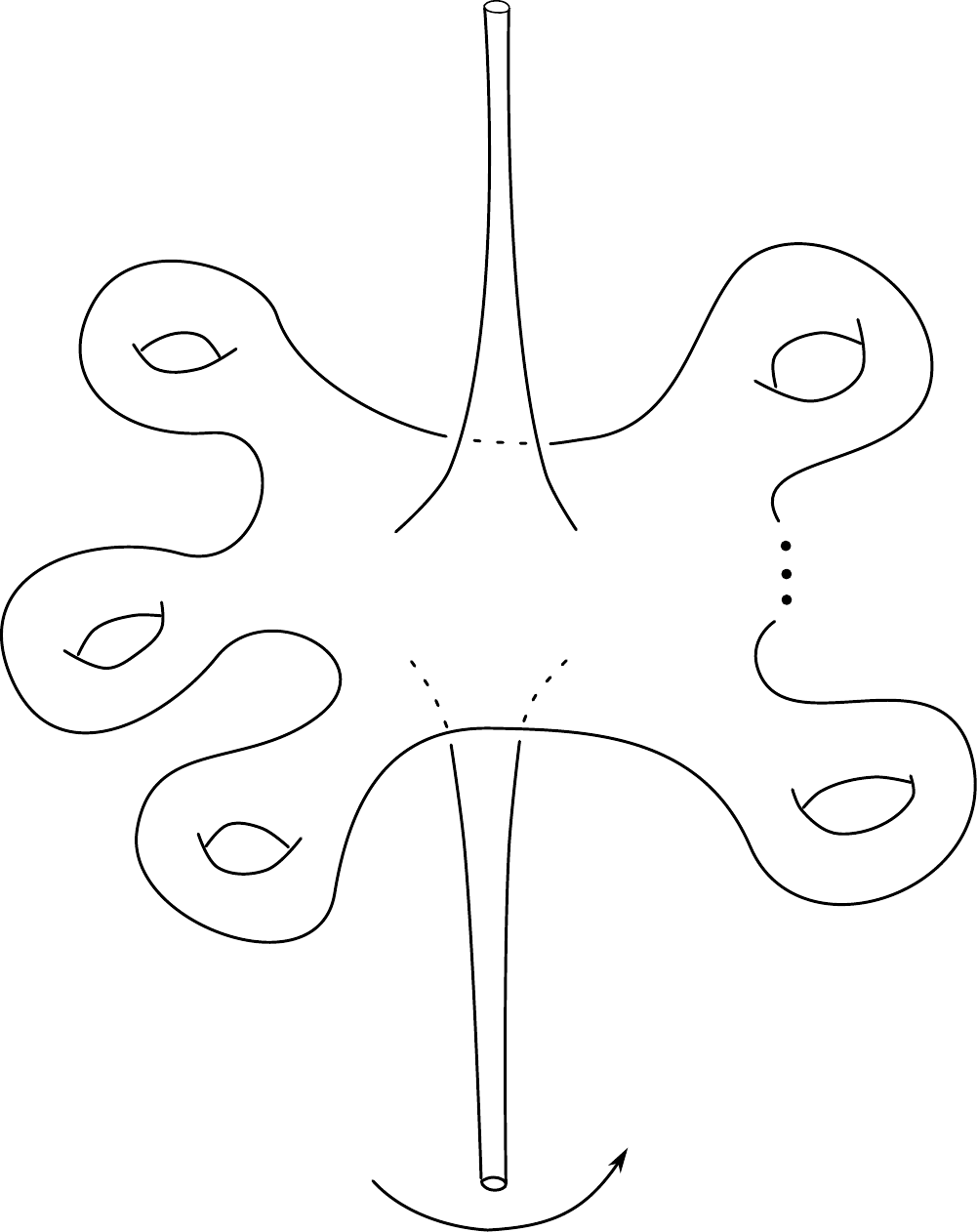}
\caption{The rotation on $\Sigma_{g,2}$ that gives a cover of $\Sigma_{1,2}$.
}\label{co}
\end{figure}

Then, as in part $(a)$, we obtain that $\hat\alpha_n,\hat\beta_n$ are stratified and 
$$vol(N_{(\hat\alpha_n,\hat\beta_n)})\asymp d_{\mathcal C}(\alpha_n,\beta_n)\asymp \log (\ell_{X}(\alpha_n)+\ell_{X}(\beta_n)),$$
where $X$ is the pull-back metric under the $g$-fold cover.

For the general case, we just start with $\Sigma_{1,2+m}$ with the sequence from part (a) and take the slit cover corresponding to $\delta$ as before. Then, the $g$-fold cover has $2+g\,m$ punctures and the same asymptotics.

\end{enumerate}

Finally, let $W = \Sigma_{g_W,k_W}$ be an essential incompressible subsurface of $S$, then
 $$\sharp\{\text{hom. classes of }  \alpha_n \text{-arcs in } W\}\leq \sharp\{\text{hom. classes of disjoint proper simple arcs in } W\}.$$
A simple Euler characteristic argument bounds this by $6g_W - 6 + 3 k_W \leq 3(\kappa(W)+3)$. Since $W$ is incompressible $\kappa(W) \leq \kappa(S)$. Lastly, as any subsurface decomposition of $S$ has at most $\kappa(S)$ subsurfaces, we obtain
$$\sum_{\Sigma\in \mathcal W}\sharp\{\text{homotopy classes of } \alpha_n \text{ and } \beta_n\text{-arcs in }\Sigma\} < 6(\kappa(S)+3)\kappa(S) \quad \text{for all} \quad n \in \mathbb N.$$
\end{proof}

\begin{customthm}{D}\label{length}Let $S = \Sigma_{1,1}$ or $\Sigma_{0,4}$, $N = PT(S)$, and let $(\al_n, \beta_n)$ be a sequence of filling pairs of essential simple closed curves in minimal position on $S$. Then there is a constant $C > 0$ and for any hyperbolic metric $X$ on $S$ there is a constant $C_X > 0$ such that:
\[  \limsup_{n \to \infty} \frac{vol(N_{(\hat{\al_n}, \hat{\beta_n})})}{\log(\iota(\al_n, \beta_n))} \leq C \quad \text{and} \quad \limsup_{n \to \infty} \frac{vol(N_{(\hat{\al_n}, \hat{\beta_n})})}{\log(\ell_{X}(\al_n) + \ell_{X}(\beta_n))} \leq C_X\]
Further, there are sequences where both equalities are attained.
\end{customthm}

\bpf By Theorem \ref{punctorus}, we know that $vol(N_{(\hat{\al}, \hat{\beta})}) \leq 2\, K_1 d_{\mathcal{C}(S)}(\al, \beta) + 2\, K_0$. A classical result of Hempel \cite{Hempel2001} and Lickorish \cite{Lickorish1962} states that $d_{\mathcal{C}(S)}(\al, \beta) \leq 2 \log_2(\iota(\al, \beta)) + 2$. Further, by a result of Basmajian \cite{Basmajian2013}, we have the bound $\iota(\al, \beta) \leq K_X (\ell_X(\al) + \ell_X(\beta))^2$. Combining the two bounds gives the two asymptotics results. Lastly, note that by Theorem \ref{ce}, that there are sequences that do attain some constants $C > 0$ and $C_X > 0$.
\epf

 \subsubsection{Examples from foliations}
 
Here, we construct examples of stratified canonical lifts coming from foliations.  Let $\mathcal F$ be a singular foliation on $\Sigma_{g,0}$ with $k$ singularities having $p_1, \ldots, p_k$ prongs, respectively. Let $\eps = +1$ if $\mathcal F$ is orientable and $-1$ otherwise.  We will allow for $p_i = 2$, which corresponds to a removable singularity. Away from the singularities, $\mathcal F$ gives rise to a section $\Sigma_{g,k} \to PT(\Sigma_{g,k})$ by looking at the tangent line field to the leaves of $\mathcal F$. For us, these sections will be stratifying surfaces for canonical lifts of non-singular closed leaves of $\mathcal F$, if they exist.

Fix $0 \leq m \leq k$. To build examples, we will work with foliations arising from quadratic differentials for a complex structure $X$ on $\Sigma_{g,m}$ and foliations associated to a pseudo-Anosov diffeomorphism. Let $Q(X)$ denote the space of holomorphic quadratic differentials with simple poles at the $m$ punctures. For $q \in Q(X)$, there are two singular transverse foliations $\mathcal F_{h,q} = \ker \mathrm{Im}\sqrt{q}$ and $\mathcal F_{v,q} =  \ker \mathrm{Re}\sqrt{q}$. For a pseudo-Anosov diffeomorphisms $\phi$ of $\Sigma_{g,m}$, the Nielsen-Thurston classification theorem produces two transverse foliations $\mathcal F_{s,\phi}$ and $\mathcal F_{u, \phi}$ that are preserved under the action of $\phi$. These are called the stable and unstable foliations, respectively. Note, in both cases, we think of the foliations as living on $\Sigma_{g,0}$, where the $m$ punctures count as singularities.

Since we will be puncturing the singularities of these foliations, we need a realizability criteria on the number of singularities.

\bthm[{\cite{MS93}}] There is a pseudo-Anosov homeomorphism $\phi$ of a $\Sigma_{g, m}$ such that $\mathcal F_{s,\phi}$ and $\mathcal F_{u,\phi}$ realise the data $(p_1, \ldots, p_k; \eps)$ with $p_i \in \mathbb N$ and $\eps = \pm 1$, if and only if:
\begin{enumerate}
\item $\sum_{i = 1}^k (p_i - 2) = 4(g-1)$
\item $\eps = -1$ is any $p_i$ is odd
\item $(p_1, \ldots, p_k; \eps) \neq (6;-1), (3,5; -1), (1, 3; -1)$ or $(\; ; -1)$
\item $\sharp \{i \mid p_i = -1\} \leq m$
\end{enumerate}
Further, we can assume that $\mathcal F_{s,\phi} = \mathcal F_{h,q}$ for some $q \in Q(X)$.
\ethm

Note, part (4) above references the convention that the $m$ punctures of $\Sigma_{g,m}$ are realized by 1-prong singularities. Notice that for any $\Sigma_{g,k}$, we can pick $0 \leq m \leq k$ and find a pseudo-Anosov or a quadratic differential for $\Sigma_{g,m}$ whose foliations have exactly $k$ singularities.
 
\begin{customthm}{H}\label{qdiff} Let $X$ be a complex structure on $\Sigma_{g,m}$. Then every $q \in Q(X)$ with $k$ singularities gives rise to a collection $cyl(q)$ of essential simple closed curves on $\Sigma_{g,k}$ such that $\{ \Gamma \subset cyl(q) \mid \Gamma \text{ is finite and } \hat{\Gamma} \text{ is stratified and  hyperbolic in } PT(\Sigma_{g,k})\}$ is infinite.
\end{customthm}

\bpf Recall that $q$ defines a finite-area, singular flat metric on $\Sigma_{g,0}$ with $k$ cone points and a natural preferred direction corresponding to $\mathcal F_{h, q}$. Rotating this direction gives a family of foliations $\mathcal F_{\theta, q} = \ker \mathrm{Im}(e^{i \theta} \sqrt{q})$ for $0 \leq \theta \leq 2 \pi$, which all share the singular points.

Let $\Theta(q) = \{ \theta \in [0, \pi) \mid \mathcal F_{\theta, q} \text{ has a non-singular closed leaf}\}$. Masur \cite{M86} showed that $\Theta(q)$ is dense in $[0, \pi)$. In particular, this gives a one to many map from $f_q : \Theta(q) \to \mathcal S$, the set of simple closed curves in $\Sigma_{g,k}$, by mapping each $\theta \in \Theta(q)$ to the collection of homotopy classes of non-singular closed leaves of $\mathcal F_{\theta, q}$ as realised on $\Sigma_{g,k}$ by puncturing all the singularities. Notice that every $\gamma \in f_q(\theta)$ has $\hat\gamma$ living on the tangent field to $\mathcal F_{\theta,q}$, which is a proper incompressible surface in $PT(\Sigma_{g,k})$. These surfaces are disjoint for distinct $\theta$ and will serve as our stratifying surfaces.

Let $cyl(q) = f_q(\Theta(q))$. If we can show that $cyl(q)$ is infinite and filling, then there are infinitely many finite filling sub-collections $\Gamma \subset cyl(q)$. Since each element of $cyl(q)$ is geodesic, each $\Gamma$ is in minimal position. By Lemma \ref{hypfill}, it follows that $\hat\Gamma$ is hyperbolic in $PT(\Sigma_{g,k})$ and we have already shown that $\hat\Gamma$ is stratified.

 \textbf{Claim 1:} $cyl(q)$ is infinite.

\bpfc Recall that the  singular flat metric given by $q$ on $\Sigma_{g,0}$ is non-positively curved, so each simple closed curve in $\Sigma_{g,0}$ has a geodesic representative. These geodesic representatives are either unions of saddle connections or closed non-singular geodesic loops. Recall that a \emph{saddle connection} is a geodesic arc in the flat metric connecting two cone points with no cone points in its interior.  Even though geodesic representatives of simple closed curves are not unique, any two that represent the same closed curve must jointly bound a flat cylinder with no cone points in its interior. In particular, they must be parallel. Since every element of $cyl(q)$ arrises as a closed non-singular geodesic loop and $\Theta(q)$ is infinite, then so is $cyl(q)$.
\epfc

{\bf Claim 2:} $cyl(q)$ is filling. 

\bpfc It is enough to show that $cyl(q)$ is filling on $\Sigma_{g,0}$ \rev{(as opposed to $\Sigma_{g,m}$)} and that it separates any two cone points. Since every simple closed curve on $\Sigma_{g,0}$ has a representative that is a union of saddle connections and any two cone points are connected by a saddle connection, it is enough to show that every saddle connection is cut by an element of $cyl(q)$. Let $L$ be a saddle connection in the flat metric defined by $q$. 
\rev{Fix} $\theta$ that is \rev{ not the direction of} $L$.  The closure of the leaves of $\mathcal F_{\theta, q}$ that meet $L$ is a subsurface $Y_\theta$ of $\Sigma_{g,0}$. If $Y_\theta$ has boundary, then some component of $\partial Y_\theta$ bounds a cylinder, which must intersect $L$ and we are done. Otherwise, there is an interval $I \subset [0, \pi)$ such that $Y_\theta$ has no boundary, i.e. $Y_\theta = \Sigma_{n,k}$ for $\theta \in I$. Further, all leaves of $\mathcal F_{\theta,q}$ must intersect $L$ for $\theta \in I$. By \cite{M86}, $\mathcal F_{\theta,q}$ has a non-singular closed leaf for some $\theta \in I$, which must intersect $L$. \epfc

Claims 1 and 2 complete the proof.\epf

\begin{figure}
\begin{overpic}[scale=.8]{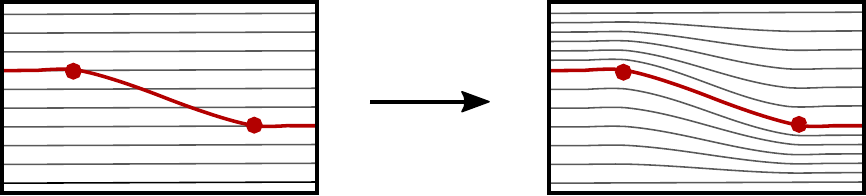}
\end{overpic}
\caption{Surgery to make the foliation parallel to the $\delta_i$ arcs.}
\label{fig:fol_surger}
\end{figure}

\brem[pseudo-Anosov Examples]\label{psuedo-anosovex} Given a pseudo-Anosov element $\phi\in \mathrm{Mod}(S_{g,k})$ we have two transverse measured foliations $\mathcal F_{s, \phi}$ and $\mathcal F_{u,\phi}$ defining a singular flat metric on $S_{g,m}$ that agrees with the transverse measures and the foliations are orthogonal with $k$-singularities. By using \rev{this, it is} possible to construct infinitely many pairs of essential simple closed loops $\gamma_1$ and $\gamma_2$ such that $\gamma_1$ is ``almost'' carried by $\mathcal F_{s,\phi}$ and $\gamma_2$ by $\mathcal F_{u,\phi}$, \rev{see the left hand-side of Figure \ref{fig:fol_surger} and the construction below}. Then, one can show that for large enough $n\in\N$ the pair: $\Gamma_n=(\phi^n(\gamma_1),\gamma_2)$ is filling and $\hat\Gamma_n$ is stratified in $PT(S_{g,k})$. By using our volume bounds, Theorem \ref{pair}, and \cite[Prop. 4.2]{MM99} one gets that there are constants $A > 1$ and $B > 0$, depending on $S_{g,k}$, such that
\[vol (N_{\hat{\Gamma}_n})\geq \frac{1}{A} d_{\mathcal C(\Sigma_{g,k})}(\phi^{n}(\gamma_1),\gamma_2)-B \]

To build $\gamma_1$, we start with a point $x_1$ on a non-singular leaf $\ell_1$ of $\mathcal F_{s,\phi}$. Pick $\eps_1 > 0$ such that the $4 \eps_1 \times \eps_1$ rectangular box $B_1$ centered at $x_1$ does not contain any singularities. Since $\mathcal F_{s,\phi}$ is uniquely ergodic, flowing along $\ell_1$ in some direction from $x_1$ is guaranteed to return to $B_1$. On our first return to $B_1$, we close up the loop by a flat geodesic segment $\delta_1$ to obtain a simple closed curve $\gamma_1$. Note that $\delta_1$ has slope between $-1/2$ and $1/2$ by the choice of $B_1$. Since the measure of $\gamma_1$ is non-zero with respect to $\mathcal F_{u,\phi}$, it must be essential. 

We construct $\gamma_2$ in the same way by starting with a point $x_2$ on a non-singular leaf $\ell_2$ of $\mathcal F_{u,\phi}$ using some $\eps_2 > 0$ and a $\eps_2 \times 4 \eps_2$ box to close up along a flat geodesic segment $\delta_2$. This time, the slope of $\delta_2$ is smaller than $-2$ or bigger than $2$. Notice that the choices of $x_i$ and $\eps_i$ will produce infinitely many distinct homotopy classes of pairs $(\gamma_1, \gamma_2)$. This can be seen by shrinking $\eps_i$ to force a longer $\gamma_i$ with respect to the transverse measures.

Smooth out the corners for each $\delta_i$ and let $\Gamma_n = (\phi^{n}(\gamma_1), \gamma_2)$.
By a surgery on $\mathcal F_{s,\phi}$ in the box $B_1$, we build a new foliation $\mathcal F_s'$ parallel to $\delta_1$, see Figure \ref{fig:fol_surger}. Note, this involves locally compressing and expanding leaves of $\mathcal F_{s,\phi}$, so the transverse measure is gone. The tangent field to $\mathcal F_s'$ lifts to an incompressible surface $S_1$ in $PT(\Sigma_{g,k})$, where we remove the $k$ singularities. Doing the same construction for $\gamma_2$ gives a surface $S_ 2 \subset PT(\Sigma_{g,k})$. Notice that because of the slopes of $\delta_i$ are in distinct slope intervals, the surfaces are disjoint.

Because $\gamma_1$ and $B_1$ avoid all singularities and $\phi(\mathcal F_{s,\phi}) = \frac{1}{\lambda} \mathcal F_{s,\phi}$, we see that $\phi^n(\mathcal F_s')$ is still transverse to $\mathcal F_{u,\phi}$. Lifting the image foliation gives a stratification of $\hat\Gamma_n$ for all $n$. Further, any immersed bigon between $\phi^{n}(\gamma_1)$ and $\gamma_2$ would have to contain punctures, so $\Gamma_n$ is in minimal position. Lastly, by \cite[Prop. 4.2]{MM99}, there is a constant $c>0$, depending only on $\Sigma_{g,k}$, such that $d_{\mathcal C(\Sigma_{g,k})}(\phi^n(\gamma_1), \gamma_1) \geq c \, n$ for all $n \in \mathbb N$. Therefore, if we choose $K \geq (3 + d_{\mathcal C(\Sigma_{g,k})}(\gamma_1, \gamma_2))/c$, then $d_{\mathcal C(\Sigma_{g,k})}(\phi^{n}(\gamma_1),\gamma_2) \geq 3$ for all $n \geq K$. Any two curves in $\mathcal C(\Sigma_{g,k})$ of distance 3 or more must be filling, so $\Gamma_n$ is filling for $n \geq K$.  Now that we have hyperbolicity, Theorem \ref{pair}, and \cite[Prop. 4.2]{MM99} give the lower bound.

\erem

 \subsection{Unstratifiable pairs of simple closed curves in higher complexity surfaces}

Recall that by Proposition \ref{canonical}, adding a null homologous simple closed curve to any collection of essential closed curves $\Gamma$ will always make $\hat\Gamma$ unstratifiable and our volume bounds cannot be directly applied. Here we will look at consequences of Lemma \ref{pucture-strat} and construct unstratifiable filling pairs of simple closed curves.

\begin{customprop}{I}\label{nogen}
Let $S$ be a hyperbolic punctured surface different from  $\Sigma_{1,1}$ or $\Sigma_{0,4}$. Then there exists a pair $(\alpha, \beta)$ of essential simple closed curves where $(\hat\alpha, \hat \beta)$ is unstratifiable.
\end{customprop}
\begin{proof}
Let $S=\Sigma_{g,k}$. By Lemma \ref{pucture-strat}, it is enough to show that there is a filling pair of simple closed curves $(\alpha,\beta)$ on $S$ such that some simply connected component of $S\setminus  ( \alpha\cup \beta)$ is not a rectangle.

\begin{itemize}
\item[$(a)$]  \textbf{Case} $g>1$. Consider a sequence of filling pairs of essential simple closed curves $(\alpha_n,\beta_n)$ on $\Sigma_{g,0}$ such that $\Sigma_{g,0}\setminus (\alpha_n\cup\beta_n)$ has at least $n+1$ connected components. Since the Euler characteristic of $\Sigma_{g,0}$ is negative then there exist at least one connected component $D$ in $\Sigma_{g,0}\setminus  ( \alpha_n\cup \beta_n)$ with $2m+6$ edges for some $m \in \mathbb N$. Taking $n \geq k$, add one puncture to $k$ distinct components of $\Sigma_{g,0}\setminus (\alpha_n\cup \beta_n)$, excluding $D$.
 \begin{figure}[h]
\begin{overpic}[scale=0.5] {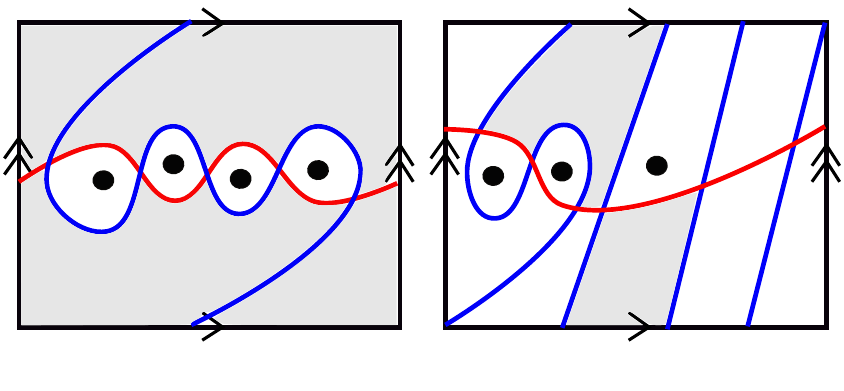}
\put(21,0){$(i)$}
\put(73,0){$(ii)$}
\end{overpic}
\caption{A pair of simple closed filling curves in $\Sigma_{1,k}$ where $k$ is even in $(i)$ and $k$ is odd ($\geq 3$) in $(ii)$.
}\label{g1}
\end{figure}
\item[$(b)$] \textbf{Case} $g=1$ and $k>1$.  Consider the filling pairs in Figure \ref{g1}. Notice that for each $k$, we have a filling pair of essential simple closed curves such that if $k$ is even, then there is a simply connected component with $4+2k$ edges and, if $k$ is odd, then there is a simply connected component with $3+k$ edges.

 \begin{figure}[h]
\begin{overpic}[scale=0.5]{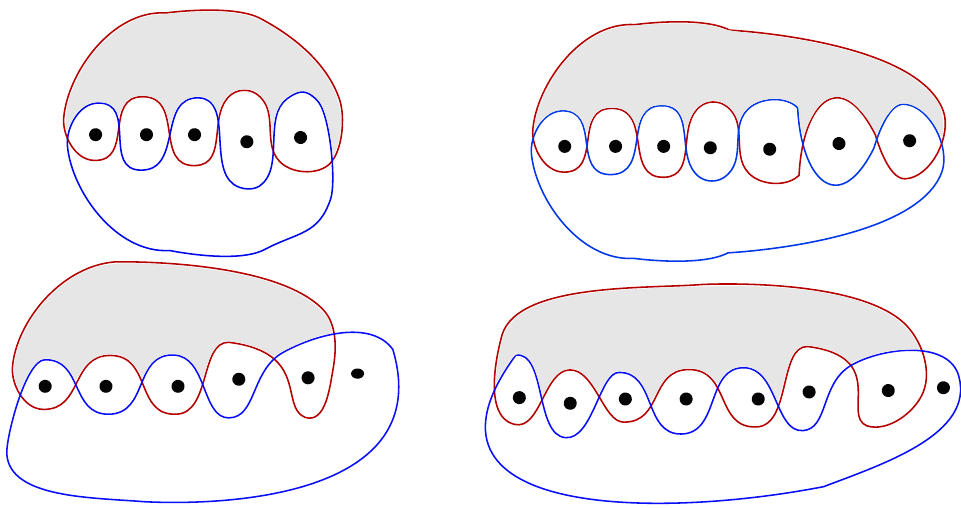}
\put(-10,37){$(i)$}
\put(-10,10){$(ii)$}
\end{overpic}
\caption{A pair of simple closed filling curves in  $(i)$ $\Sigma_{0,6}$ (left) and $\Sigma_{0,8}$ (right), and in $(ii)$ $\Sigma_{0,7}$ (left) and $\Sigma_{0,9}$ (right); where the last puncture is at infinity.}\label{n0}
\end{figure}
\item[$(c)$] \textbf{Case} $g=0$ and $k\geq 6$. Consider the filling pairs in Figure \ref{n0}. Notice that for each $k$, we get a filling pair of essential simple closed curves such that if $k$ is even, then there is a simply connected component with $k$ edges and, if $k$ is odd, then there is a simply connected component with $k-1$ edges.
\end{itemize}

\begin{figure}[h]
\begin{overpic}[scale=0.4] {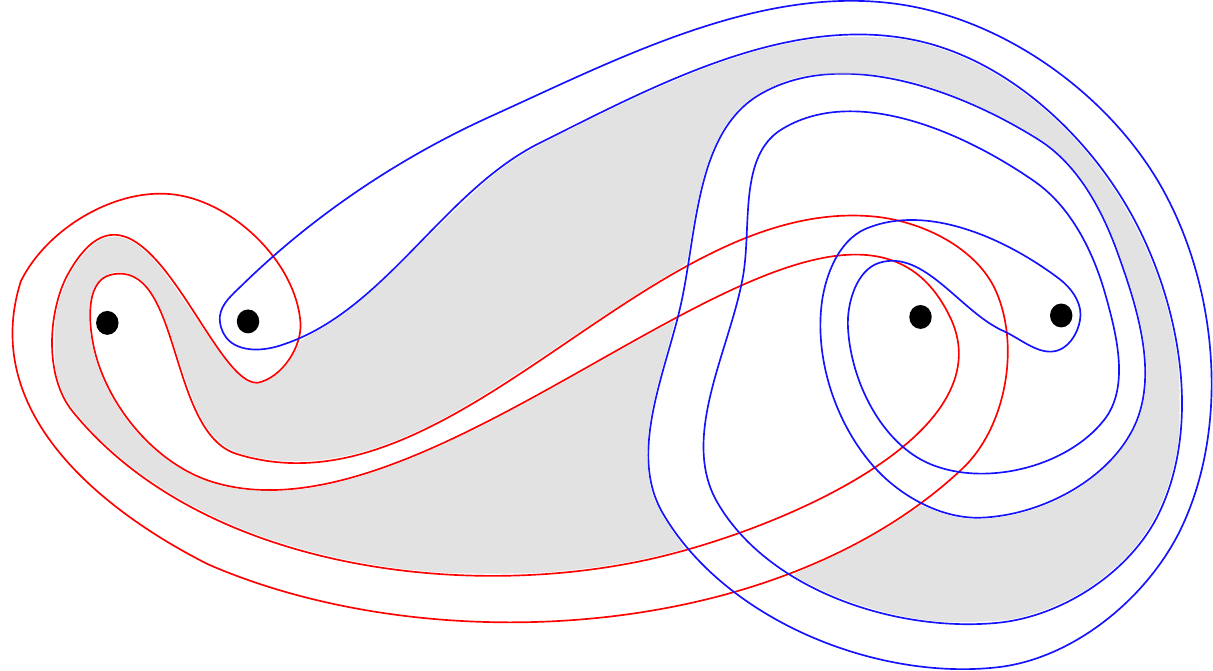}
\end{overpic}
\caption{A pair of simple closed filling curves in $\Sigma_{0,5}$ where the last puncture is at infinity.}\label{5puncsphere}
\end{figure}
Finally, for the case where $g=0$ and $k=5,$ consider the filling pair of essential simple closed curves in Figure \ref{5puncsphere}. Notice that the shaded region is simply connected with $6$ edges.\end{proof}

\nocite{BP1992,Ha2002,He1976,Sh1975,Th1978,Ja1980,MT1998,MT1998}

\thispagestyle{empty}
{\small
\markboth{References}{References}
\bibliographystyle{alpha}
\bibliography{mybib}{}
}

	\bigskip
	\small
	{\setstretch{0.8}\noindent Department of Mathematics, University of Southern California, Kaprelian Hall.

\noindent 3620 S. Vermont Ave, Los Angeles, CA 90089-2532.

email: \texttt{cremasch@usc.edu}

\vspace{0.1cm}

\noindent  Mathematisches Institut der Universit\"at M\"unchen

 \noindent Theresienstr. 39 D-80333 M\"unchen
 
 email: \texttt{migueles@math.lmu.de}

\vspace{0.1cm}
\noindent Department of Mathematics, Princeton University, Fine Hall.

\noindent Washington Road, Princeton NJ 08544-1000.

 \noindent email: \texttt{yarmola@princeton.edu}}

\end{document}